%% file: paper_DM_2015.tex
\documentclass[12pt]{article}
\usepackage[english]{babel}
\usepackage[utf8]{inputenc}

\pdfoutput=1  % pour que ça marche bien avec hal et arxiv

\usepackage{amsmath}
\usepackage{amssymb}

\usepackage{numprint}
\usepackage{graphicx}
\usepackage{subfigure}
\usepackage{color}
\newcommand\Kn{\text{Kn}}

\newcommand\fTD{F}               % Distribution function
\newcommand\xTD{\vec{x}}         % Spacial coordinates
\newcommand\uTD{\vec{u}}         % Macroscopique velocity
\newcommand\vTD{\vec{v}}         % Microscopique velocity
\newcommand\nTD{\vec{n}}         % normal vector
\newcommand\qTD{\vec{q}}         % Heat flux
\newcommand\wTD{\vec{w}}         % vector w
\newcommand\uwTD{\uTD_w}         % Boundary velocity
\newcommand\MaxwTD{\mathcal{M}}  % Maxwellian
\newcommand\pTD{\overline{\overline{\Sigma}}}              % stress tenssor
\newcommand\dxTD{\,\text{d}V}                              % volume derivative
\newcommand\dvTD{\,\text{d}v_x\text{d}v_y\text{d}v_z}      % microscopique velocity derivative
\newcommand\feqTD{{\cal E}}                              % equilibrium distribution function

\newcommand\Vin{\mathcal{V_{\text{in}}}}                   % incoming velocities
\newcommand\Vout{\mathcal{V_{\text{out}}}}                 % outgoing velocities
\newcommand\rhow{\phi}                                     % boundary density
\newcommand\Tw{T_w}                                        % boundary temperature
\newcommand\ThetaSol{\theta}
\newcommand\dotThetaSol{\dot{\theta}}
\newcommand\mass{m}
\newcommand\inertia{J}
\newcommand\Torque{T}
\newcommand\nwTD{\nTD_w}                                   % normal vector

\newcommand\Bound{\Gamma}

\newcommand\fDD{f}               % Distribution function
\newcommand\gDD{g}
\newcommand\xDD{\mathbf{x}}      % Spacial coordinates
\newcommand\uDD{\mathbf{u}}      % Macroscopique velocity
\newcommand\vDD{\mathbf{v}}      % Microscopique velocity
\newcommand\nDD{\mathbf{n}}      % normal vector
\newcommand\qDD{\mathbf{q}}      % Heat flux
\newcommand\wDD{\mathbf{w}}      % vector w
\newcommand\MaxwDD{M}  % Maxwellian
\newcommand\pDD{\Sigma}          % stress tenssor
\newcommand\dxDD{\,\text{d}S}    % surface derivative
\newcommand\dvDD{\,\text{d}v_x   % microscopique velocity derivative
				\text{d}v_y}
\newcommand\feqDD{\textsc{f}}    % equilibrium distribution function
\newcommand\geqDD{\textsc{g}}

\newcommand\uwDD{\uDD_w}         % Boundary velocity
\newcommand\nwDD{\nDD_w}         % normal vector
\newcommand\xwDD{\mathbf{r}^n_{i,j}}
\newcommand\pwDD{\Sigma_w}
\newcommand\xSolDD{\mathbf{c}}
\newcommand\uSolDD{\dot{\xSolDD}}
\newcommand\FSolDD{\mathbf{F}}

\newcommand\fBoundDD{f_w(t^n,\xwDD,\vDD_{\iv})}

\newcommand\dxUD{\,\text{d}l}    % line derivative

\newcommand\iv{p}                % velocity index

\newcommand{\virt}{\overline{\Omega}}
\newcommand\VirtDD{\virt_{i,j}}         % virtual cell
\newcommand\VirtDDim{\virt_{i-1,j}}         
\newcommand\VirtDDjp{\virt_{i,j+1}}   
\newcommand\VirtDDp{\virt_{i',j'}}
\newcommand\SVirtDD{\overline{s}_{i,j}}                  % virtual cell area
\newcommand\SVirtDDp{\overline{s}_{i',j'}}

\newcommand\CVolDD{\mathcal{C}_{i,j}}                    % control volume
\newcommand\CVolDDjpun{\mathcal{C}_{i,j+1}}                    % control volume

\newcommand\SCVolDD{s_{i,j}}                             % control volume area
\newcommand\SCVolDDp{s_{i',j'}}

\newcommand\fVirtDD{\overline{f}_{i,j,\iv}}              % Distribution function
\newcommand\fVirtDDp{\overline{f}_{i',j',\iv}}
\newcommand\fCVolDD{f_{i,j,\iv}}
\newcommand\FluxiDD[1]{\mathcal{F}_{i{#1}\frac{1}{2},j,\iv}}
\newcommand\FluxjDD[1]{\mathcal{F}_{i,j{#1}\frac{1}{2},\iv}}
\newcommand\FluxwDD{\mathcal{F}_{i,j,\iv}}
\newcommand\FluxiDDp[1]{\mathcal{F}_{i'{#1}\frac{1}{2},j',\iv}}
\newcommand\FluxjDDp[1]{\mathcal{F}_{i',j'{#1}\frac{1}{2},\iv}}
\newcommand\FluxwDDp{\mathcal{F}_{i',j',\iv}}

\newcommand\DxDD[1]{|L_{i{#1}\frac{1}{2},j}^n|}
\newcommand\DyDD[1]{|L_{i,j{#1}\frac{1}{2}}^n|}
\newcommand\DwDD{|L_{i,j}^n|}
\newcommand\LxDD[1]{L_{i{#1}\frac{1}{2},j}}
\newcommand\LyDD[1]{L_{i,j{#1}\frac{1}{2}}}
\newcommand\LwDD{L_{i,j}}
\newcommand\LxDDp[1]{L_{i'{#1}\frac{1}{2},j'}}
\newcommand\LyDDp[1]{L_{i',j'{#1}\frac{1}{2}}}
\newcommand\LwDDp{L_{i',j'}}

\newcommand\VirtTD{\virt_{i,j,k}}         % virtual cell

\newcommand\CVolTD{\mathcal{C}_{i,j,k}}                    % control volume

\newcommand\VVirtTD{\overline{V}_{i,j,k}}                  % virtual cell area
\newcommand\VVirtTDp{\overline{V}_{i',j',k'}}
\newcommand\VCVolTD{V_{i,j,k}}                             % control volume area

\newcommand\fVirtTD{\overline{F}_{i,j,k,\iv}}              % Distribution function
\newcommand\fVirtTDp{\overline{F}_{i',j',k',\iv}}
\newcommand\fCVolTD{F_{i,j,k,\iv}}
\newcommand\FluxiTD[1]{\mathcal{F}_{i{#1}\frac{1}{2},j,k,\iv}}
\newcommand\FluxjTD[1]{\mathcal{F}_{i,j{#1}\frac{1}{2},k,\iv}}
\newcommand\FluxkTD[1]{\mathcal{F}_{i,j,k{#1}\frac{1}{2},\iv}}
\newcommand\FluxwTD{\mathcal{F}_{i,j,k,\iv}}
\newcommand\FluxiTDp[1]{\mathcal{F}_{i'{#1}\frac{1}{2},j',k',\iv}}
\newcommand\FluxjTDp[1]{\mathcal{F}_{i',j'{#1}\frac{1}{2},k',\iv}}
\newcommand\FluxkTDp[1]{\mathcal{F}_{i',j',k'{#1}\frac{1}{2},\iv}}
\newcommand\FluxwTDp{\mathcal{F}_{i',j',k',\iv}}

\newcommand\DxTD[1]{|S_{i{#1}\frac{1}{2},j,k}^n|}
\newcommand\DyTD[1]{|S_{i,j{#1}\frac{1}{2},k}^n|}
\newcommand\DzTD[1]{|S_{i,j,k{#1}\frac{1}{2}}^n|}
\newcommand\DwTD{|S_{i,j,k}^n|}
\newcommand\SxTD[1]{S_{i{#1}\frac{1}{2},j,k}}
\newcommand\SyTD[1]{S_{i,j{#1}\frac{1}{2},k}}
\newcommand\SzTD[1]{S_{i,j,k{#1}\frac{1}{2}}}
\newcommand\SwTD{S_{i,j,k}}

\newcommand\xwTD{\vec{r}^n_{i,j,k}}
\newcommand\pwTD{\pTD_w}

\newcommand\fBoundTD{F_w(t^n,\xwTD,\vTD_{\iv})}

\usepackage{tikz}

\newcommand\Lvset{\phi}
\newcommand\LvsDD[2]{\phi_{{#1},{#2}}}
\newcommand\IDD[2]{\delta_{{#1},{#2}}}
\newcommand\xeDD[2]{\xDD_{{#1},{#2}}}
\newcommand\IwDD[4]{\delta_{{#1},{#2}}^{{#3},{#4}}}
\newcommand\xwiDD[1]{\xDD^w_{{#1},j}}
\newcommand\xwjDD[1]{\xDD^w_{i,{#1}}}

\begin{document}

% \begin{flushright}
%   {\it draft, \today}
% \end{flushright}

\begin{center}
{\bf A Cartesian Cut Cell Method for Rarefied Flow Simulations around Moving Obstacles}

\vspace{1cm}

G. Dechrist\'e$^1$, L. Mieussens$^2$

\bigskip
$^1$Univ. Bordeaux, IMB, UMR 5251, F-33400 Talence, France.\\
CNRS, IMB, UMR 5251, F-33400 Talence, France. \\
({\tt Guillaume.Dechriste@math.u-bordeaux1.fr})

\bigskip
$^2$Univ. Bordeaux, IMB, UMR 5251, F-33400 Talence, France.\\
CNRS, IMB, UMR 5251, F-33400 Talence, France. \\
Bordeaux INP, IMB, UMR 5251, F-33400 Talence, France. \\
INRIA, F-33400 Talence, France.\\
({\tt Luc.Mieussens@math.u-bordeaux1.fr})

\end{center}

\begin{abstract}
  For accurate simulations of rarefied gas flows around moving
  obstacles, we propose a cut cell method on Cartesian grids: it allows
  exact conservation and accurate treatment of boundary
  conditions. Our approach is designed to treat Cartesian cells and
  various kind of cut cells by the same algorithm, with no need to
  identify the specific shape of each cut cell. This makes the
  implementation quite simple, and allows a direct extension to 3D
  problems. Such simulations are also made possible by using an
  adaptive mesh refinement technique and a hybrid parallel
  implementation.  This is illustrated by several test cases,
  including a 3D unsteady simulation of the Crookes radiometer.
\end{abstract}

Keywords: kinetic equations, deterministic
method, immersed boundaries, cut cell method, rarefied gas dynamics
\\

%-----------------------------%
%                             %
%         Introduction        %
%                             %
%-----------------------------%
\section{Introduction}

% Gas is considered rarefied when the order of magnitude of a flow
% characteristic length is similar to the one of gas molecules mean free
% path.  Developments of rarefied gas flows computations have been
% initially dedicated to aerodynamics applications, such as re-entry
% problems, and have been mainly treated with DSMC
% simulations~\cite{Bird-1994}.  In the past decades, deterministic
% methods for simulations of Boltzmann equation and its related
% simplified kinetic models have come into being.  They are yet able to
% deals accurately with complex 3D problems: we could cite for instance
% the work of Chigullapalli et
% al.~\cite{Chigullapalli-Alexeenko-2011b,Chigullapalli-Alexeenko-2011a}
% concerning micro devices.

In gas dynamic problems, the rarefied regime appears when the mean
free path of the molecules of the gas is of the same order of
magnitude as a characteristic macroscopic length. The flow has to be
modeled by the Boltzmann equation of the kinetic theory of gases. Most
of numerical simulations for rarefied flows are made with the
stochastic DSMC method~\cite{Bird-1994}, especially for aerodynamical
flows in re-entry problems. In the past few years, several
deterministic solvers have been proposed, that are based on
discretizations of the Boltzmann equation or simplified models, like
BGK, ES-BGK, or Shakhov models~\cite{luc_rgd_2014}. They are efficient
for accurate simulations, multi-scale problems, or transitional flows,
for instance.

A recent issue is the account of solid boundary motion in rarefied
flow simulations.  This is necessary to simulate flows around moving
parts of micro-electromechanical systems
(MEMS)~\cite{Gad-2001,Karniadakis-2000}, as well as flows inside
vacuum pumps.  A fascinating illustration of rarefied flows with
moving boundaries is the Crookes radiometer, subject of many debates
from the late 19th to early 20th
century~\cite{Ketsdever-Gimelshein-Selden-2012}. Recent deterministic
simulations help to understand the origin of the radiometric
forces~\cite{Selden-Ngalande-Gimelshein-Ketsdever-2009,TK_2012,Taguchi-Aoki-2012b,Chen-Xu-Lee-2012}. The
numerical simulation of the Crookes radiometer is difficult because
the motion of the vanes is induced by gas/solid interaction (like
thermal creep), which means that an accurate prediction of the flow in
the vicinity of the boundary is needed in order to predict the correct
velocity of the vanes.

There are several numerical methods for moving boundary problems
designed for computational fluid dynamics: some of them have recently
been extended to deterministic discretizations of kinetic models, and
can be divided in two main categories. 

First, with body fitted methods, the mesh is adapted at each time step
so that the boundary of the computational domain always fit with the
physical boundary: moving mesh~\cite{Tang-2005} and ALE
methods~\cite{Hirt-1971,Hirt-Amsden-Cook-1974} fall into this
category.  Despite their extensive use in computational fluid
dynamics, very few similar works have been reported in kinetic theory,
except by Chen et al.~\cite{Chen-Xu-Lee-Cai-2012}.  Methods of the
second category are based on Cartesian grid computations and are
usually referred to as immersed boundary
methods~\cite{Mittal-Iaccarino-2005}.  The mesh does not change during
computations, and hence does not fit with the physical boundary.
Special treatment is applied on mesh cells that are located close to
the boundary in order to take its motion into account.  Various
extensions of these methods to kinetic theory have been proposed by
several authors
in~\cite{Arslanbekov-Kolobov-Frolova-2011,Pekardan-Chigullapalli-Sun-Alexeenko-2012,Dechriste-Mieussens-2014,
  Bernard-2015}.  Two recent variants are the inverse Lax-Wendroff
immersed boundary method proposed by Filbet and
Yang~\cite{Filbet-Chang-2013} and the Cartesian grid-based unified gas
kinetic scheme of Chen and Xu~\cite{Chen-Xu-2014}: the boundary motion
is not taken into account in these two works, but these methods could
in principle be extended to this kind of problem.
%The main idea of all these method is to interpolated a ghost value in the cell that are in the solid and that take into account the boundary condition.
We also mention the Lagrangian method: while it falls into the first
category in CFD, it does not in kinetic theory.  Indeed, whatever the
motion of the mesh, the distribution function has to be interpolated
at the foot of the characteristic for each microscopic velocity.  The
accuracy of these methods have been shown
in~\cite{Russo-Filbet-2009,Tsuji-Aoki-2013} for one dimensional
problems.  Finally, we mention that moving
boundary flows can also be treated with DSMC solvers: see, for
instance,~\cite{Ohwada-Kunihisa-2003,Rader-Gallis-Torczynshi-2011,Shrestha-Tiwari-Klar-Hardt-2014,Stefanov-Gospodinov-Cercignani-1998}.

In this paper, we try to mix the advantages of body fitted and
Cartesian methods: we present a cut cell method for computing rarefied
gas flows around moving obstacles. The cut cell method belongs to the
Cartesian grid based methods and has been widely used in computational
fluid dynamics~\cite{Ingram-Causon-Mingham-2003}.  However, this is
the first extension to moving boundary problems in kinetic theory
(complex 3D stationary DSMC simulations have already been investigated
in~\cite{LeBeau-1999,Zhang-Schwartzentruber-2012}).  This approach is
well suited to deterministic approximations of the Boltzmann equation
and is easy to implement because of the Cartesian structure of the
mesh.  Moreover, this is, up to our knowledge, the only immersed
boundary method to be conservative.  The versatility and robustness of
the technique is illustrated by various 2D flows, and by the
simulation of the unsteady rotation of the vanes of a 3D Crookes
radiometer.  This article is an extended version of our work announced
in~\cite{Dechriste-Mieussens-2014}.  Here, the Boltzmann collision
operator is replaced by BGK like models, which is approximated by a
discrete velocity method. However, this is not a restriction: other
collision operators could be used, and any velocity approximation
(like the spectral method) could be used.

Generally, the problem of cut cell methods is that it is difficult to
take into account the various shapes of cells that are cut by the solid
boundary: for instance, in 2D, a cut cell can be a triangle, a
quadrangle, or a pentagon, and this is worse in 3D. Here, we propose
a simple representation of these cells by using the notion of virtual
cells that are polygons (or polyedrals) with possibly degenerated edges (or
faces). This makes the treatment of any cut cell completely generic: in
the implementation, the different kinds of cut cells and the non cut
cells are treated by the same algorithm. This makes the extension of
the method to 3D problems very easy. However, to make large scale 3D
simulations possible, we also use an adaptive mesh refinement (AMR)
technique and a special parallel implementation.

The outline of our paper is as follows.  In section~\ref{sec:RGD}, we
give the governing equations of rarefied gas flows and introduce some
notations.  Our cut cell method is presented in
section~\ref{SecCutCellDD} for 2D problems. It is validated on three
different numerical examples in section~\ref{SubSecExamplesDD}.  Then,
in section~\ref{sec:3D}, our algorithm is extended to 3D simulations,
and a 3D unsteady simulation of the Crookes radiometer is
presented. Finally, some conclusions and perspectives are discussed in
section~\ref{sec:concl}. Technical details like computations of
geometric parameters of the cells are presented in the Appendix.

%-----------------------------%
%                             %
%   Rarefied gas dynamics     %
%                             %
%-----------------------------%
\section{Rarefied gas dynamics}\label{sec:RGD}
\subsection{Boltzmann equation}\label{SubSecBoltzmannTD}
%When a gas is far from thermodynamical equilibrium state, the basic macroscopic conservation laws are no more valid.
In rarefied regimes, a monoatomic gas is described by the Boltzmann equation:
\begin{equation}\label{BoltzmannTD}
	\frac{\partial \fTD}{\partial t} + \vTD\cdot \nabla \fTD = Q(\fTD).
\end{equation}
The distribution function $F(t,\xTD,\vTD)$ is the mass density of 
molecules at time $t$ that are located at the space coordinate
$\xTD\in\mathbb{R}^3$ and that have a velocity $\vTD\in \mathbb{R}^3$.
For our approach, it is more relevant to look at the integral form
of~(\ref{BoltzmannTD}) in a time dependent volume $V(t)$.
The Reynolds transport theorem leads to:
\begin{equation}\label{IntBoltzmannTD}
	\frac{\partial}{\partial t}\int_V \fTD\dxTD + 
		\int_{\partial V} (\vTD-\wTD)\cdot\nTD \fTD \dxDD =
		\int_V Q(\fTD) \dxTD,
\end{equation}
where $\partial V(t)$ is the surface of the volume $V(t)$.
Let $\xTD$ be a point of this surface: it is moving at a velocity
$\wTD(t,\xTD)$ and the vector $\nTD(t,\xTD)$ is the outward normal
vector to the surface at this point.

The density $\rho$, momentum $\rho\uTD$, total energy $E$, stress
tensor $\pTD$ and heat flux $\vec{q}$, are computed by the first
moments of the distribution function with respect to the velocity:
\begin{equation}\label{MacroscopicTD}\begin{aligned}
	\left[\begin{array}{c}
		\rho\\ \rho\uTD\\ E
		\end{array}\right] &=
		\int_{\mathbb{R}^3} \left[\begin{array}{l}
		1 \\ \|\vec{v}\| \\ \frac{1}{2}\|\vTD\|^2
		\end{array} \right]
		\fTD(t,\xTD,\vec{v}) \dvTD,\\
	\pTD &=\int_{\mathbb{R}^3} (\vTD-\uTD)\otimes(\vTD-\uTD)
		 \,\fTD(t,\xTD,\vTD) \dvTD,\\
	\qTD &= \int_{\mathbb{R}^3} 
		\frac{1}{2}(\vTD-\uTD) \|\vTD-\uTD\|^2
		\,\fTD(t,\xTD,\vTD) \dvTD,
\end{aligned}\end{equation}
where the norm is defined by $\|\vTD\|^2=v_x^2+v_y^2+v_z^2$.
The temperature $T$ of the gas is related to the 
the energy by the relation
$E=\frac{1}{2}\rho\|\uTD\|^2+\frac{3}{2}\rho RT$,
where $R$ is the gas constant defined as the ratio between the Boltzmann constant and the molecular mass of the gas.
Moreover, the pressure is computed with the standard equation of state for ideal gases: $P=\rho RT$.

When a gas is at rest, which means in equilibrium state, the molecules are uniformly distributed around the macroscopic velocity and
the distribution function is a Gaussian function called Maxwellian:
\begin{equation}\label{MaxwellianTD}
	\MaxwTD[\rho,\uTD,T] (\vec{v})=\frac{\rho}{(2\pi RT)^{3/2}}
	\exp\left(-\frac{\|\vTD-\uTD\|^2}{2RT}\right).
\end{equation}
The Boltzmann equation implies that the total variation of $\fTD$ comes from the collisions between particles, that are modeled by the Boltzmann operator $Q(\fTD)$.
This operator is computationally expensive and several models have been introduced to make its computation easier.
Simplest models  consists in a relaxation of the distribution function towards a corresponding equilibrium function $\feqTD$:
\begin{displaymath}
	Q(\fTD) = \frac{1}{\tau} (\feqTD-\fTD),
\end{displaymath}
where $\tau$ is a relaxation time (see various definitions in section~\ref{SubSecExamplesDD}).
Bhatnagar, Gross and Krook~\cite{Bhatnagar-Gross-Krook-1954} proposed to take $\feqTD$ equal to the local equilibrium state, that is 
\begin{displaymath}
	\feqTD=\MaxwTD[\rho,\uTD,T].
\end{displaymath}
The corresponding BGK model conserves mass, momentum and total energy since the first three moments of the Maxwellian are the same as those of the distribution function.
A Chapman-Enskog expansion~\cite{Chapman-Cowling-1939} gives relations between viscosity $\mu$, heat conduction $\kappa$ and relaxation time.
In this case, the expansion yields
$\mu=\tau\,P$ and $\kappa=\frac{5}{2}R\,\tau\,P$.
It may be seen that the BGK model necessarily leads to Prandtl number $\Pr=\frac{5}{2}\mu R/\kappa$ equal to 1.
However for most of gases, it is physically found to be less than this.
More complex functions $\feqTD$ such that in Holway~\cite{Holway-1966} and Shakhov~\cite{Shakhov-1968} models have been developed to obtain the correct Prandtl number.
For instance, for the Shakhov model, $\feqTD$ reads
\begin{displaymath}
	\feqTD=\MaxwTD[\rho,\uTD,T]\left[1+(1-\Pr)(\vTD-\uTD)\cdot\qTD
	\left(\frac{\|\vTD-\uTD\|^2}{RT}-5\right)\Big/(5PRT)\right].
\end{displaymath}

In this article, solid wall interactions are taken into account by the
standard fully diffuse reflection.  This model states that all
particles that collide with a boundary are absorbed by the wall and re
emitted with a Maxwellian distribution:
\begin{equation}\label{BoundaryTD}
		F(t,\xTD\in\Bound,\vTD\in\Vin) = \rhow\MaxwTD[1,\uwTD,\Tw],
\end{equation}
where $\Tw$ and $\uwTD$ are the temperature and the velocity of the boundary $\Bound$ at position $\xTD$.
The coefficient $\rhow$ is computed in order to set the net mass flux across the wall to zero:
\begin{displaymath}
	\rhow = -\left(\int_{\vTD\in\Vout}(\vTD-\uwTD)\cdot\nwTD\,\fTD\dvTD\right)\Big/
		\left(\int_{\vTD\in\Vin}(\vTD-\uwTD)\cdot\nwTD\,\MaxwTD[1,\uwTD,\Tw]\dvTD\right).
\end{displaymath}
The set of incoming velocities is defined by $\Vin=\{\vTD \text{ such that }(\vTD-\uwTD)\cdot\nwTD<0\}$, where $\nwTD$ is the normal vector to the boundary, pointed to the wall.
Similarly, the set of outgoing velocities is $\Vout=\{\vTD \text{ such that }(\vTD-\uwTD)\cdot\nwTD>0\}$.
Note that the boundary condition is defined only for the relative incoming microscopic velocities.
\subsection{Reduced model}\label{SubSecBoltmannDD}
For plane flows, the computational complexity of the Boltzmann equation can be decreased by the use of a standard reduced distribution technique~\cite{Chu-1965}.
This classical method has been extensively used for numerical
computations of BGK and Shakhov models.
First note that in plane flows, the third component of the macroscopic velocity $u_z$ is equal to zero, as well as $q_z$, $\Sigma_{xz}$, and $\Sigma_{yz}$.
From now on, we define the two dimensional variables
\begin{displaymath}
	 \xDD=(x,y),\quad\vDD=(v_x,v_y),\quad \uDD=(u_x, u_y),\quad \qDD=(q_x, q_y), \quad
 	 \pDD= \left(\begin{matrix} \Sigma_{xx} & \Sigma_{xy}
  					\\ \Sigma_{yx} & \Sigma_{yy} \end{matrix}\right).
\end{displaymath} 
Let $\fDD$ and $\gDD$ be the reduced distribution functions defined by
\begin{displaymath}
	 \fDD(\vDD) = \int_{\mathbb{R}} \fTD(\vTD)\text{d}v_z,\quad
	 \gDD(\vDD) = \int_{\mathbb{R}} \frac{1}{2}v_z^2\, \fTD(\vTD)\text{d}v_z.
\end{displaymath}
The macroscopic quantities can  be computed from $\fDD$ and $\gDD$.
Indeed, the set of equations~(\ref{MacroscopicTD}) readily becomes
\begin{equation}\label{MacroscopicDD}\begin{aligned}
	\left[\begin{array}{c}
		\rho\\ \rho \uDD\\ E	
		\end{array}\right] &=
		\int_{\mathbb{R}^2} \left[\begin{array}{c}
		1\\ \vDD\\ \frac{1}{2}\| \vDD\|^2
		\end{array} \right] \fDD(\vDD)\dvDD
		+ \int_{\mathbb{R}^2} 
		\left[\begin{array}{c}
		0\\ \mathbf{0}\\ 1	
		\end{array}\right]
		\gDD(\vDD)  \dvDD,\\
	\pDD &=\int_{\mathbb{R}^2} (\vDD-\uDD)\otimes (\vDD-\uDD)
		\fDD(\vDD)\dvDD, \\
	\qDD &= \int_{\mathbb{R}^2}(\vDD-\uDD)\left(\frac{1}{2}\|\vDD-\uDD\|^2 \fDD(\vDD)
		+\gDD(\vDD)\right)\dvDD,
\end{aligned}\end{equation}
where $\|\vDD\|^2=v_x^2+v_y^2$.
Multiplying by $(1,\frac{1}{2}v_z^2)$ the Boltzmann equation~(\ref{IntBoltzmannTD})
%, as well as its integral form~(\ref{IntBoltzmannTD}), 
and then integrating the result with respect to $v_z$  gives the following set of equations:
\begin{equation}\label{IntBoltzmannDD}
%	\begin{aligned}
%		\frac{\partial \fDD}{\partial t} + \vDD\cdot \nabla \fDD = 
%			\frac{1}{\tau} (\feqDD-\fDD),\\
%		\frac{\partial \gDD}{\partial t} + \vDD\cdot \nabla \gDD = 
%			\frac{1}{\tau} (\geqDD-\gDD),
%	\end{aligned}
%	\qquad\text{ and }\qquad
	\begin{aligned}
		\frac{\partial}{\partial t}\int_{S} \fDD \dxDD 
			+ \int_{\partial S} (\vDD-\wDD)\cdot \nDD\fDD\dxUD = 
			\int_S \frac{1}{\tau} (\feqDD-\fDD)\dxDD, \\
		\frac{\partial}{\partial t}\int_{S} \gDD \dxDD 
			+ \int_{\partial S} (\vDD-\wDD)\cdot \nDD\gDD\dxUD = 
			\int_S \frac{1}{\tau} (\geqDD-\gDD)\dxDD.
	\end{aligned}
\end{equation}
In this case, $\partial S(t)$ is the contour of the surface $S(t)$.
Each point $\xDD\in\partial S(t)$ is moving at a velocity $\wDD(t,\xTD)$ and $\nDD(t,\xTD)$ is the outward normal vector to the contour at this point.
The reduced equilibrium functions are defined by
\begin{equation}\label{BGKDD}
%\begin{aligned}
%	\feqDD &= \MaxwDD[\rho,\uDD,T],\\
%	\geqDD &= \frac{RT}{2}\MaxwDD[\rho,\uDD,T],
%\end{aligned}
\feqDD=\MaxwDD[\rho,\uDD,T]\quad\text{ and }\quad\geqDD=\frac{RT}{2}\MaxwDD[\rho,\uDD,T],
\end{equation}
for the BGK model and by
\begin{equation}\label{ShakhovDD}\begin{aligned}
	\feqDD &= \MaxwDD[\rho,\uDD,T]\left[1+(1-\Pr)(\vDD-\uDD)\cdot\qDD
			\left(\frac{\|\vDD-\uDD\|^2}{RT}-4\right)\Big/(5PRT)\right],	\\
	\geqDD &= \frac{RT}{2}\MaxwDD[\rho,\uDD,T]\left[1+(1-\Pr)(\vDD-\uDD)\cdot\qDD
			\left(\frac{\|\vDD-\uDD\|^2}{RT}-2\right)\Big/(5PRT)\right],
\end{aligned}\end{equation}
for the Shakhov model, where $\MaxwDD[\rho,\uDD,T]$ is the reduced Maxwellian given by
\begin{displaymath}
	 \MaxwDD[\rho,\uDD,T](\vDD)=\frac{\rho}{2\pi RT}\exp\left(-\frac{\|\vDD-\uDD\|^2}{2RT}\right).
\end{displaymath}
To close this section, the boundary conditions~(\ref{BoundaryTD}) are written with the reduced distribution functions as:
\begin{equation}\begin{aligned}\label{BoundaryDD}
	&\fDD(t,\xDD\in\Bound,\vDD\in\Vin) = \rhow\MaxwDD[1,\uwDD,\Tw], \\
	&\gDD(t,\xDD\in\Bound,\vDD\in\Vin) = \rhow\frac{RT}{2}\MaxwDD[1,\uwDD,\Tw], 
\end{aligned}\end{equation}
where $\rhow$ is computed by
\begin{displaymath}
	\rhow = -\left(\int_{\vDD\in\Vout}(\vDD-\uwDD)\cdot\nwDD\,\fDD\dvDD\right)\Big/
		\left(\int_{\vDD\in\Vin}(\vDD-\uwDD)\cdot\nwDD\,\MaxwDD[1,\uwDD,\Tw]\dvDD\right). 
\end{displaymath}
In this formula, $\uwDD$ and $\Tw$ are the velocity and temperature of
the point $\xDD$ that belongs to the boundary~$\Bound$, and $\nwDD$ is the normal vector to the boundary pointed to the wall.
The relative incoming and outgoing velocities at this point are
therefore defined by $\Vin=\{\vDD \text{ such that
}(\vDD-\uwDD)\cdot\nwDD<0\}$ and $\Vout=\{\vDD  \text{ such that }(\vDD-\uwDD)\cdot\nwDD>0\}$.
%
%-----------------------------%
%                             %
%   Numerical discretization  %
%                             %
%-----------------------------%
\section{The cut-cell method for two dimensional problems}
\label{SecCutCellDD}
In this section, we present a numerical method to simulate plane flows with moving boundaries.
The governing equations are detailed in
section~\ref{SubSecBoltmannDD}. The discretization of each variable
(velocity, space, and time) is presented in separate sections.

\subsection{Discrete velocity approximation}\label{SubSecVelocity}
The velocity space is discretized by a Cartesian grid.
Let $\vDD_{\min}\in \mathbb{R}^2$ and $\vDD_{\max}\in \mathbb{R}^2$ be the lower-left and upper-right  corners of this grid.
The number of discrete velocities is $N^2$, the velocity step is
denoted by $(\Delta v_x, \Delta v_y) = (\vDD_{\max}-\vDD_{\min})/N$,
and the $\iv^{th}$ velocity is $\vDD_{\iv}=\vDD_{\min}+(\iv_1\Delta v_x,\iv_2\Delta v_y)$, such that $\iv=\iv_2 N+\iv_1$ for all $(\iv_1,\iv_2)\in [0,N-1]^2$.
The approximation of the distribution function is defined by $\fDD_{\iv}(t,\xDD)= \fDD(t,\xDD,\vDD_{\iv})$.
The set of equations~(\ref{IntBoltzmannDD}) is discretized with
respect to $\vDD$ by the following set of $2N^2$ equations:
\begin{equation}\label{DiscreteIntBoltzmannDD} \begin{split}
	\frac{\partial}{\partial t}\int_{S} \fDD_{\iv} \dxDD 
		+ \int_{\partial S} (\vDD_{\iv}-\wDD)\cdot\nDD\,\fDD_{\iv}\dxUD = 
		\int_S \frac{1}{\tau} (\feqDD_{\iv}-\fDD_{\iv})\dxDD, \\
	\frac{\partial}{\partial t}\int_{S} \gDD_{\iv} \dxDD 
		+ \int_{\partial S} (\vDD_{\iv}-\wDD)\cdot\nDD\,\gDD_{\iv}\dxUD = 
		\int_S \frac{1}{\tau} (\geqDD_{\iv}-\gDD_{\iv})\dxDD.
\end{split}\end{equation}
The macroscopic quantities are computed with~(\ref{MacroscopicDD}), where the integrals over $\mathbb{R}^2$ are approximated by a sum over the $N^2$ discrete velocity points.
They are therefore given by
\begin{equation}\label{MacroscopicDDv}\begin{split}
	\left[\begin{array}{c}
		\rho\\ \rho \uDD\\ E
		\end{array}\right] &=
		\sum_{\iv=0}^{N^2-1} \left[\begin{array}{c}
		1\\ \vDD_{\iv}\\ \frac{1}{2}||\vDD_{\iv}||^2
		\end{array} \right] \fDD_{\iv}\,\Delta v_x\Delta v_y
		+ \sum_{\iv=0}^{N^2-1} \left[\begin{array}{c}
		0\\ \mathbf{0}\\ 1	\end{array}\right] \gDD_{\iv} 
		\,\Delta v_x\Delta v_y,\\
	\pDD &= \sum_{\iv=0}^{N^2-1} (\vDD_{\iv}-\uDD)\otimes(\vDD_{\iv}-\uDD) \fDD_{\iv}
		 \, \Delta v_x\Delta v_y,\\
	\qDD &= \sum_{\iv=0}^{N^2-1}	(\vDD_{\iv}-\uDD)
		\left(\frac{1}{2} ||\vDD_{\iv}-\uDD||^2 \fDD_{\iv} + \gDD_{\iv}\right)\, \Delta v_x\Delta v_y.
\end{split}\end{equation}
Finally, the equilibrium functions $\feqDD_{\iv}$ and $\geqDD_{\iv}$
are computed either with~(\ref{BGKDD}) or with~(\ref{ShakhovDD}).
However, the triplet $(\rho,\uDD,T)$ used in these formulas is
obtained with a Newton algorithm that preserves the discrete
conservation of Boltzmann equation, rather than with direct
computation~(\ref{MacroscopicDDv}) of the macroscopic quantities.
Note that instead of using the algorithm of~\cite{luc_jcp} which is
based on entropic variables, we use the algorithm of
Titarev~\cite{Titarev-2009}.
%In order to be conservative, the three first moments of collision term has to be equal to zero.
%Although equation~() complies with this condition at the continuous level, this is no longer valid at the discrete level.
%Indeed the three first discrete moments of the Maxwellian~() differs from those of the distribution function.
%A formulation of a discrete conservative collision term has been first proposed by~\cite{Mieussens-2000a} in the framework of BGK and ES-BGK models computation.
%An extension of this approach to Shakhov model computation is given in~\cite{}.

\subsection{Space discretization}

\subsubsection{Cartesian grid and cut cells}

Let $\Omega = [x_{\min},x_{\max}]\times[y_{\min},y_{\max}]$ denote the
space computational domain.
It is discretized by a Cartesian grid of $(N_x+1)\times(N_y+1)$ points.
Their coordinates are computed for all $(i,j)\in[0,N_x]\times[0,N_y]$ by 
$\xDD_{i+\frac{1}{2},j+\frac{1}{2}}=\xDD_{\min}+(i\Delta x,j\Delta y)$,
where $(\Delta x, \Delta y)=([x_{\max}-x_{\min}]/N_x,[y_{\max}-y_{\min}]/N_y)$ and $\xDD_{\min}=(x_{\min},y_{\min})$.
The computational mesh is therefore made up by $N_x\times N_y$
rectangular cells: each cell is denoted by $\Omega_{i,j}$ and its
center is $\xDD_{i,j}$.

\begin{figure}[h]
\begin{center}
\subfigure
	{\resizebox{0.44\textwidth}{!}{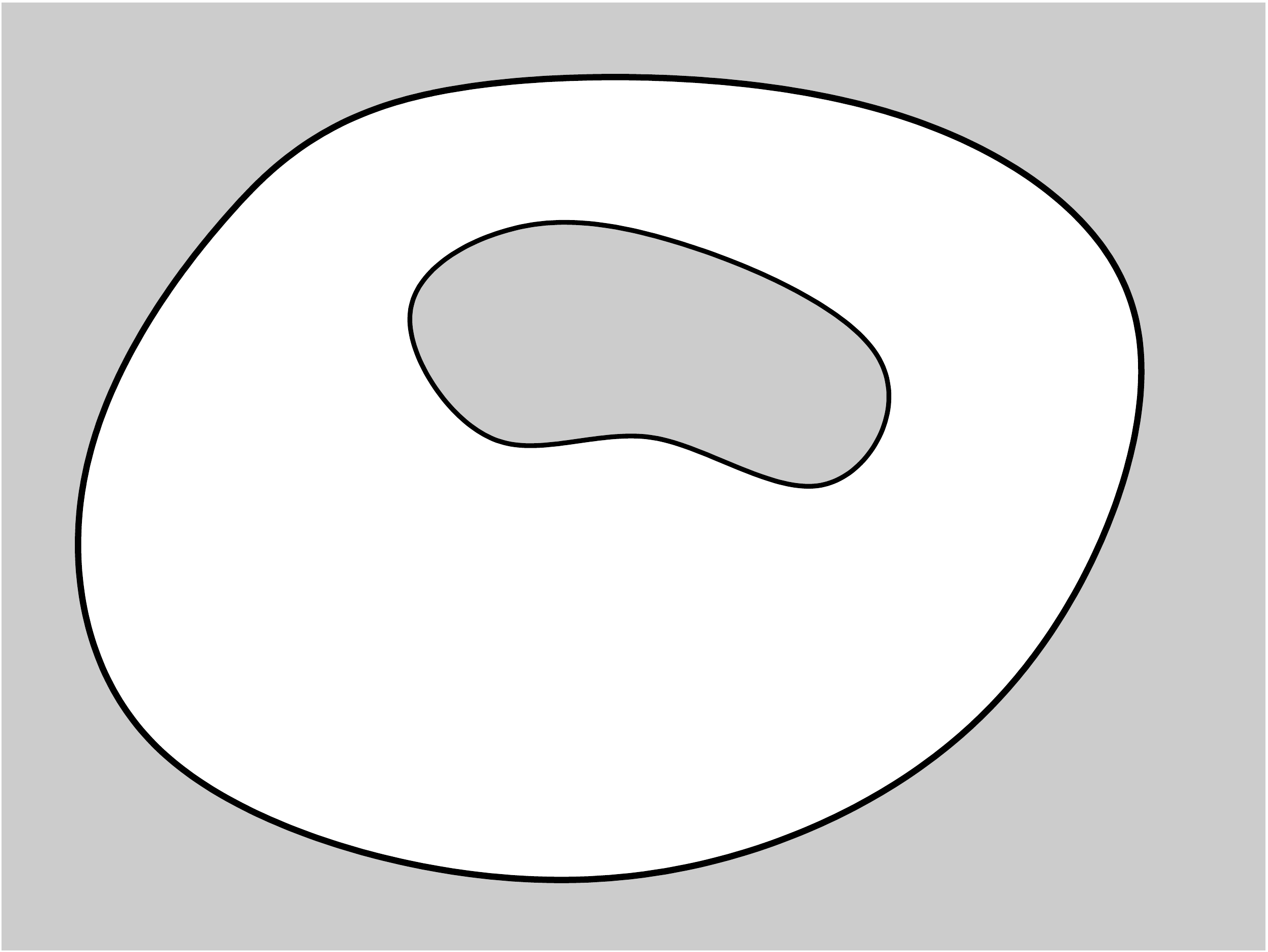}}
\hspace*{0.1\textwidth}	
\subfigure
	{\resizebox{0.44\textwidth}{!}{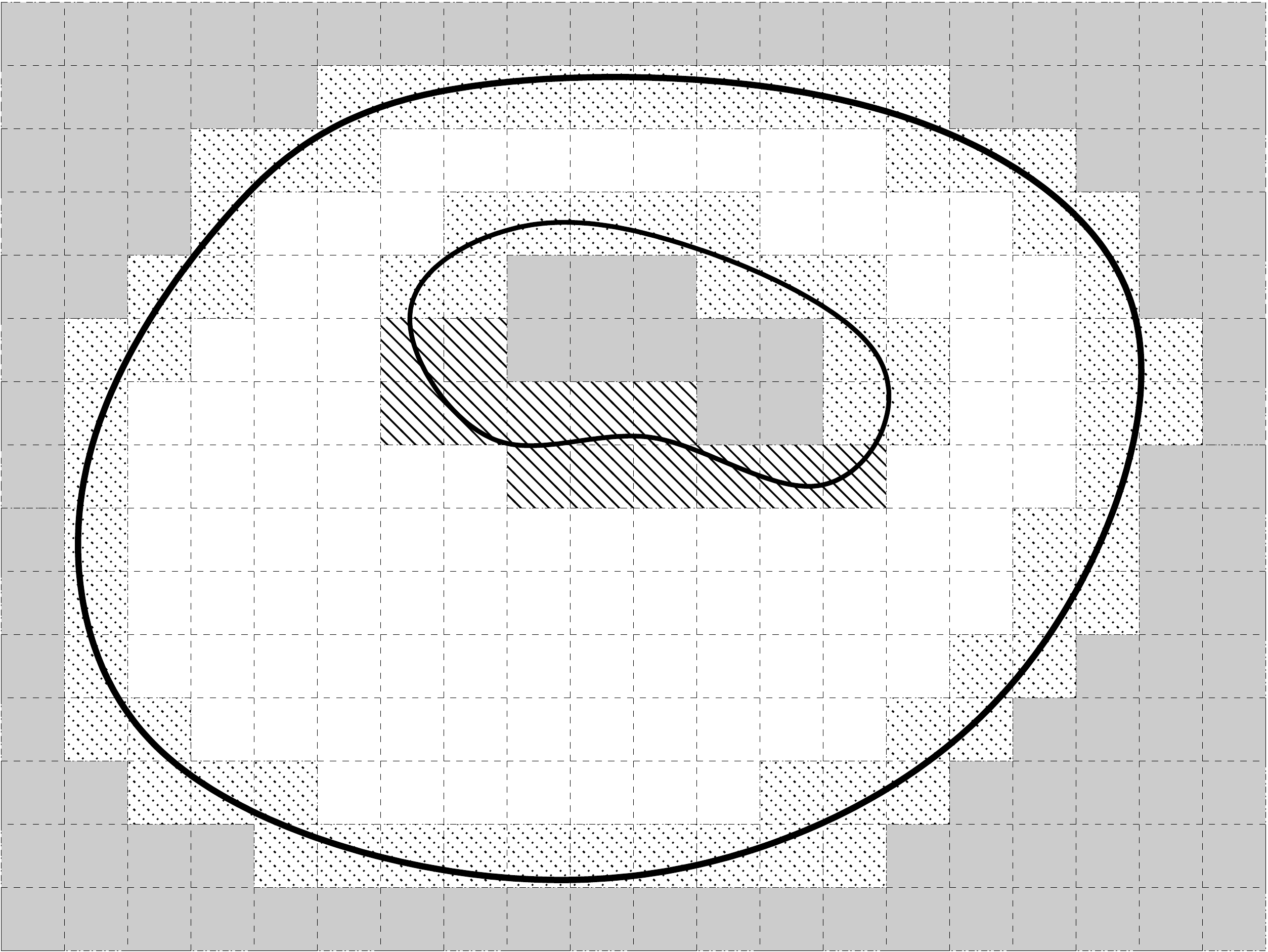}}
\end{center}
	\caption{Cells of the computational domain are classified in three categories: 
	gas cells are represented in white, solid cells are shaded and cut cells are hatched.}
	\label{Categories}
\end{figure}

Since the computational domain is rectangular, physical boundaries do not necessarily fit with the mesh boundary.
In order to simulate arbitrary shaped objects, solid and gaseous domains $\Omega_s(t)$ and $\Omega_g(t)$ are introduced.
They correspond to the solid and gaseous parts of the computational domain and hence $\Omega=\Omega_s(t)\cup\Omega_g(t)$.
We point out that while $\Omega_s$ and $\Omega_g$ are time dependent, $\Omega$ is not.
At time~$t>0$, a rectangular cell $\Omega_{i,j}$ can be
in one of these three different states only:
\begin{itemize}
	\item $\Omega_{i,j}$ is a \textit{gas cell} if it is completely contained in the gaseous domain: $\Omega_{i,j}\cap\Omega_s(t)=\emptyset$.
	\item $\Omega_{i,j}$ is a \textit{solid cell} if it is completely contained in the solid domain: $\Omega_{i,j}\cap\Omega_g(t)=\emptyset$.
	\item $\Omega_{i,j}$ is a \textit{cut cell} if it is partially contained in the gaseous domain and partially contained in the solid domain:
	$\Omega_{i,j}\cap\Omega_s(t)\neq\emptyset$ and $\Omega_{i,j}\cap\Omega_g(t)\neq\emptyset$.
\end{itemize}
These three states of cells are shown in figure~\ref{Categories}.

\subsubsection{Virtual cells}

To each cell $\Omega_{i,j}$ is now associated a virtual cell
$\VirtDD(t)$, which is the section of $\Omega_{i,j}$ 
contained in the gaseous domain: this reads
$\VirtDD(t)=\Omega_{i,j}\cap\Omega_g(t)$. Whatever the state of the
cell (that is to say gas, solid or cut), it is defined with five
virtual edges, whose lenghts can be zero. Four of them, that are
denoted by~$\LxDD{\pm}(t)$ and $\LyDD{\pm}(t)$, fit with the lines of the
Cartesian mesh. The last one, denoted by $\LwDD(t)$, is a linear
approximation of the solid boundary.  If the virtual cell has less
than five real edges, then at least one length is zero.  Finally, we
denote by $\SVirtDD$ the area of the virtual cell, and $|L|$ will denote
the length of any edge $L$.

At a given time $t^n=n\Delta t$, all these parameters are denoted as
follows: 
% Since the shape of the solid and gaseous domains change in time, the
% state of a cell is also time dependent, and the state of a virtual
% cell as well. We anticipate \fbox{ verifier ce mot }on the section devoted to the time
% discretization (section \ref{subsec:numscheme}) to detail how the
% evolution of the virtual cell is taken into account. The time step is
%denoted by $\Delta t$ and $t^n=n \Delta t$, and we define
\begin{displaymath}
	\VirtDD^n=\VirtDD(t^n),\quad \LxDD{\pm}^n=\LxDD{\pm}(t^n),\quad 
	\LyDD{\pm}^n=\LyDD{\pm}(t^n), \quad \LwDD^n=\LwDD(t^n),\quad
	\SVirtDD^n=\SVirtDD(t^n).
\end{displaymath}

\begin{figure}[h]\begin{center}
	\subfigure%[$(\DxDD{-},\DyDD{-})=\Delta\xDD$]
		{\resizebox{0.25\textwidth}{!}{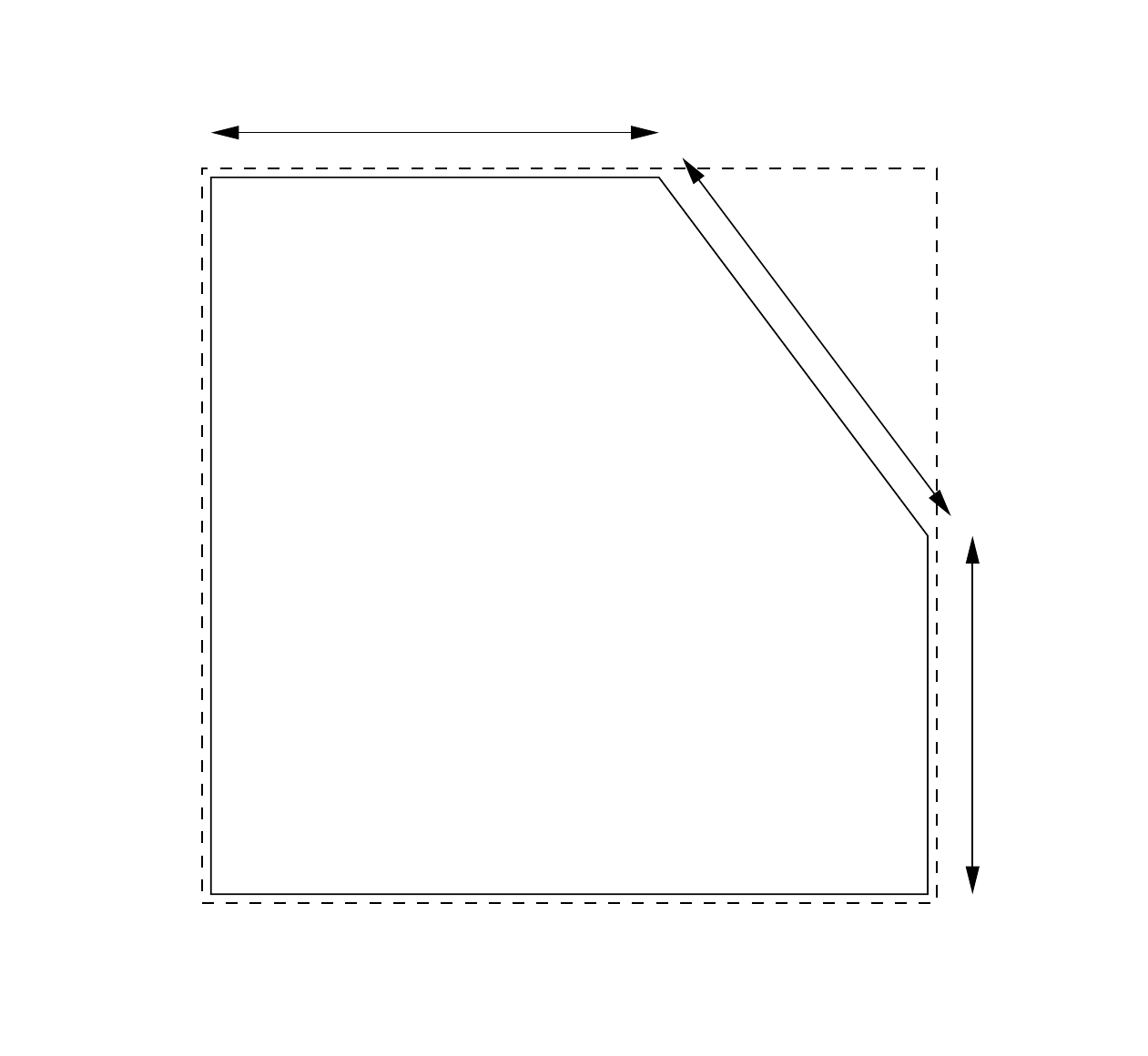}}
	\hspace*{0.1\textwidth}
	\subfigure%[$\DxDD{-}=\DyDD{-}=0$]
		{\resizebox{0.25\textwidth}{!}{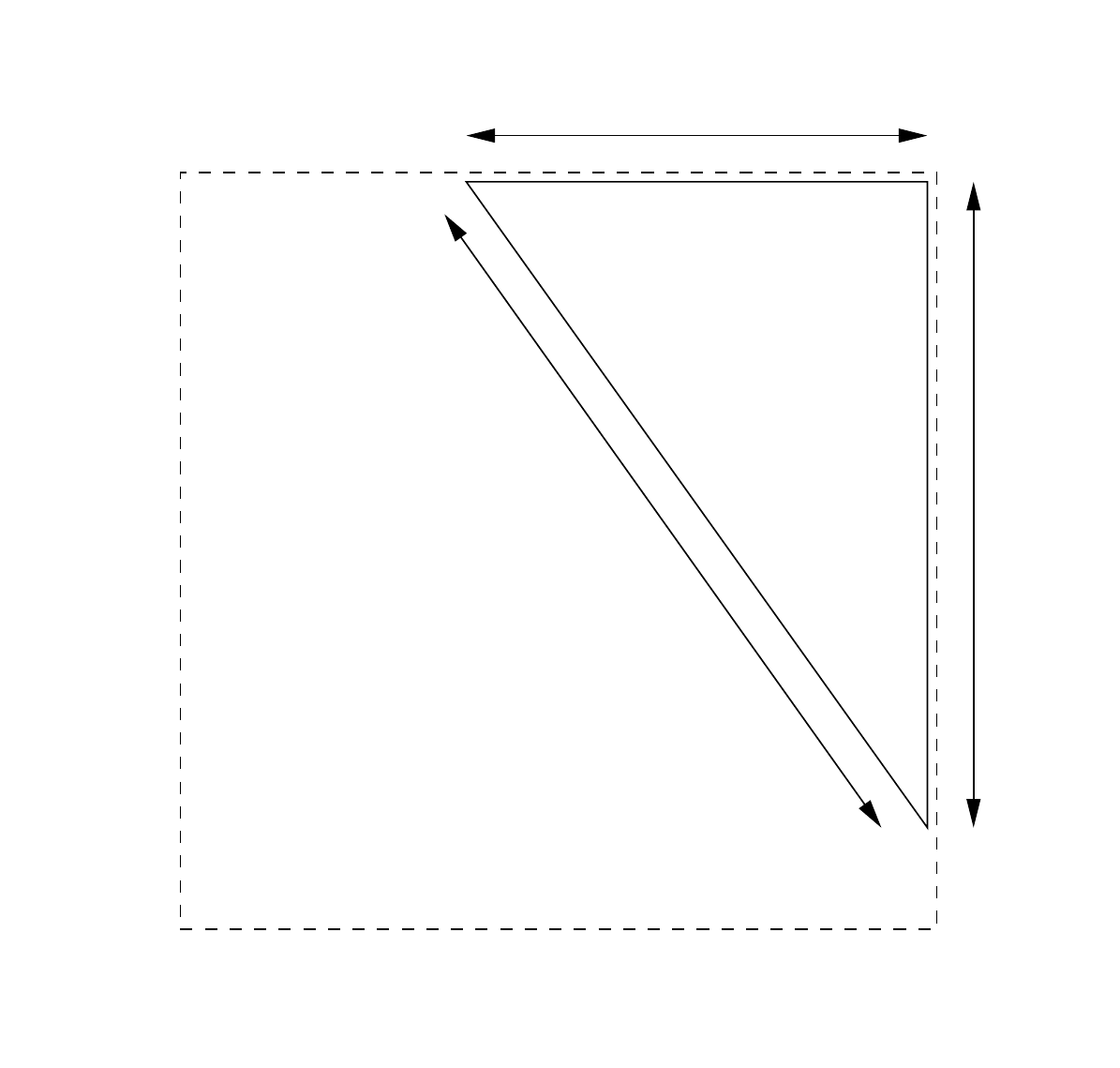}}
	\hspace*{0.1\textwidth}
	\subfigure%[$\DyDD{+}=0$, and $\DyDD{-}=\Delta y$]
		{\resizebox{0.25\textwidth}{!}{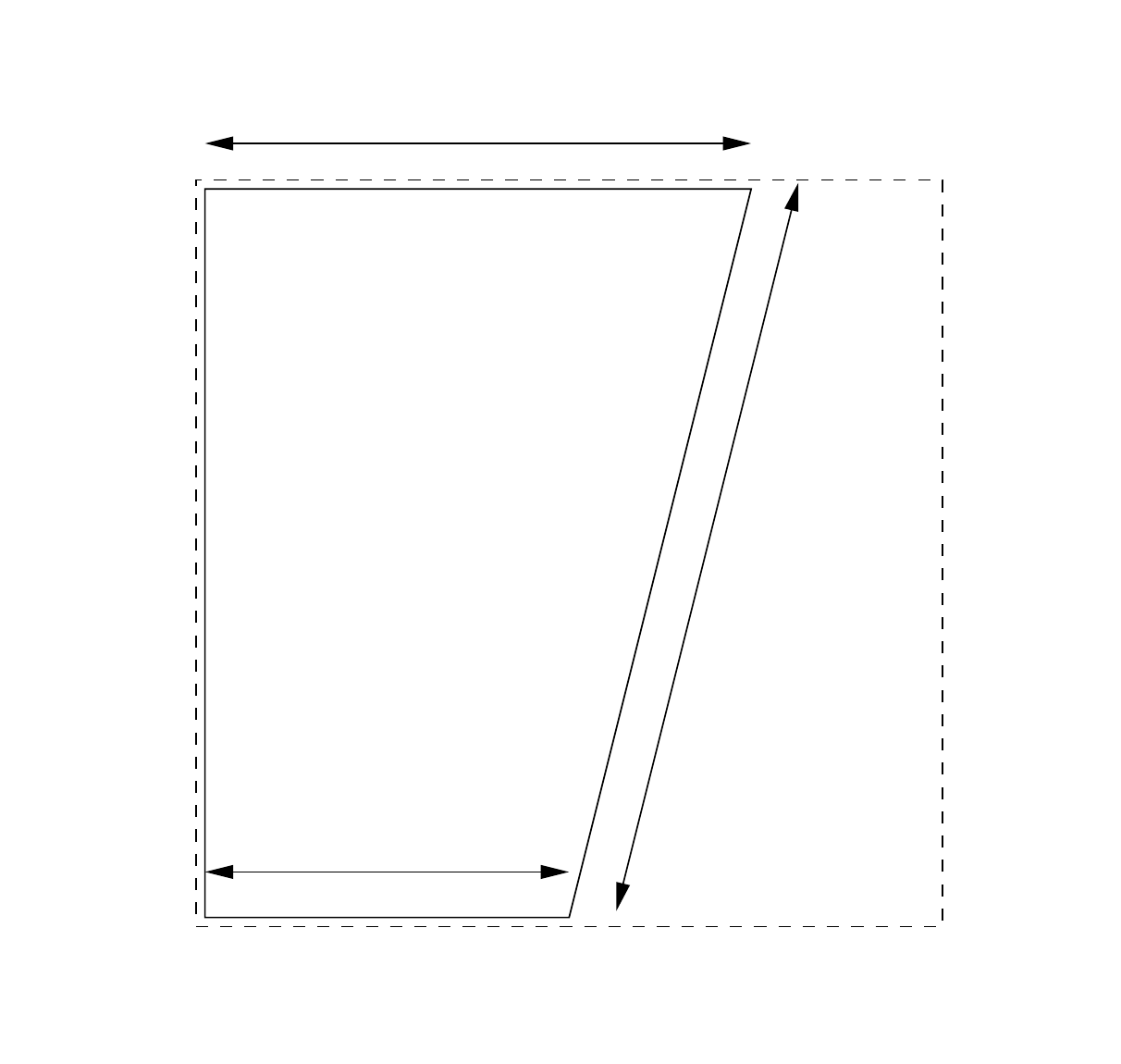}}
                \caption{Three examples of virtual cells: the cell
                  $\Omega_{i,j}$ is drawn with the dashed line, while
                  the corresponding virtual cell $\VirtDD$ is drawn
                  with the solid line. This virtual cell is a polygon
                  with at most five edges. Left: five edges. Middle:
                  three edges, while two virtual edges
                  $L_{i,j-\frac{1}{2}}^n$ and $L_{i-\frac{1}{2},j}^n$
                  have zero length. Right: four edges while the
                  virtual edge $L_{i+\frac{1}{2},j}^n$ has zero
                  length.}\label{CutCells}
\end{center}\end{figure}

Note that a difficult problem in the cut cell method is that cut cells
can take many different shapes (mainly in 3D), which can make the code
very complex. A key element of our approach is that all these
different shapes are treated generically by using this notion of
virtual cell with its 5 virtual edges. Indeed, all the cells are
treated in the same way, whatever their state (gas, solid, cut cell)
or shape. The different parameters of the three cell states are
summarized below, and we refer to figure~\ref{CutCells} for three
examples of cut cells:
\begin{itemize}
  \item gas cell: $\VirtDD^n=\Omega_{i,j}$, $\DxDD{\pm}=\Delta y$,
    $\DyDD{\pm}=\Delta x$, $\DwDD=0$, $\SVirtDD^n=\Delta x \Delta y$,
  \item solid cell: $\VirtDD^n=\emptyset$, $\DxDD{\pm}=0$,
    $\DyDD{\pm}=0$, $\DwDD=0$, $\SVirtDD^n=0$,
  \item cut cell: $\VirtDD^n$ is a part of $\Omega_{i,j}$, and the
    five lengths of the virtual cell $\DxDD{\pm}$, $\DyDD{\pm}$, and
    $\DwDD$ can take any value between 0 and $\Delta y$, $\Delta x$,
    and $\sqrt{\Delta x^2+\Delta y^2}$, respectively (see
    figure~\ref{CutCells}).
\end{itemize}
All these parameters are computed with a levelset technique, as explained in
appendix~\ref{LevelSet}.

\subsubsection{Control volumes}
\label{subsubsec:CV}

The notion of control volume
is essential to avoid the use of very small virtual cells that would
lead to prohibitively small time steps. The idea is to merge small
virtual cells with larger neighboring cells when their areas are smaller
than half of the area of a Cartesian cell.

\begin{figure}[h]\begin{center}
	\resizebox{0.95\textwidth}{!}{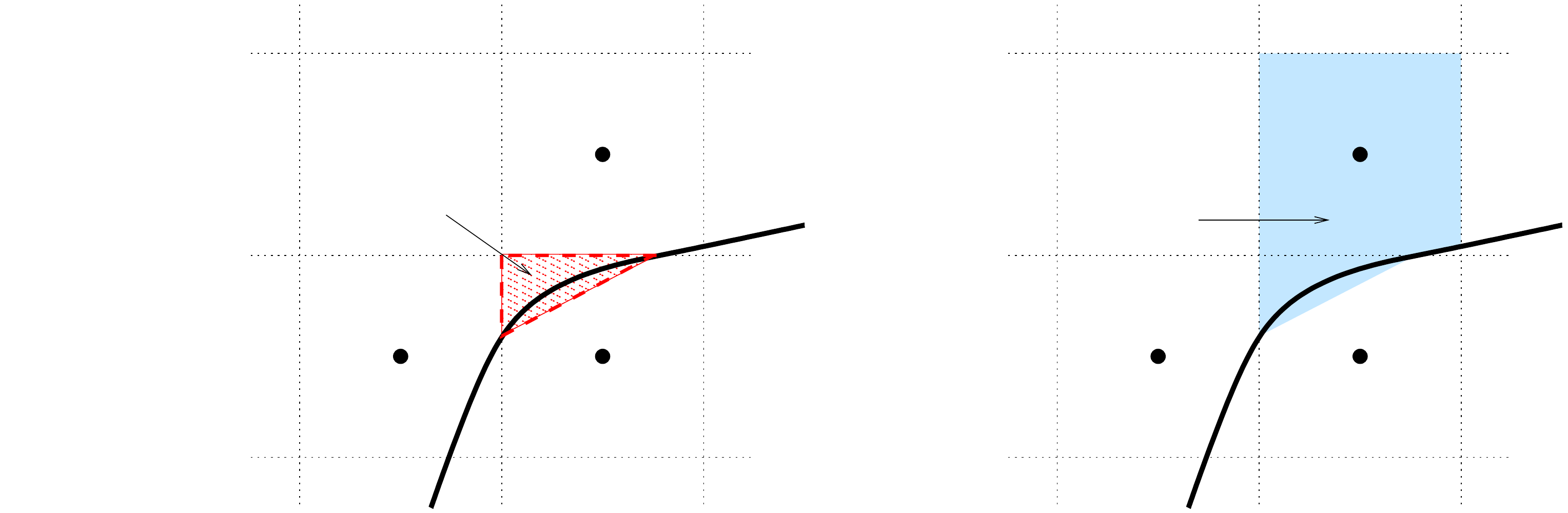}

\bigskip

	\resizebox{0.95\textwidth}{!}{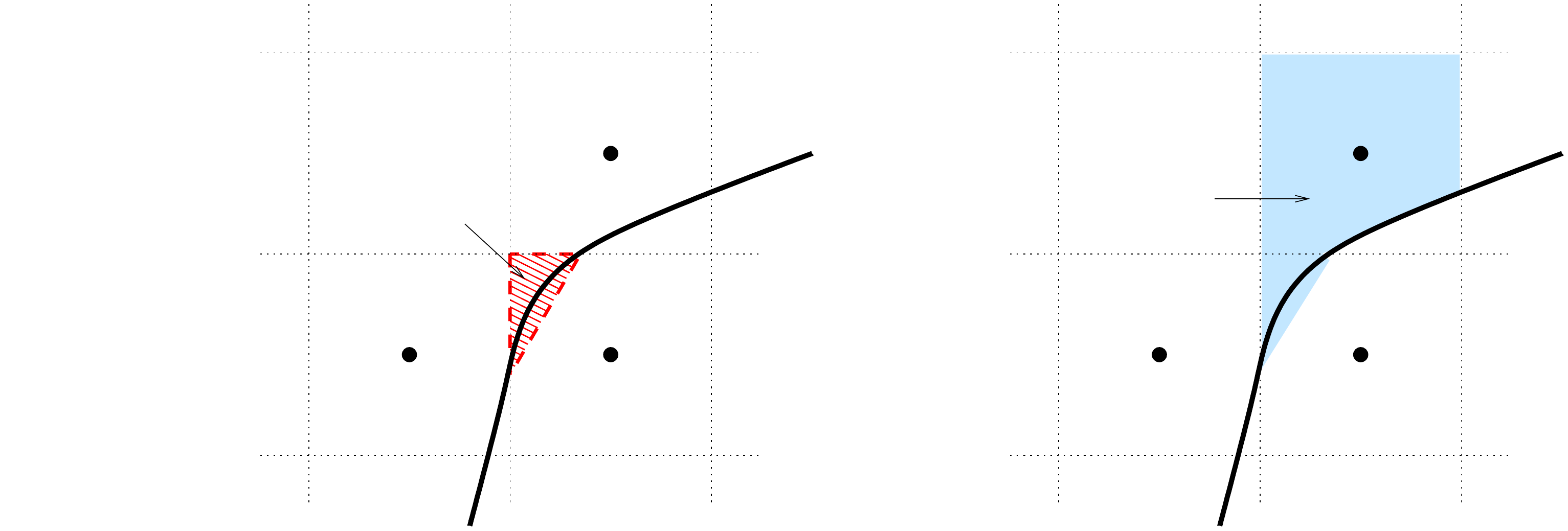}

\bigskip

	\resizebox{0.55\textwidth}{!}{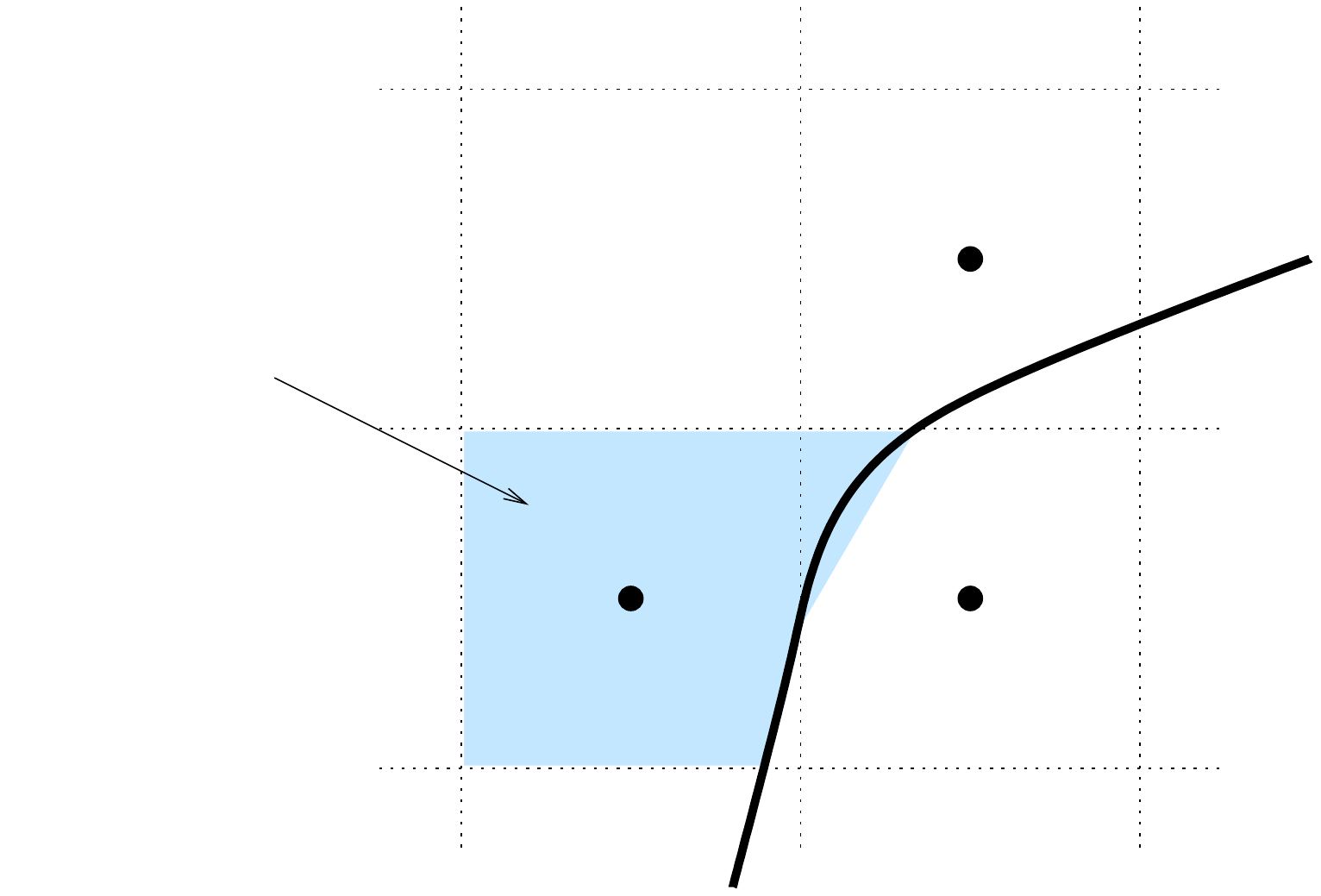}
\caption{Top-Left: 
the virtual cell $\VirtDD^n$ (hatched)
has to merge with $\VirtDDjp^n$
  because $|\LyDD{+}^{n}|>|\LxDD{-}^{n}|$. Top-Right: this results in the
  control volume $\CVolDD^n(t^n)$, in blue. Middle-Left: at time $t^{n+1}$,
  the solid
  boundary has slightly turned counter clock-wise and the
  virtual cell $\VirtDD^{n+1}$ is smaller. Middle-Right: the control volume
  at this time now is $\CVolDD^n(t^{n+1})$, in blue. Bottom: since the larger
  edge of the virtual cell $\VirtDD^{n+1}$ now is $|\LxDD{-}^{n+1}|$,
it has to merge with $\VirtDDim^{n+1}$, and this results in the new
control volume $\CVolDD^{n+1}(t^{n+1})$, which is different from
$\CVolDD^n(t^{n+1})$. \label{fig:ctrlvol}
}
\end{center}\end{figure}

The control volume is constructed by recursion: we look at a given cut
cell $\Omega_{i,j}$ whose center is inside the solid domain. The
corresponding virtual cell $\VirtDD^n$ necessarily has an area smaller than
$\frac12\Delta x\Delta y$, and it has to be merged with one of its non
solid neighboring cells. This cell is chosen by looking at the
largest non solid edge of $\VirtDD^n$: the neighboring cell that
shares the same edge is chosen for merging (for instance $\VirtDDjp^n$
in figure~\ref{fig:ctrlvol}, top). If the corresponding neighboring cell
has its center inside the gas domain, then the algorithm is stopped
and the resulting control volume contains two virtual cells. It
happens sometimes that the neighboring virtual cell is also too
small (its center is inside the solid domain too): in this case, the same
algorithm is used recursively for this virtual cell. This merging
procedure ensures that the area of the control volume is always
greater that $\frac12\Delta x\Delta y$.

It is convenient to denote by $\sigma_{i,j}^n$ the set of indices
$(i',j')$ such that all the virtual cells $\VirtDDp(t)$ are merged
together. For example, if $\VirtDD^n$ and $\VirtDDjp^n$ merge,
then $\sigma_{i,j}^n = \lbrace (i,j),(i,j+1)\rbrace$.

The previous algorithm defines the control volume at time $t^n$. For
$t > t^n$, the virtual cells change (since the solid boundary moves),
and the control volume as well. For $t$ between $t^n$ and $t^{n+1}$
the time dependent control volume $\CVolDD^n(t)$ is defined as
follows: 
\begin{equation}\label{CvolDD}
	\CVolDD^n(t) = \bigcup_{(i',j')\in \sigma_{i,j}^n} \VirtDDp(t).
\end{equation}
In other words, the set of virtual cells selected at time $t^n$ for
merging defines the control volume up to $t^{n+1}$. We point out that
if the shape of the virtual cells (and hence of the control volume)
can vary in time, the set $\sigma_{i,j}$ is fixed for
$t\in[t^n,t^{n+1}[$.

At time $t^{n+1}$, we have to take into account that there are new
virtual cells, some others have disappeared, and the shape of all of
them have changed: therefore, a new control volume, denoted by
$\CVolDD^n(t^{n+1})$, has to be constructed (by the previous recursive
algorithm).  We refer to figure~\ref{fig:ctrlvol} for an illustration
of this algorithm.

While the previous procedure might look complicated, note that most
of the virtual cells do not merge, and hence $\CVolDD^n(t)=\VirtDD(t)$ for most
of them.  

The area of the control volumes $\CVolDD^n(t^n)$ and
$\CVolDD^n(t^{n+1})$ are computed easily with
\begin{displaymath}
	\SCVolDD^n=\sum_{(i',j')\in \sigma^n_{i,j}} \SVirtDDp^{n} \quad\text{ and }
	\quad\SCVolDD^{n+1,\ast}=\sum_{(i',j')\in\sigma^n_{i,j}} \SVirtDDp^{n+1}
\end{displaymath}

Note that there is a kind of redundancy with this approach: indeed, in
the previous example, since $\VirtDD^n$ and $\VirtDDjp^n$ belong to
the same control volume, then the control volumes $\CVolDD^n(t)$ and
$\CVolDDjpun^n(t)$ are the same. However, this makes the
implementation much simpler, while the overhead of the computational
time is very small: indeed, the number of merged cut cells is very small
as compared to the number of gas cells.

\subsection{Numerical scheme}
\label{subsec:numscheme}
From now on, calculations are detailed with the reduced distribution function~$\fDD$.
The same analysis can be done with~$\gDD$.
The cut cell method is based on a finite volume scheme.
One time iteration (which will be divided into three steps) consists
in computing the average value~$\fVirtDD^{n+1}$ of the distribution
function over the virtual cell~$\VirtDD^{n+1}$ from the average value~$\fVirtDD^n$, defined by
\begin{equation}\label{fVirtDD}
	\fVirtDD^n = \frac{1}{\SVirtDD^n}
		\int_{\VirtDD^n} \fDD_{\iv}(t^n,\xDD)\dxDD.
\end{equation}
Similarly, $\fCVolDD^n$ and $\fCVolDD^{n+1,\ast}$ stand for the average values of the distribution function over the control volumes~$\CVolDD^n(t^n)$ and $\CVolDD^n(t^{n+1})$.
They are defined by
\begin{equation}\label{fCVolDD}
	\fCVolDD^n=\frac{1}{\SCVolDD^n}\int_{\CVolDD^n(t^n)} \fDD_{\iv}(t^n,\xDD)\dxDD
	\qquad\text{and}\qquad
	\fCVolDD^{n+1,\ast}=\frac{1}{\SCVolDD^{n+1,\ast}}
		\int_{\CVolDD^n(t^{n+1})} \fDD_{\iv}(t^n,\xDD)\dxDD.
\end{equation}

The first step of the method is the computation of~$\fCVolDD^n$
through the average values of $f$ over the virtual cells included in
$\CVolDD^n(t^n)$. Definitions~(\ref{CvolDD}), (\ref{fVirtDD}) and
(\ref{fCVolDD}) readily lead to
\begin{equation}\label{FusionDD}
	\fCVolDD^n := \frac{1}{\SCVolDD^n}
		\sum_{(i',j')\in\sigma^n_{i,j}} \SVirtDDp^n \fVirtDDp^n.
\end{equation}
% This procedure is essential because as mentioned before, virtual cells $\VirtDD^n$ can be infinitely small, and applying a finite volume scheme on these cells would lead to an infinitely restrictive CFL condition.

The second step is the time integration of the integral form of
the Boltzmann equation~(\ref{IntBoltzmannDD}) between $t^n$ and $t^{n+1}$:
this relation is applied by choosing the surface $S(t)$ as the control
volume $\CVolDD^n(t)$ and by using~(\ref{fCVolDD}). This gives
%The surface $S(t)$ in~(\ref{IntBoltzmannDD}) is defined by the control volume $\CVolDD^n(t)$ which yields
\begin{equation}\label{SchemeDD}
	\SCVolDD^{n+1,\ast}\,\fCVolDD^{n+1,\ast} - \SCVolDD^n\,\fCVolDD^n :=
	\int_{t^n}^{t^{n+1}} \Big( T_{i,j,\iv}(t)+Q_{i,j,\iv}(t)
        \Big)\, dt,
\end{equation}
where the transport and collision terms are defined by
\begin{equation} \label{eq-intTQ} 
\begin{split}
&	T_{i,j,\iv}(t)=-\int_{\partial \CVolDD^n(t)}
		(\vDD_{\iv}-\wDD(t))\cdot\nDD(t)\,f_{\iv}(t,\xDD)\dxUD,\quad
\\
&		Q_{i,j,\iv}(t) = \int_{\CVolDD^n(t)}
		\frac{1}{\tau(t,\xDD)}(\feqDD_{\iv}(t,\xDD)-f_{\iv}(t,\xDD))\dxDD.
\end{split}
\end{equation}
The transport integral can be computed as follows. First,
definition~(\ref{CvolDD}) implies that the integral over $\partial
\CVolDD^n(t)$ is the sum of the integrals over the contours of all the
virtual cells~$\VirtDDp(t)$ that merge into the control
volume~$\CVolDD^n(t)$.  Moreover, the velocity $\wDD\cdot\nDD$ is zero
for the four edges that fit with the Cartesian mesh lines, while this
velocity $\wDD$ is $\uwDD$ for the last edge of the cell, since it fits
with the solid boundary.  Finally, the transport integral is written
as:
\begin{equation}\label{TransportDD}\begin{split}
	T_{i,j,\iv}(t)
		&=-\sum_{(i',j')\in\sigma^n_{i,j}} \int_{\partial \VirtDDp(t)}
			(\vDD_{\iv}-\wDD(t))\cdot\nDD(t)\,f_{\iv}(t,\xDD)\dxUD\\
		&=-\sum_{(i',j')\in\sigma^n_{i,j}}\left[
			\int_{\LwDDp(t)}(\vDD_{\iv}-\uwDD(t))\cdot\nDD(t)\,f_{\iv}(t,\xDD)\dxUD
			 +\int_{L(t)}\vDD_{\iv}\cdot\nwDD(t)\,f_{\iv}(t,\xDD)\dxUD \right],
\end{split}\end{equation}
where $L=\LxDDp{+}\cup\LxDDp{+}\cup\LyDDp{-}\cup\LyDDp{-}$ is the
union of the four edges that fit with the Cartesian mesh lines. The
collision integral of~(\ref{eq-intTQ}) is readily approximated by
$\frac{\SCVolDD(t)}{\tau_{i,j}(t)}(\feqDD_{i,j,\iv}(t)-\fCVolDD(t))$.
The approximation of the time integral in the right-hand side
of~(\ref{SchemeDD}) will be detailed in the following sections.

The third step of the method is the computation of $\fVirtDD^{n+1}$,
the average value of $f$ in the new virtual cell $\VirtDD^{n+1}$.
This is done by distributing the value $\fCVolDD^{n+1,\ast}$ given
by~(\ref{SchemeDD}) to the cells $\VirtDD^{n+1}$ merged into the
control volume $\CVolDD^n(t^{n+1})$:
\begin{equation}\label{UpdateDD}
	\fVirtDD^{n+1}:=\fCVolDD^{n+1,\ast}.
\end{equation}

\begin{figure}[h]\begin{center}
	\subfigure[The boundary is moving from the right to the left.]
		{\resizebox{0.9\textwidth}{!}{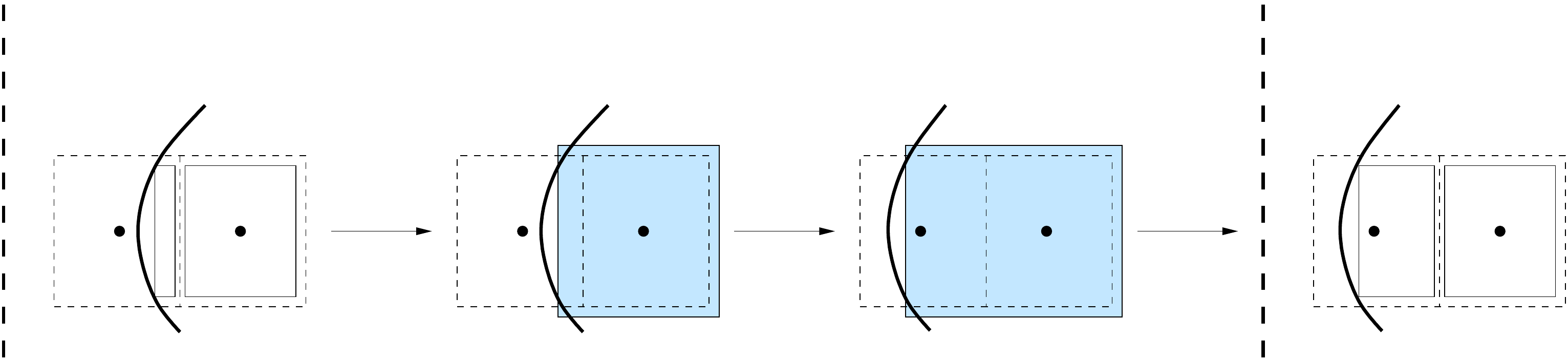}}
	\subfigure[Appearing cut cell: the boundary is moving from the right to the left.]
		{\resizebox{0.9\textwidth}{!}{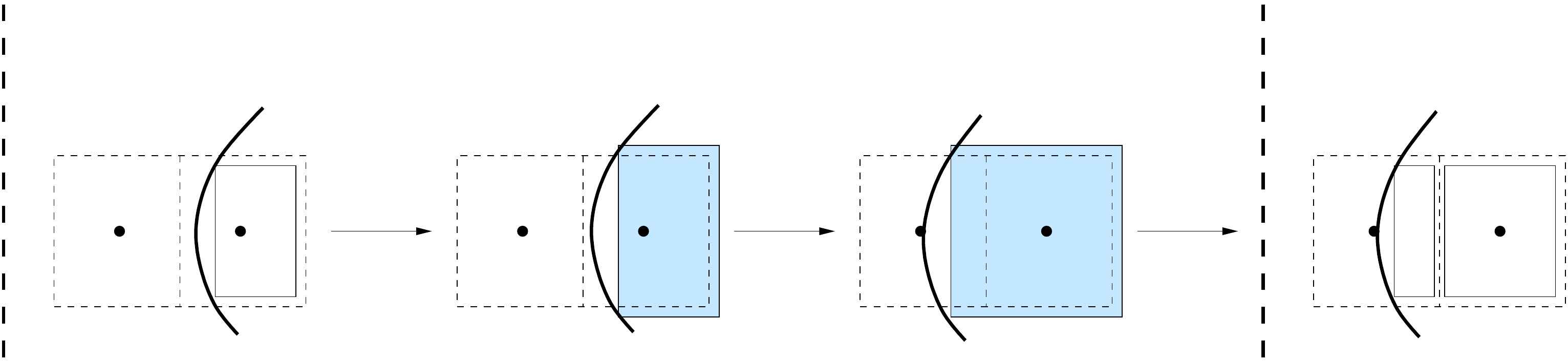}}
	\subfigure[Disappearing cut cell: the boundary is moving from
        the left to the right.]
		{\resizebox{0.9\textwidth}{!}{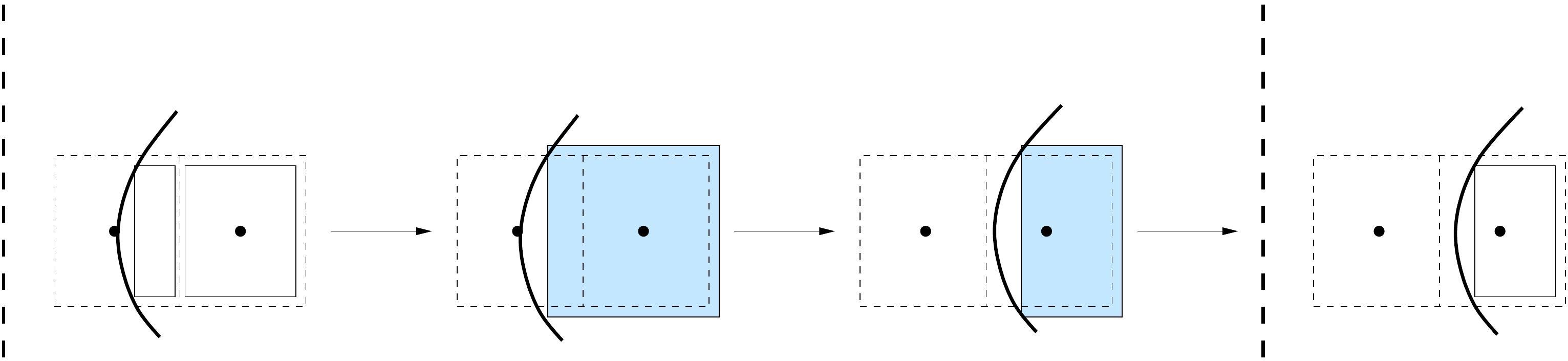}}
	\caption{Illustration of the three steps of the cut cell method. The cells $\Omega_{i,j}$ are drawn with the dashed line and the virtual cells $\VirtDD(t)$ with the solid line.
	The control volume $\CVolDD(t)$ is shaded.}\label{Steps}
\end{center}\end{figure}

A first summary of the cut cell method is given below:
\begin{enumerate}
\item The virtual cells~$\VirtDD(t)$ merge into some control
  volumes~$\CVolDD(t)$ and the values~$\fCVolDD^n$ are computed
  with~(\ref{FusionDD}).
	\item The numerical scheme~(\ref{SchemeDD}) is applied in order to computed the values~$\fCVolDD^{n+1,\ast}$.
	\item The values $\fVirtDD^{n+1}$ are updated with formula~(\ref{UpdateDD}).
\end{enumerate}
The three steps of the method are illustrated in figure~\ref{Steps}
for various situations. 
Note that because of the merging procedure, there is no issue of
appearing/disappearing gas cells: in other words, a small virtual cell
necessarily merges with a larger cell before it disappears, and
conversely, the average value of $f$ in a new appearing virtual cell
is naturally defined through steps 2 and 3. This ensures that the
method is conservative, see section~\ref{subsec:prop}.

It remains to explain how the time integral is approximated
in~(\ref{SchemeDD}) for step 2, which is done in
section~\ref{SubSubSecExplicit} and~\ref{subsubsec:semi}, and to
explain how the motion of the solid body is taken into account: this is
done in section~\ref{SubSecUw}. The complete scheme is summarized in
section~\ref{subsubsec:summary}

\subsubsection{First order explicit scheme}\label{SubSubSecExplicit}
A backward Euler method is applied in order to get a first order
approximation of the Boltzmann equation~(\ref{SchemeDD}).  This means
that the time integral in~(\ref{SchemeDD}) is approximated by the
rectangle rule. Using~(\ref{TransportDD}) and~(\ref{eq-intTQ}), we
find that relation~(\ref{SchemeDD}) becomes:
\begin{equation}\label{ExSchemeDD}\begin{aligned}
 	\fCVolDD^{n+1,\ast} = \frac{\SCVolDD^n}{\SCVolDD^{n+1,\ast}}
 	\fCVolDD^n &-\frac{\Delta t}{\SCVolDD^{n+1,\ast}}\sum_{(i',j')\in\sigma_{i,j}^n}\Big[
			\left(\FluxiDDp{+}^n-\FluxiDDp{-}^n\right) \\
			&\hspace*{80pt}+\left(\FluxjDDp{+}^n-\FluxjDDp{-}^n\right)
			+\FluxwDDp\Big] \\
 	& +\frac{\SCVolDD^n}{\SCVolDD^{n+1,\ast}} 
 	 	\frac{\Delta t}{\tau_{i,j}^n}(\feqDD^n_{i,j,\iv}-\fCVolDD^n),
\end{aligned}\end{equation}
where $\FluxiDD{\pm}$, $\FluxjDD{\pm}$ and $\FluxwDD$ are the
numerical fluxes across the five edges of the virtual cell that are
computed with a standard upwind scheme:
\begin{equation}\label{FluxDD}\begin{split}
    &\FluxiDD{+}^n := \DxDD{+}
    \bigl(\min(v_{\iv_1},0)\,\fDD_{i+1,j,\iv}^n+\max(v_{\iv_1},0)\,\fDD_{i,j,\iv}^n\bigr),\\
    &\FluxjDD{+}^n := \DyDD{+}
    \bigl[\min(v_{\iv_2},0)\,\fDD_{i,j+1,\iv}^n+\max(v_{\iv_2},0)\,\fDD_{i,j,\iv}^n\bigr],\\
    &\FluxwDD^n := \DwDD
    \bigl[\min([\vDD_p-\uwDD(t^n,\xwDD)]\cdot\nwDD(t^n,\xwDD),0)\,\fBoundDD
    \\
    &
    \hspace{14ex}+\max([\vDD_p-\uwDD(t^n,\xwDD)]\cdot\nwDD(t^n,\xwDD),0)\,\fDD_{i,j,\iv}^n\bigr],
\end{split}\end{equation}
where $\xwDD$ is the center of $\LwDD^n$.
It is recalled that $v_{\iv_1}$ and $v_{\iv_2}$ are the coordinates of
the $\iv^\text{th}$ microscopic velocity,
i.e. $\vDD_{\iv}=(v_{\iv_1},v_{\iv_2})$.  Moreover, $\nwDD(t^n,\xwDD)$ is the
outward normal to the edge $\LwDD^n$, that fit with the physical
boundary.  The computation of the velocity~$\uwDD(t^n,\xwDD)$ of the boundary is
detailed in section~\ref{SubSecUw}.
%%%% Eventuellement, dire "To simplify the notations, the dependance of $\uwDD$
%%%% and $\nwDD$ on $t^n$ and $\xwDD$ will be omitted in the following"
% Note that this velocity, as well
% as the normal vector $\nwDD$, are actually taken in the center $\xwDD$
% of the edge $\LwDD^n$, to simplify the notations, the dependance on
% $t^n$ and $\xwDD$ will be \fbox{pas utilise} in the following.
Finally the discrete boundary condition is similar to its continuous
form~(\ref{BoundaryDD}), that is to say:
\begin{equation}\label{BoundaryDDv}\begin{aligned}
	&f_w(t^n,\xwDD\in\Bound,\vDD_{\iv}\in\Vin) = \rhow\MaxwDD[1,\uwDD(t^n,\xwDD),\Tw], \\
	&g_w(t^n,\xwDD\in\Bound,\vDD_{\iv}\in\Vin) = \rhow\frac{RT}{2}\MaxwDD[1,\uwDD(t^n,\xwDD),\Tw], 
\end{aligned}\end{equation}
where $\rhow$ is given by
\begin{displaymath}
	\rhow = -\frac{\sum_{\vDD_{\iv}\in\Vout}(\vDD_{\iv}-\uwDD(t^n,\xwDD))\cdot\nwDD(t^n,\xwDD)
		\fDD_{i,j,\iv}^n\Delta v_x\Delta v_y}{
		\sum_{\vDD_{\iv}\in\Vin}(\vDD_{\iv}-\uwDD(t^n,\xwDD))\cdot\nwDD(t^n,\xwDD)
		\MaxwDD[1,\uwDD(t^n,\xwDD),\Tw]\Delta v_x\Delta v_y}. 
\end{displaymath}
Note that the boundary condition is defined only for the velocities
$\vDD_{\iv}\in \Vin=\{\vDD|(\vDD-\uwDD)\cdot\nwDD<0\}$, which is
compatible with the definition of the numerical boundary flux $\FluxwDD^n$ 
(see~(\ref{FluxDD})).

\subsubsection{First order semi-implicit scheme} \label{subsubsec:semi}

When the Knudsen number is very small, the previous
scheme~(\ref{ExSchemeDD}), which is explicit, is too expensive: the
CFL condition induces a time step which is of the order of the
relaxation time. It is now standard in kinetic theory to use instead
an implicit/explicit scheme (see for
instance~\cite{Pieraccini-Puppo-2007} for the BGK equation
and~\cite{DP_2014} for other methods). The idea is to
use an implicit scheme for the stiff collision part, while the
transport part is still approximated by an explicit scheme.
This kind of scheme can be easily extended to the cut cell method. For
instance, the simplest first order semi-implicit scheme is
\begin{equation}\label{ImpSchemeDD}\begin{aligned}
 	\fCVolDD^{n+1,\ast} = \frac{\SCVolDD^n}{\SCVolDD^{n+1,\ast}}
 	\fCVolDD^n &-\frac{\Delta t}{\SCVolDD^{n+1,\ast}}\sum_{(i',j')\in\sigma_{i,j}^n}\Big[
			\left(\FluxiDDp{+}^n-\FluxiDDp{-}^n\right) \\
			&\hspace*{80pt}+\left(\FluxjDDp{+}^n-\FluxjDDp{-}^n\right)
			+\FluxwDDp\Big] \\
 	& +\frac{\Delta t}{\tau_{i,j}^n}(\feqDD^{n+1,\ast}_{i,j,\iv}-\fCVolDD^{n+1,\ast}),
\end{aligned}\end{equation}
where the numerical fluxes are still computed with~(\ref{FluxDD}).
The equilibrium function $\feqDD^{n+1,\ast}_{i,j,\iv}$ depends on the
macroscopic quantities that have to be computed before
$\fCVolDD^{n+1,\ast}$:  this can be done by summing~(\ref{ImpSchemeDD})
over the discrete velocities, so that the collision term vanishes,
which gives an explicit relation for these macroscopic quantities.

%\subsubsection{Second order explicit scheme}

\subsubsection{Motion of the solid body}\label{SubSecUw}

The motion of the solid body is taken into account in the scheme by
the variation of the area of the control volume (from $\SCVolDD^n$ to
$\SCVolDD^{n+1,\ast}$) and by the velocity $\uwDD(t^n,\xwDD)$ of the
solid boundary, see~(\ref{ExSchemeDD})
and~(\ref{FluxDD}). These quantities are computed as follows.

Let $\xSolDD(t)$ and $\ThetaSol(t)$ be the coordinates of the center
of mass and the inclination of the solid body.  Its
translational and rotational velocities are then denoted by
$\uSolDD(t)$ and $\dotThetaSol(t)$. The motion of the solid body,
with mass $\mass$ and moment of inertia $\inertia$, is
modeled by the Newton's laws of motion that are discretized as
follows:
\begin{align}
& \xSolDD^{n+1} = \xSolDD^n +\Delta t \,  \uSolDD^n
	\quad \text{ and } \quad
\ThetaSol^{n+1} = \ThetaSol^{n} + \Delta t \,   \dotThetaSol^n,  
\label{motionDD_posit}\\
& 
\uSolDD^{n+1} = \uSolDD^n +\Delta t  \,  \FSolDD^n/\mass
\quad \text{ and } \quad
\dotThetaSol^{n+1} = \dotThetaSol^{n} +\Delta t \,   \Torque^n/\inertia, \label{motionDD_vit}
\end{align}
where $\FSolDD$ and
$\Torque$ are
the force and torque exerted by the gas on the solid body. They can be
computed by using the stress tensor $\pwDD$ at the boundary with the
formula
\begin{equation*}
\FSolDD=\int_{\partial \Omega_g}\Sigma_w\nwDD \, dl \quad \text{ and }
\quad 
\Torque=\int_{\partial \Omega_g}(\xDD-c)\times (\Sigma_w\nwDD)\, dl.
\end{equation*}
These relations can be approximated by any quadrature formula, and we
find it convenient to use a summation over all the cells of the
computational domain to avoid too many tests. This yields:
\begin{align}
	\FSolDD^n & =\sum_{i=1}^{N_x}\sum_{j=1}^{N_y}
		\pwDD(t^n,\xwDD) \nwDD(t^n,\xwDD) \DwDD,\label{eq-force}  \\ 
	\Torque^n & = \sum_{i=1}^{N_x}\sum_{j=1}^{N_y}
		\left((\xwDD-\xSolDD^n)\times(\pwDD(t^n,\xwDD)\nwDD(t^n,\xwDD))\right)
		\DwDD\label{eq-torque},
\end{align}
Since $\DwDD$ is non-zero only for solid edges of cut cells, these
formula are consistent approximations of the previous definition.

Moreover, while 
% We remind that $\pwDD$ and $\nwDD$ are time and space dependent and
% should be written $\pwDD(t^n,\xwDD)$ and $\nwDD(t^n,\xwDD)$ here.
the stress tensor is defined by~(\ref{MacroscopicDDv}), the boundary
condition has to be taken into account to define the distribution of
incoming velocities, and we set
\begin{equation}\label{StressTensorDD}\begin{split}
	\pwDD(t^n,\xwDD)&=\sum_{\vDD_{\iv}\in\Vin}(\vDD_{\iv}-\uwDD(t^n,\xwDD))\otimes
						(\vDD_{\iv}-\uwDD(t^n,\xwDD))\fBoundDD \,\Delta v_x\Delta v_y\\
			& \quad +\sum_{\vDD_{\iv}\in\Vout}(\vDD_{\iv}-\uwDD(t^n,\xwDD))\otimes
						(\vDD_{\iv}-\uwDD(t^n,\xwDD))\fCVolDD^n \,\Delta v_x\Delta v_y.
\end{split}\end{equation}
The boundary condition $\fBoundDD$ is defined by~(\ref{BoundaryDDv}).
The new velocity of the wall is finally computed with
\begin{equation}\label{uwDD}
  \uwDD(t^{n+1},\mathbf{r}^{n+1}_{i,j}) = \uSolDD^{n+1} + (\mathbf{r}^{n+1}_{i,j}-\xSolDD^{n+1})^\perp \dotThetaSol^{n+1},
\end{equation}
where $a^\perp$ is the vector obtained after a rotation of 90 degrees
of any vector $a$ in the counter-clockwise sense.

\subsubsection{Summary of the numerical scheme}
\label{subsubsec:summary}

For the convenience of the reader, the different steps of the
complete numerical scheme are summarized below. 

We assume that, at time $t^n$, all the following quantities are known:
the average value of the distribution function $\fVirtDD^n$ in each
virtual cell $\VirtDD^{n}$, the parameters of position
($\xSolDD^n$,$\ThetaSol^n$) and velocity ($\uSolDD^{n}$,
$\dotThetaSol^{n}$) of the solid body, and hence the wall velocity
$\uwDD(t^n,\xwDD)$.  One time iteration of the numerical scheme is
decomposed into the following 7 steps:
\begin{enumerate}

\item The position  of the solid body at $t^{n+1}$ is computed with~(\ref{motionDD_posit}).

\item The virtual cells are arranged into control volumes
  $\CVolDD^{n}(t_{n})$ following the rule given in
  section~\ref{subsubsec:CV}. The distribution function is averaged
  over the control volumes $\CVolDD^n(t^n)$, see~(\ref{fCVolDD}).

\item The virtual cells and control volumes are moved according to the
  new position computed at step 1. The areas $\SCVolDD^{n}$ and
  $\SCVolDD^{n+1,\ast}$ are computed by using a level set method (see
  appendix~\ref{LevelSet} for more details).

\item The value of the distribution at each solid boundaries is
  computed through boundary condition~(\ref{BoundaryDDv}).

\item  The stress tensor $\Sigma_w$ at each solid boundaries is computed
  with~(\ref{StressTensorDD}), which gives force $\FSolDD^n$ and
  torque $\Torque^n$ with~(\ref{eq-force}) and~(\ref{eq-torque}). Then
the translational and rotational velocities are computed at time
$t^{n+1}$ by the discrete Newton laws~(\ref{motionDD_vit}), and
finally the new wall velocity is computed with~(\ref{uwDD}).

\item Scheme~(\ref{ExSchemeDD}) is used to pass
  from $\fCVolDD^n$ (the average value of $f$ at time $t^n$ in the
  control volume $\CVolDD^n(t^n)$) to $\fCVolDD^{n+1,\ast}$ (the
  average value of $f$ at time $t^{n+1}$ in the control volume
  $\CVolDD^n(t^{n+1})$).

\item The values $\fVirtDD^{n+1}$ of $f$ at time $t^{n+1}$ in each
  virtual cell that are merged into the control volume $\CVolDD^n(t^{n+1})$ are
  updated with~(\ref{UpdateDD}).

\end{enumerate}

%=========
   \subsection{Properties of the scheme}
   \label{subsec:prop}
%=========

\subsubsection{Positivity}

   Standard arguments show that the explicit scheme~(\ref{ExSchemeDD})
   preserves the positivity of the solution if $\Delta t$ satisfies
   the following CFL condition:
\begin{equation}  \label{eq-CFL}
	\Delta t \left( 
		\max_{i,j}\left[\frac{1}{\tau_{i,j}^n}\right] +
		\max_{i,j,\iv}\left[\frac{\phi_{i,j,\iv}^n}{\SCVolDD^n}\right]
		\right)\leq 1,
\end{equation}
where 
\begin{equation}\label{phiCFL}
	\begin{split}
	\phi_{i,j,\iv}&=\sum_{(i',j')\in\sigma_{i,j}^n}
\left(|\LxDDp{+}|v_{\iv_1}^+-|\LxDDp{-}|v_{\iv_1}^-  %\right. \\
 + |\LyDDp{+}|v_{\iv_2}^+-|\LyDDp{-}|v_{\iv_2}^- \right.\\
& \left. \quad \qquad \qquad + |\LwDDp|((\vDD_p-\uwDD(t^n,\xwDD))\cdot \nwDD(t^n,\xwDD))^+
\right)
\end{split}
\end{equation}
and $v^+=\max(v,0)$, $v^-=\min(v,0)$. For a correct description of
the flow, it is necessary that the discrete velocity grid contains the
solid velocity $\uwDD$ at any time (since the diffuse boundary
condition produces particles with $\uwDD$ as mean velocity). Under
this assumption, it can be proved that $\phi_{i,j,\iv}\leq
C(v_{x,\max}/\Delta x + v_{y,max}/\Delta y)$,
where $C=1/2$. This gives a simpler CFL condition, and in
practice, it is relaxed by taking $C=0.9$ without any stability problems

For the semi-implicit scheme~(\ref{ImpSchemeDD}), a similar condition
can be found, which is independent of $\tau^n_{i,j}$.

\subsubsection{Conservation}

Since our scheme is a finite volume method in which a conservative
reflexion boundary condition is applied to compute the numerical
fluxes at each solid edge, it is naturally conservative. This is
proved below for the explicit scheme~(\ref{ExSchemeDD}). 

Let $M^n$ be
the total mass of gas in the gas domain $\Omega_g$ at time $t^n$. 
It is convenient to write 
the total mass at time $t^{n+1}$ as
\begin{equation}\label{def-masse}
	M^{n+1} = 
		\sum_{i,j=1}^{N_x,N_y}\delta_{i,j}^n\SCVolDD^{n+1,\ast}
		\sum_{p=0}^{N^2-1}\fCVolDD^{n+1,\ast}\Delta v_x\Delta v_y,
\end{equation}
where  $\delta^n_{i,j}=1$ if $\xDD_{i,j}\in \Omega_g$ and $0$ else:
this function allows to take into account merged cells of a same
control volume only once. Indeed, there is only one couple of indices
$(i',j')$ in $\sigma_{i,j}^n$ for which
$\xDD_{i',j'}$ is inside the gas domain.

Then~$\fCVolDD^{n+1,\ast}$ is replaced by its value given
by~(\ref{ExSchemeDD}) and we get
\begin{displaymath}\begin{split}
	M^{n+1}  
%= M^n + \sum_{i,j=1}^{N_x,N_y}\delta_{i,j}^n\sum_{i'\in\sigma_i(t^n)}
% 		 	\sum_{p=0}^{N^2-1} \Delta t\,\mathcal{F}^n_{i',j,k}\Delta v_x\Delta v_y			%	
		= M^n + \sum_{i,j=1}^{N_x,N_y}\sum_{p=0}^{N^2-1} \Delta t\, 
			\mathcal{F}^n_{i,j,k}\Delta v_x\Delta v_y.
\end{split}\end{displaymath}
Indeed, opposite fluxes across same Cartesian edges cancel out, the
velocity sum of the collision operator is zero, and there remains only
the numerical fluxes $\mathcal{F}^n_{i,j,k}$ across the solid edge of
cut cells. By using the boundary condition, the velocity sum of such
fluxes gives
\begin{equation*}
\begin{split}
	\sum_{p=0}^{N^2-1}\mathcal{F}^n_{i,j,p}\Delta v_x\Delta v_y =& 
	\sum_{\vDD_{\iv}\in\Vin}
        (\vDD_{\iv}-\uwDD(t^n,\xwDD))\cdot\nwDD(t^n,\xwDD)\,\fBoundDD\DwDD
        \Delta v_x\Delta v_y \\
	 & + \sum_{\vDD_{\iv}\in\Vout} (\vDD_{\iv}-\uwDD(t^n,\xwDD))\cdot\nwDD(t^n,\xwDD)\,
	 \fCVolDD^n \DwDD\Delta v_x\Delta v_y\\ 
 & = 0.
\end{split}
\end{equation*}
This shows that $M^{n+1}=M^n$ and concludes the proof.

Finally, note that it is standard from the positivity and conservation
properties to conclude that the scheme if $L^1$ stable.

\section{Numerical results}\label{SubSecExamplesDD}

\subsection{Translational motion under radiometric effect}
\begin{figure}[h]\begin{center}
	\subfigure[Experimental set up of the moving plates. The computational domain is drawn with a doted line.]
		{\resizebox{0.4\textwidth}{!}
		{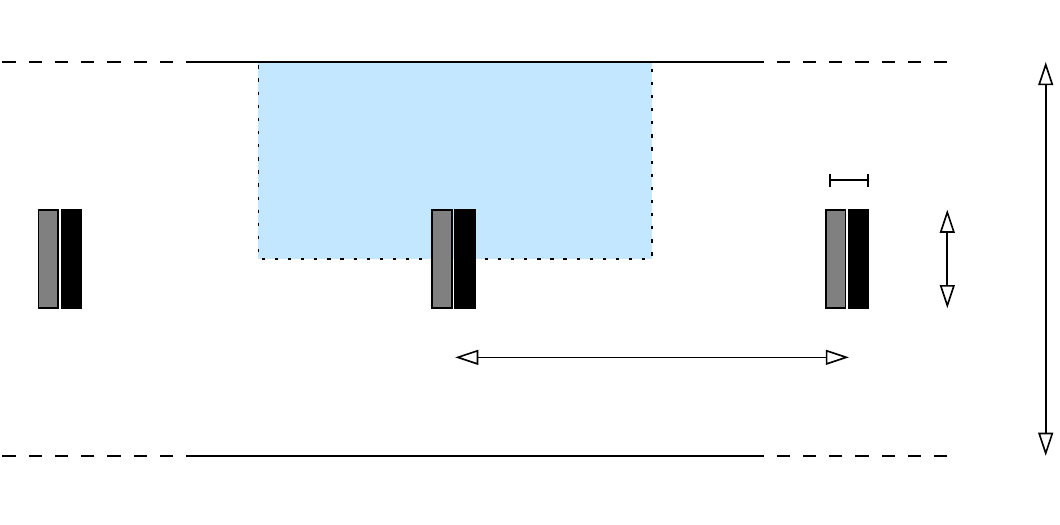}\label{Taguchi_1}}
	\hspace*{0.09\textwidth}
	\subfigure[Evolution of the velocity of the plates for three
        Knudsen numbers.]
		{\resizebox{0.49\textwidth}{!}
	 {\includegraphics[width=0.48\textwidth]{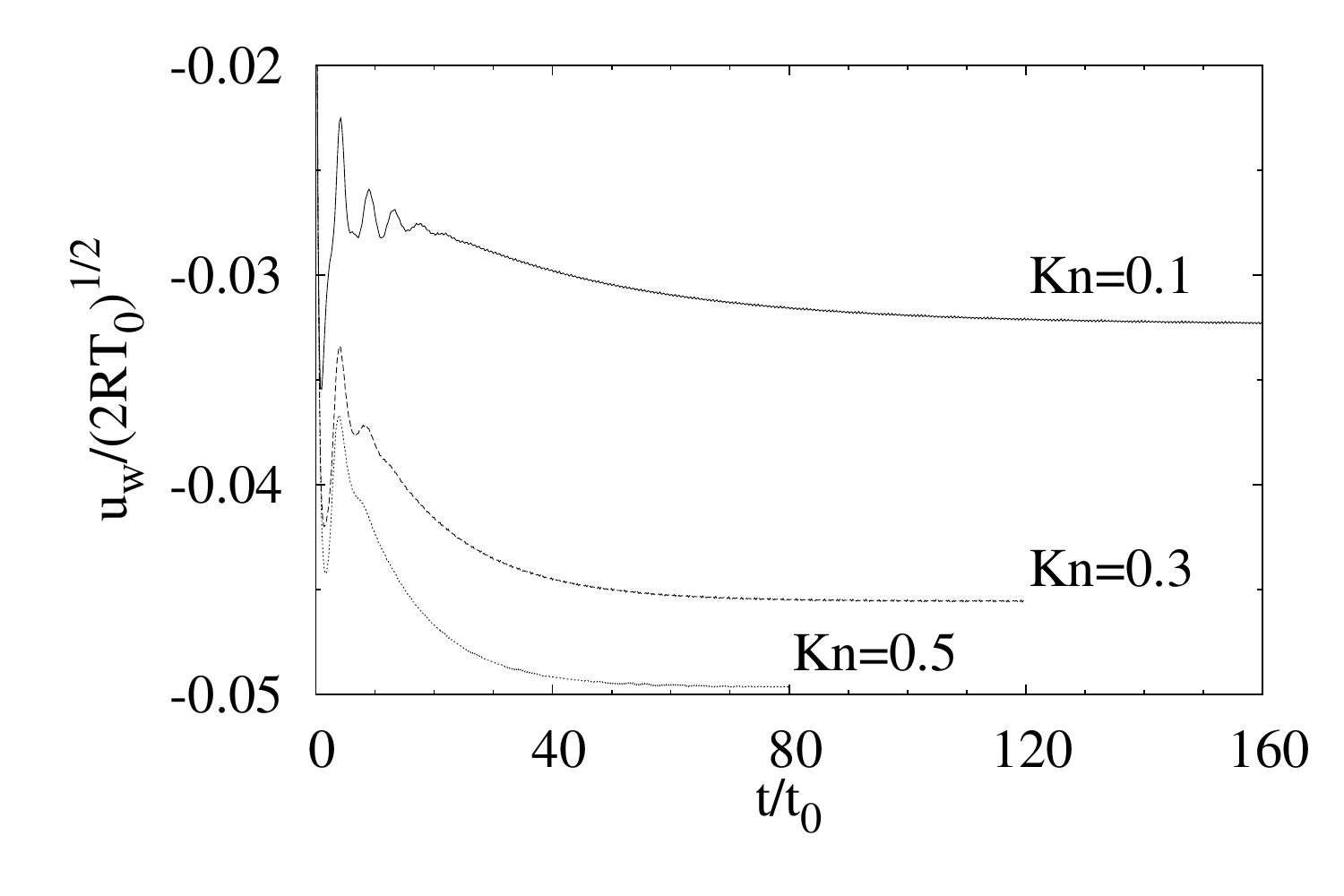}}\label{Taguchi_test}}
	\caption{Translational plates under the radiometric effect}
\end{center}\end{figure}

An infinite set of thin plates of height~$D$ is located in an infinite
channel of width~$4D$.  The distance between two plate centers is
$4D$.  This experiment has been numerically investigated
in~\cite{Taguchi-Aoki-2012b}, where the plates are infinitely thin
(their thickness is zero).  The temperature of the right side of a
plate is twice as the temperature of its left side.  This induces a
force that can be interpreted as the difference between the
radiometric force and gas friction, and make the plates move.  We
point out that all the plates move at the same velocity.  After a
while, the radiometric force is balanced with gas friction and the
total force decreases to zero.  At this time, the plates have reached
their stationary velocity.  In~\cite{Taguchi-Aoki-2012b} the ES-BGK
model is used and simulations are made in a moving reference frame in
which the channel velocity is positive while the plates are
motionless.  The total force applied on the plates is zero for a
specific channel velocity which corresponds to the stationary velocity
of the plates.  In our work, the BGK model is used to described the
behavior of the gas and plates are $D/10$ thick~(see
Fig.~\ref{Taguchi_1}).  The transient motion of the plates
(i.e. plates velocity before the stationary state is reached) can be
simulated with the cut cell method in a fixed frame of reference.  In
that case, the motion depends on the mass of a plate.  If we denote by
$\rho_0$ the initial density of the gas, this mass is set to $m=\rho_0
D^2/2$.  The computational domain is a rectangle of size $4D\times 2D$
that describes the upper part of the channel (Fig.~\ref{Taguchi_1}).
The left and right boundary conditions are periodic so as to simulate
the infinite channel.  The bottom of the computational domain fits
with the center of the channel: specular-reflection boundary condition
is applied to take into account the symmetry.  At last, the top of the
domain corresponds to the wall of the channel, the standard diffuse
boundary condition is used.  During a whole simulation, there is
always one and only one plate in the computational domain: when a
plate goes out, an other one comes in.  Note that solid and
cut cells only appear near this plate.  To take its temperature into
account, the boundary of the plate is modeled with the
diffuse-reflection condition.  Finally, the relaxation time and
Knudsen number are defined by the relations:
\begin{displaymath}
	\tau=\frac{\mu}{\rho R T_0}
	\quad\text{ and }\quad
	\Kn = \left[\frac{2}{\sqrt{\pi}}\frac{\sqrt{2RT_0}}{\rho_0RT_0/\mu}\right]\Big/ D,
\end{displaymath}
where $\mu$ is the viscosity of the gas and $T_0$ the initial temperature of the gas.

Our simulations have been done for a wide range of Knudsen numbers
from $10^{-3}$ to $1$.  Converged results are obtained for a velocity
grid that contains from $20^2$ to $40^2$ points (this depends on the
Knudsen number) and for a spatial mesh made up of $400\times 200$
cells, which means that a plate encloses 10 cells.  For coarser grids,
the number of cells enclosed in the plate is too small to capture the
shape of its edges with enough accuracy.
% As a result, force and velocity are slightly oscillating around the
% equilibrium values for coarse meshes.
As expected, the magnitude of the velocity of the plates increases
until they reach their final velocity, as illustrated for three
different Knudsen number on figure~\ref{Taguchi_test}.

\begin{figure}[h]\begin{center}
		{\resizebox{0.49\textwidth}{!}
	 {\includegraphics[width=0.48\textwidth]{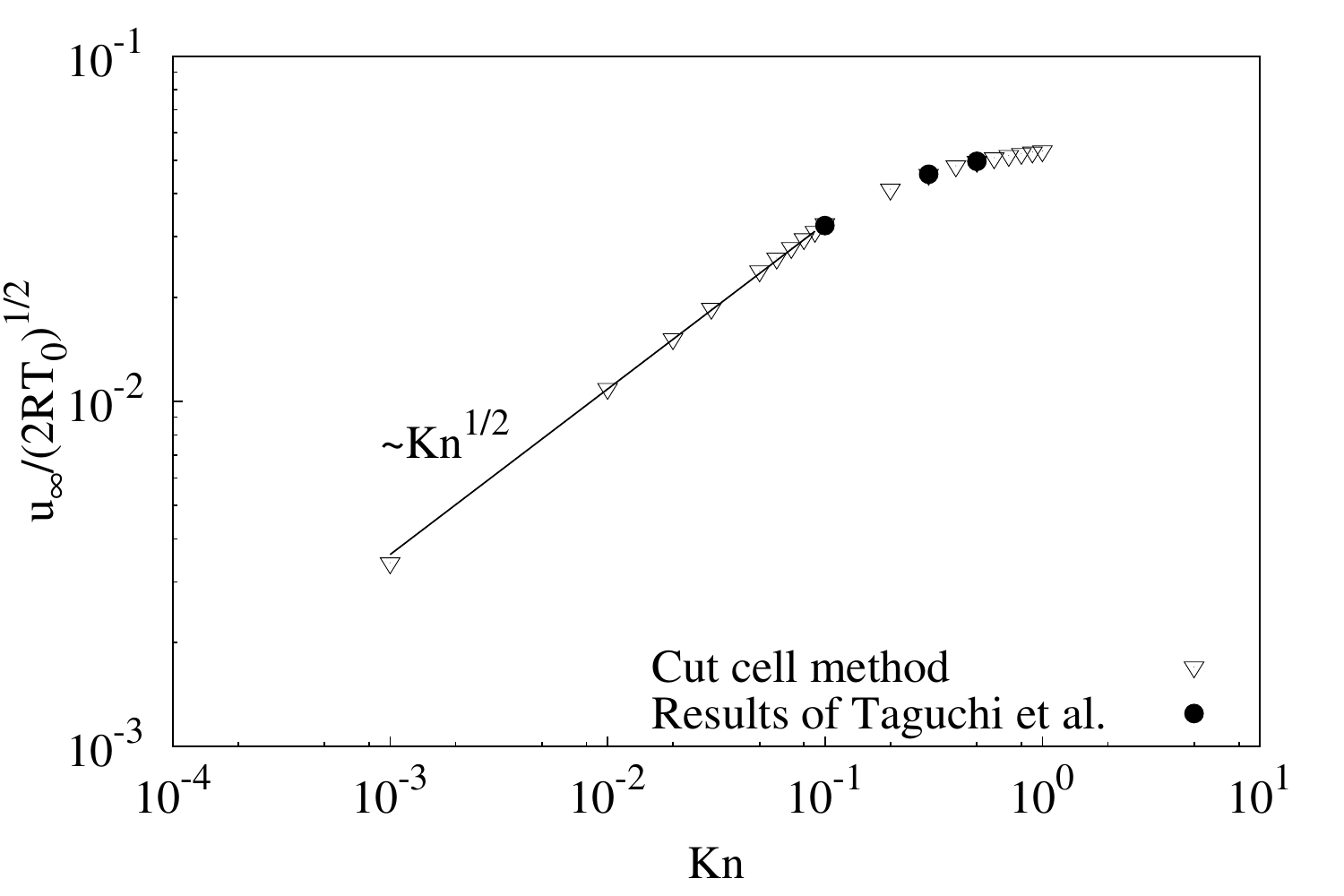}}}
	\caption{Translational plates under the radiometric effect:
          stationary velocity of the plates as a function of the
          Knudsen number, comparison with Taguchi et
          al.~\cite{Taguchi-Aoki-2012b}.  \label{Taguchi_2modif}}
\end{center}\end{figure}

The variation of the stationary velocity of the plates is plotted as
function of $\Kn$ on figure~\ref{Taguchi_2modif}.  For small Knudsen
numbers, the final velocity seems to be proportional to $\sqrt{\Kn}$
while this velocity tends to a constant for high Knudsen number.
These simulations show a good agreement between our results and the
results obtained in~\cite{Taguchi-Aoki-2012b}.

%%%%%%%%%%%%%%%%%%%%%%%%%%%%%%%%%%%%%%%%%%%%%%%%%%%%%%%%%%%%%%%%%%%%%%%%
\subsection{The Crookes radiometer}
\label{subsec:crookes}
The Crookes radiometer was invented by Crookes in
1874~\cite{Crookes-1874}: it is a glass globe containing four vanes
immersed in a low pressure gas.  Each vane has one black side and one
shiny side, and when the globe is exposed to light, the vanes
rotate. This was first understood as a rarefied gas dynamics effect by
Reynolds~\cite{Reynolds-1879}, but there are still discussions on the
order of magnitudes of the forces involved in this device. We refer to
the recent review of Ketsdever et
al.~\cite{Ketsdever-Gimelshein-Selden-2012} for historical
details. Recently, numerical simulations improved the understanding of
the radiometric effect, like
in~\cite{TK_2012,Taguchi-Aoki-2012b,
  Selden-Ngalande-Gimelshein-Muntz-Alexeenko-Ketsdever-2009}.

The dynamical acceleration process of the vanes has been recently
studied in~\cite{Chen-Xu-Lee-2012}: it uses the unified gas-kinetic
scheme combined with a moving mesh
approach~\cite{Chen-Xu-Lee-Cai-2012} to simulate the motion of the
vanes.  In this case, the moving mesh approach is very convenient
because the initial mesh just rotates without distortions.  In
this section, we show that the same results as~\cite{Chen-Xu-Lee-2012}
can be obtained with our cut cell method.

\begin{figure}[h]\begin{center}
	\subfigure[Experimental set up, the computational domain is drawn with a doted line.]
		{\resizebox{0.38\textwidth}{!}
		{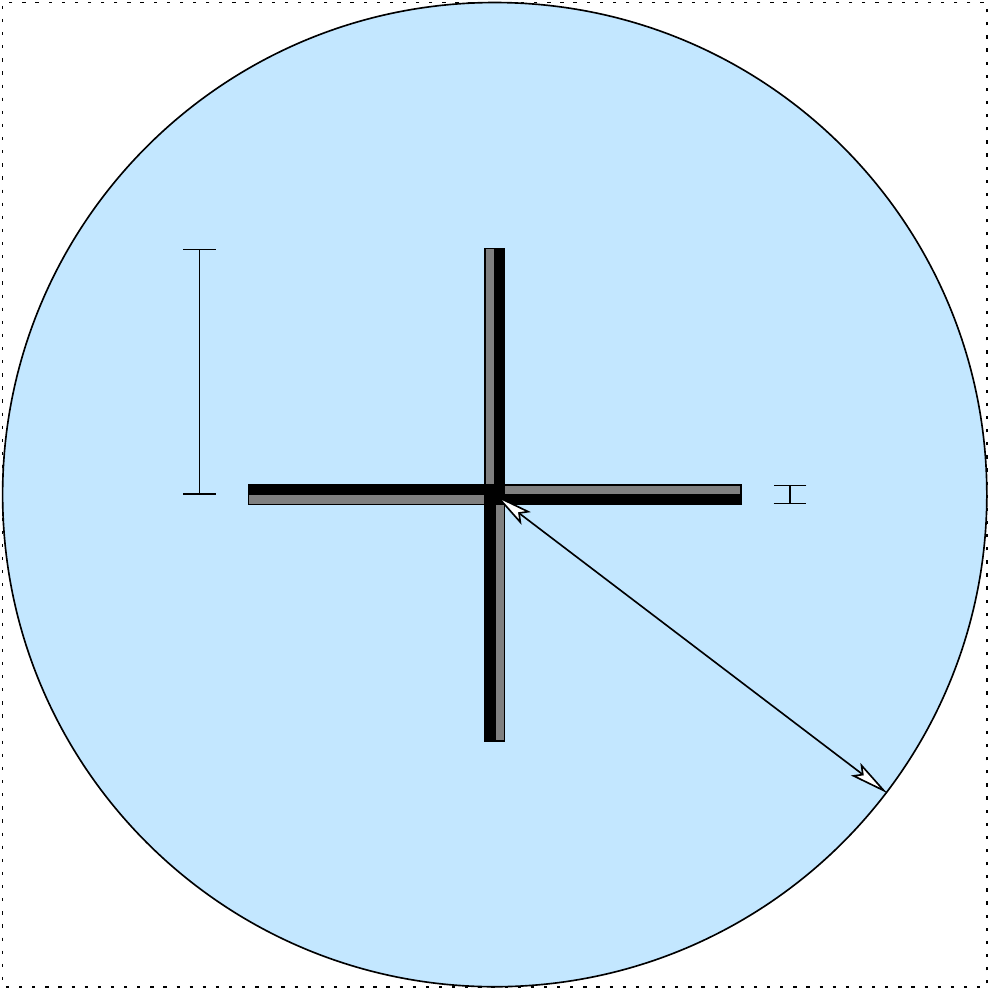}\label{Radiometer_1}}
%	\hspace*{0.08\textwidth}
	\subfigure[Radial velocity of the vanes as a function of time.]
	{\includegraphics[width=0.55\textwidth]{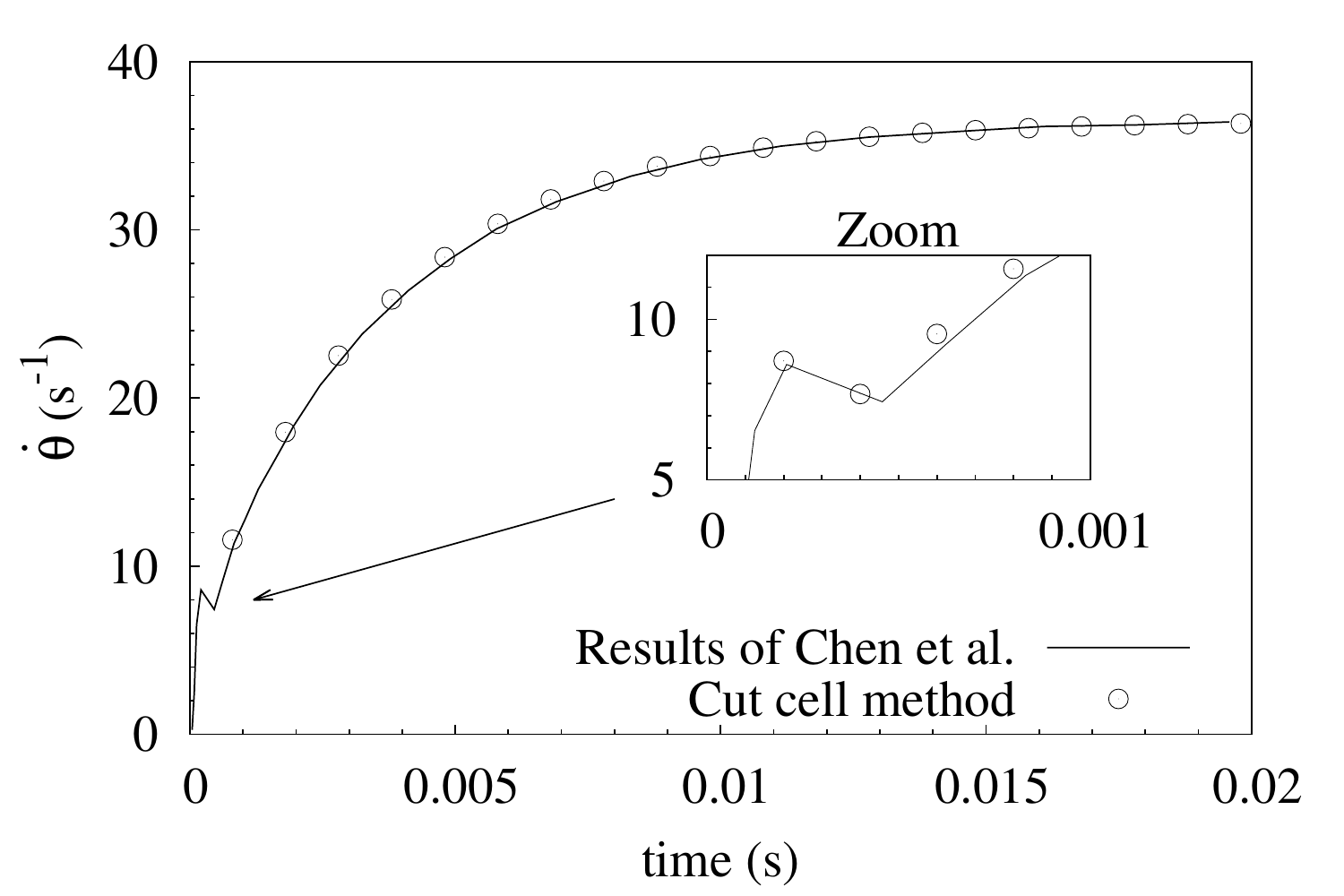}\label{Radiometer_2}}
	\caption{The 2D Crookes radiometer.}
\end{center}\end{figure}

Since this test is for illustrating our 2D method, the device
simulated here is a 2D radiometer composed of endless vanes immersed
in an unbounded cylinder of radius $R=20$cm.  The length of a vane is
$L=0.1$cm, its thickness is $l=0.01$cm (figure~\ref{Radiometer_1}) and its rotational
moment of inertia is $J=4.9\times 10^{-9}$kg$\cdot$m$^2$.  The
temperature $T_h$ of the black side of the vane is supposed to be
higher than the temperature $T_c$ of its shiny side.  These
temperatures are set to $T_h=400$K and $T_c=350$K while the
temperature $T_0$ of the globe is equal to $T_0=300$K.  Note that all
the boundary conditions are computed with diffuse reflection conditions.
In order to compare our results to~\cite{Chen-Xu-Lee-2012}, we take
the same Shakhov relaxation model with $Pr=2/3$.  Moreover, the
relaxation time is computed by the equation
\begin{equation}\label{RelaxationTime}
	\tau=\frac{\mu}{P}\cdot\left(\frac{T}{T_0}\right)^\omega,
\end{equation}
where $\mu$ is determined by the hard sphere model for argon which
yields $\mu=1.678\times 10^{-5}$Nsm$^{-2}$ and $\omega=0.68$.
Finally, we set the initial density $\rho_0=8.582\cdot
10^{-6}$kg$\cdot$m$^{-3}$ and get a Knudsen number based on the length
of a vane equal to 0.1.

Converged results are obtained with a $30^2$ points velocity grid and
$400^2$ cells spatial mesh.  This large number of spatial cells is
required to describe the shape of the vanes with enough accuracy.  We
plot in figure~\ref{Radiometer_2} the radial velocity of the vanes as
a function of time and compare our results to those obtained
in~\cite{Chen-Xu-Lee-2012}.  The results are in very good agreement.

\subsection{Roots blower}
\begin{figure}[h]\begin{center}
	\subfigure[Geometry of a lobe. The plain line is an epicycloidal and 
			the dashed line is an hypocycloidal.]
		{\resizebox{0.4\textwidth}{!}{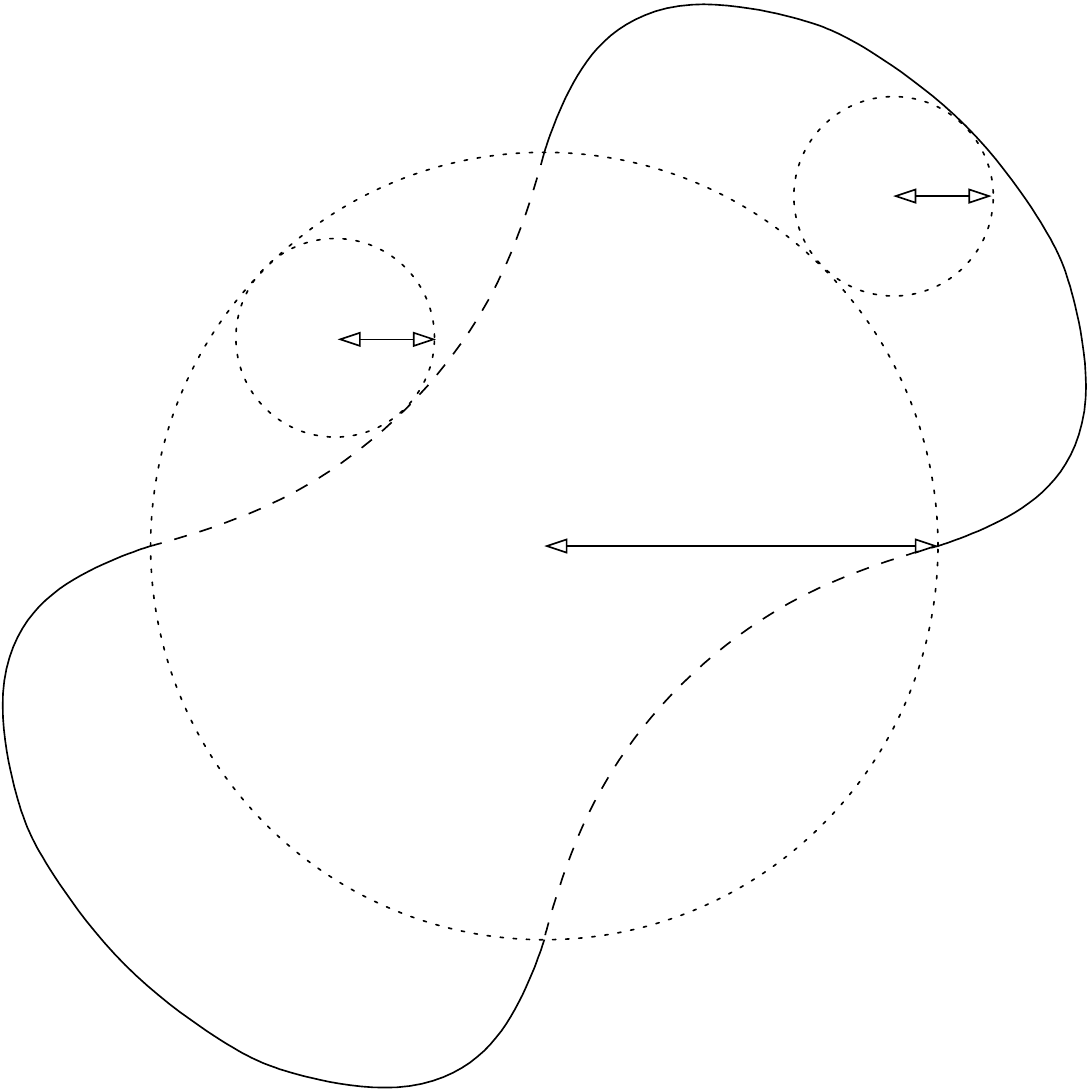}\label{Pump_1}}
	\hspace*{0.15\textwidth}
	\subfigure[Pump geometry. All the units are given in
        meter. The
        computational domain is the square ${[-0.4, 0.4]}^2$.]
		{\resizebox{0.4\textwidth}{!}{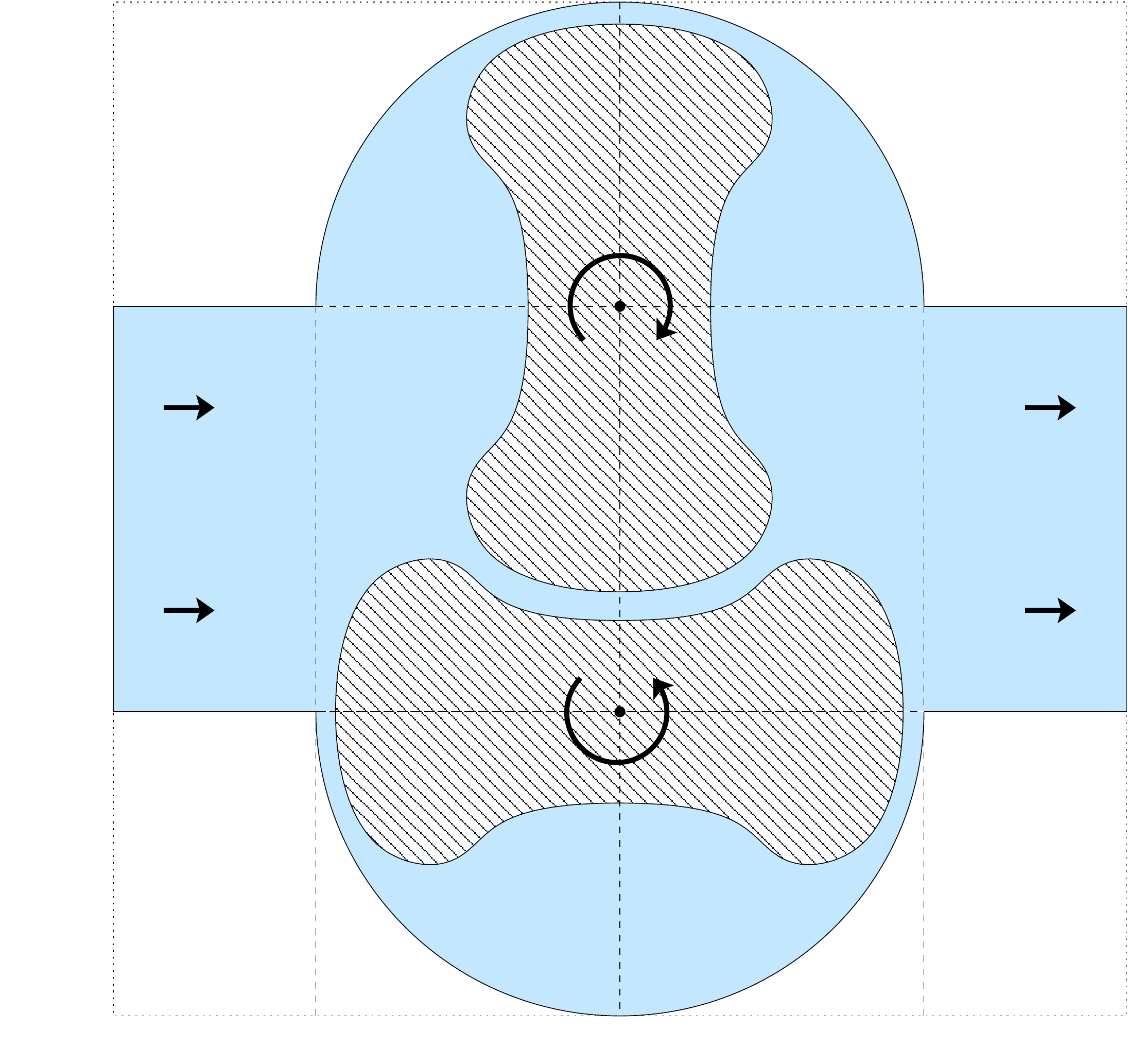}\label{Pump_2}}
	\caption{Roots blower.}
\end{center}\end{figure}

Various kinds of vacuum pumps are used in industrial processes. A common one is the
Roots blower.  It is made of several lobes that rotate simultaneously.
As a result, the gas is trapped by the lobes at one side and then
carried to the other side of the pump.  The simplest shape of Roots
blower is a two-lobed rotor.  In this case the profile of a lobe is
defined by sections of epicycloid and hypocycloid (figure~\ref{Pump_1}).  In
parametric coordinates, this profile is given for all $\theta\in
[-\pi,-\frac{\pi}{2}]\cup[0,\frac{\pi}{2}]$ by the epicycloidal
equation
\begin{displaymath}
	\left\lbrace \begin{aligned}
		&x(\theta) = 5r\cos(\theta)-r\sin(5\theta),\\
		&y(\theta) = 5r\sin(\theta)-r\sin(5\theta),		
	\end{aligned} \right.
\end{displaymath}	
and for all $\theta\in [-\frac{\pi}{2},0]\cup[\frac{\pi}{2},\pi]$ by the hypocycloidal equation
\begin{displaymath}
	\left\lbrace \begin{aligned}
		&x(\theta) = 3r\cos(\theta)+r\sin(3\theta),\\
		&y(\theta) = 3r\sin(\theta)-r\sin(3\theta),		
	\end{aligned}\right.
\end{displaymath}
where $r$ is the radius of the small generating circle that rolls on
the large circle of radius $4r$.  In our simulation, we took $r=3.8$cm
for both lobes.  The full geometry of the pump is detailed in
figure~\ref{Pump_2}.  Note that lobes are not in contact : the minimum
distance between them is $d=1.6$cm.

Since this type of pump mostly operates in atmospheric environment,
the initial conditions are given by $T_0=300$K, $P_0=10^5$Pa,
$\vec{u}_0=\vec{0}$, and the considered gas is argon.  The relaxation
time is computed with formula~(\ref{RelaxationTime}) where the
viscosity coefficient and index for argon are provided by
Bird~\cite{Bird-1994}, that is $\mu=2.117\times10^{-27}$Nsm$^{-2}$ and
$\omega=0.81$.  Because there is no friction between the lobes, a
Roots blower can proceed at a rotary speed that range from 1500rpm to
3000rpm ($\sim 150$rad$\cdot$s$^{-1}$ to $300$rad$\cdot$s$^{-1}$).
For the following simulations, the rotational velocity of the lobes is
set to $\dot{\theta}=\pm 200$rad$\cdot$s$^{-1}$. Note that since the
velocity is imposed here, step 5 of the algorithm is not used (see
section~\ref{subsubsec:summary}). In practice, the temperature of the
lobes tends to increase because of the mechanical heating due to their
high rotary speed.  However, to make it simpler, the wall of the Roots
blowers and its lobes are modeled with diffuse boundary conditions
with constant temperature $T_0$.  At the left side of the pump, it is
assumed that all the gas surrounding the computational domain is in
the same state as the gas located at the inlet.  This can be modeled
by a Neumann boundary condition.  At the outlet (right side of the
pump), we assume that the gas is released in the atmosphere, and hence
the boundary condition is given by a Maxwellian built with the initial
conditions $P_0$, $T_0$ and $\vec{u}_0$.

\begin{figure}[h]\begin{center}
	\subfigure[Pressure distribution at $t=0.03s$.]
		{\includegraphics[width=0.48\textwidth]{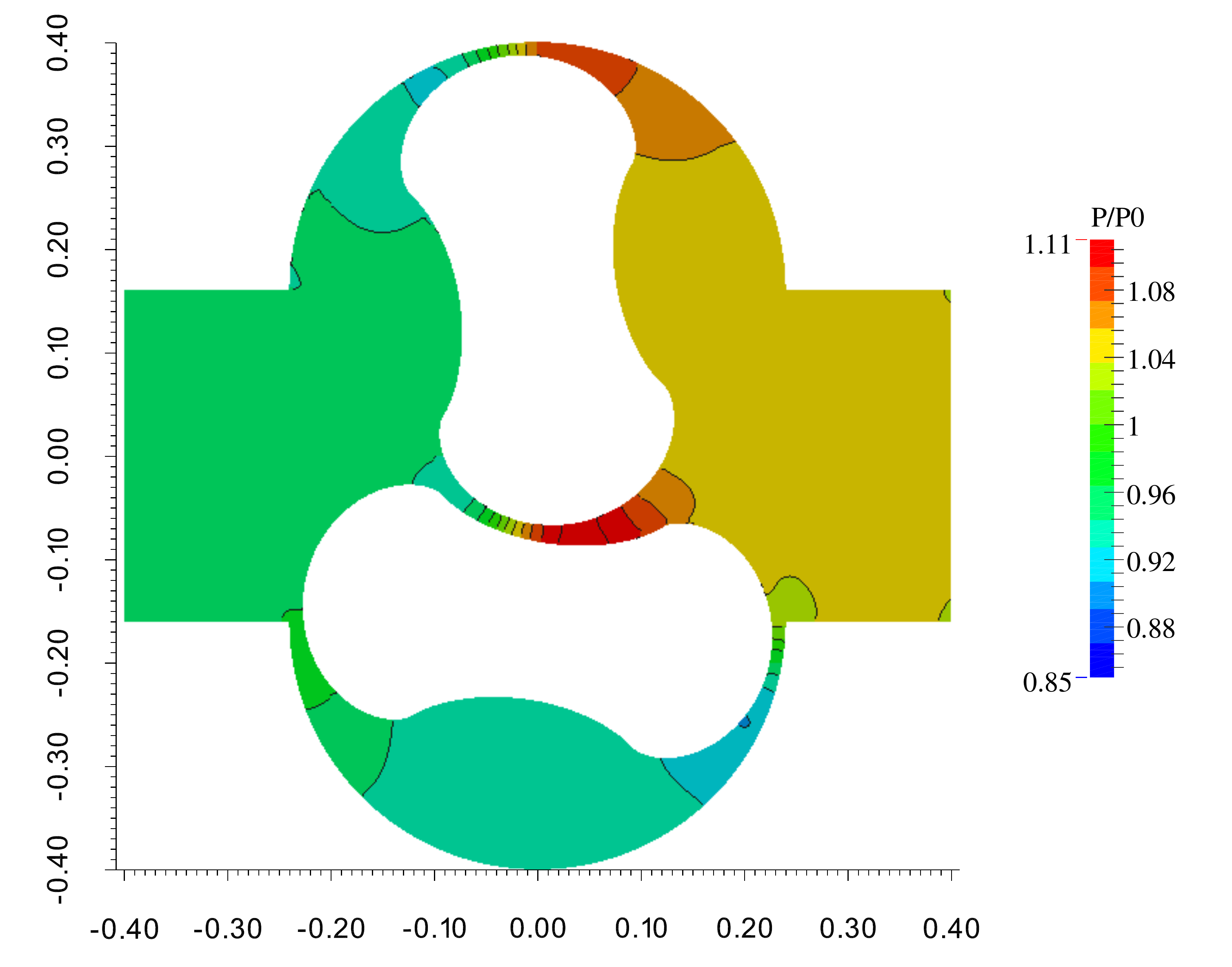}}
	\hspace*{0.02\textwidth}
	\subfigure[Pressure distribution at $t=0.06s$.]
		{\includegraphics[width=0.48\textwidth]{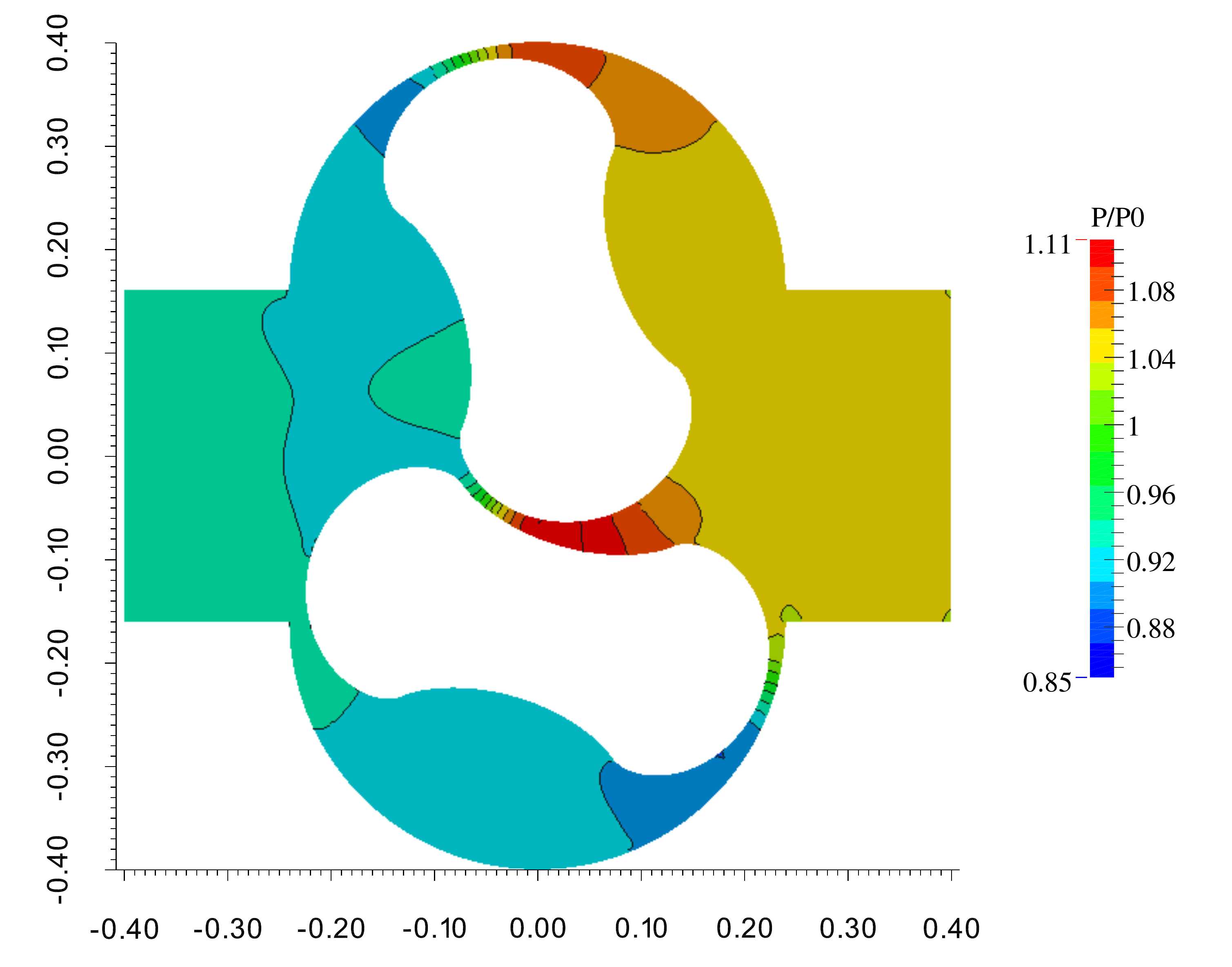}}\\
	\subfigure[Pressure distribution at $t=0.09s$.]
		{\includegraphics[width=0.48\textwidth]{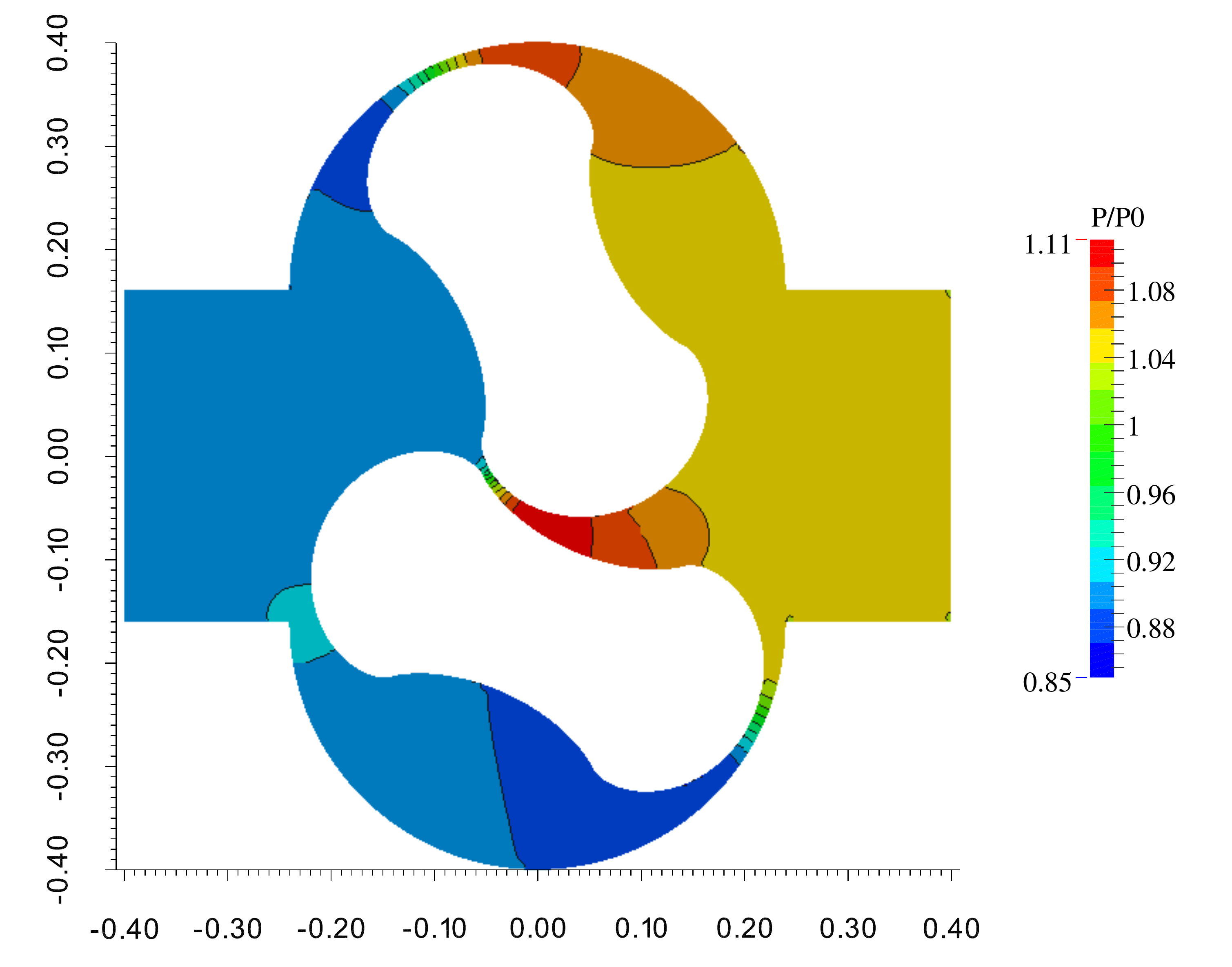}}
	\hspace*{0.02\textwidth}
	\subfigure[Inlet pressure as a function of time:
		$P_{\text{inlet}}=\int_{\text{inlet}} \frac{P(x)}{P_0}\text{d}x$.]
		{\includegraphics[width=0.48\textwidth]{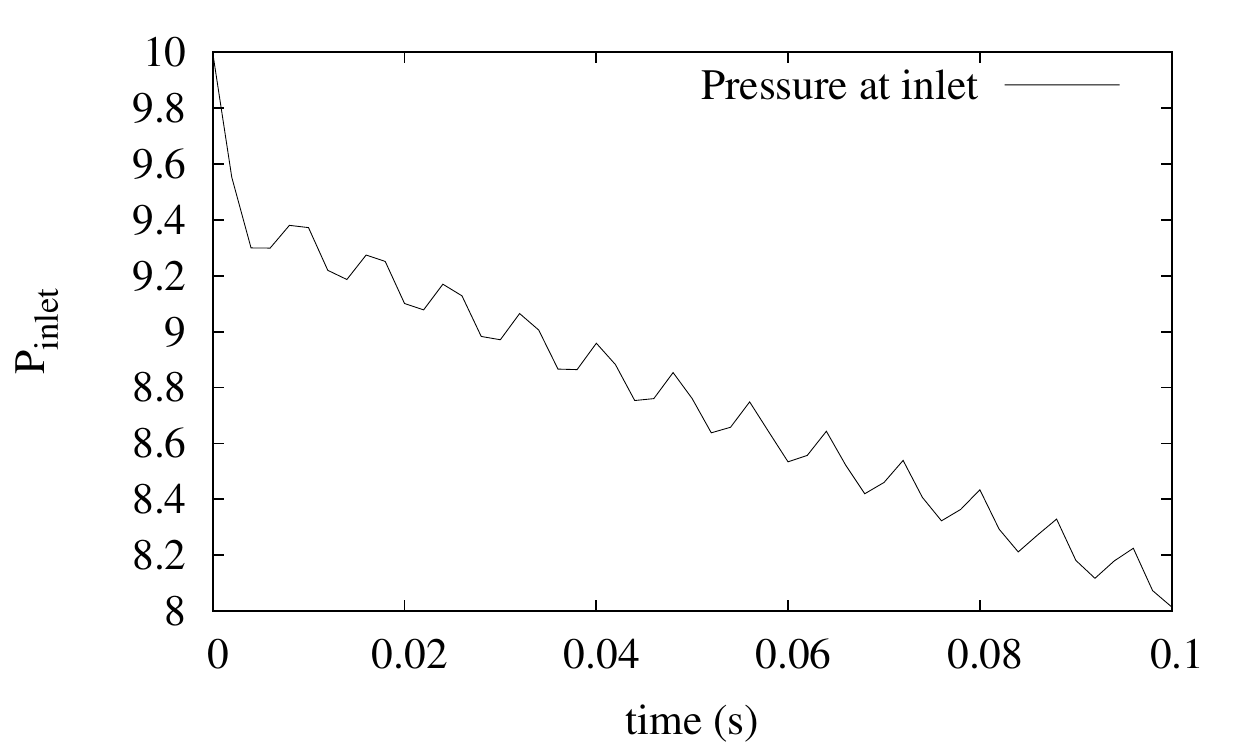}\label{PressureInlet}}
	\caption{Pressure in the pump at several times.\label{PressureProfiles}}
\end{center}\end{figure}

Pressure contours at several times are shown
figure~\ref{PressureProfiles}.  We observe that the pressure at
outlet does not change while the pressure at inlet decreases, as it is
shown in figure~\ref{PressureInlet}.  We stop
the computation at $t=0.1s$.  At this time, the inlet pressure
is $80\%$ of the initial pressure, which means that we
get a pressure drop of $20\%$.

We point out that this simulation is only a qualitative analysis.
Kinetic equations are not really relevant here because the Knudsen
number is very small: $\Kn \approx 3\times 10^{-5}$ for a reference length
equal to the distance $d$ between the two lobes (it would be larger
with a smaller distance).
Hence Navier-Stokes equations might be more relevant in this case.
However, this simulation shows that the cut cell method works well
with complex shaped objects for moderate velocities flows
($Ma\approx0.15$), while this would be more more difficult with the
moving mesh approach, for instance.

%-----------------------------%
%                             %
%   Numerical discretization  %
%                             %
%-----------------------------%
\section{Three dimensional flow simulations}
\label{sec:3D}

In this section, the cut cell method presented in
section~\ref{SecCutCellDD} is extended to 3D simulations. The approach
is only detailed here for the simulation of the Crookes radiometer,
which means that only pure rotation is considered.  The angle between
the position of a vane at time~$t$ and its initial position is denoted
by~$\ThetaSol(t)$, and~$\dotThetaSol(t)$ stand for its rotational
velocity.  The gas governing equations are that of
section~\ref{SubSecBoltzmannTD} and the velocity discretization is
similar to the one explain in section~\ref{SubSecVelocity}: there are
$N^3$ velocity points and the $\iv^{\text{th}}$ velocity is denoted by
$\vTD_{\iv}=(v_{\iv_1},v_{\iv_2},v_{\iv_3})$.
The different steps of the cut cell method -- that consists in
updating the values $\ThetaSol^n$, $\dotThetaSol^n$ and $\fVirtTD^n$
for the next time step -- are that of section~\ref{subsubsec:summary}.
They are written below to highlight the differences with respect to
the 2D case. 

We assume that at time $t^n$ we have the average value of the
distribution function $\fVirtTD^n$ in every virtual cell
$\VirtTD^n$. Like in 2D, a virtual cell is defined as the intersection
between the gaseous domain~$\Omega_g(t^n)$ and the cuboid
cell~$\Omega_{i,j,k}$ of the Cartesian mesh, but now it is a polyhedral
with seven virtual faces.  The first six faces fit with the Cartesian
mesh interfaces and are denoted by $\SxTD{\pm}^n$, $\SyTD{\pm}^n$,
$\SzTD{\pm}^n$.  The seventh one is a plane approximation of the solid
boundary denoted by $\SwTD$ and $\nwTD$ is its normal vector directed
outward.  We also assume that at time $t^n$, the angle $\ThetaSol^n$
and rotational velocity $\dotThetaSol^{n}$ of the solid are known. The
different steps are the followings:

\begin{enumerate}
	\item The position  of the solid body at $t^{n+1}$ is computed with~(\ref{motionDD_posit}).

	\item The virtual cells are arranged into
          control volumes $\CVolTD^{n}(t_{n})$ following the rule
          given in section~\ref{subsubsec:CV} which is naturally
          extended to 3D. The average value $\fCVolTD^n$ of the distribution over
          the control volume $\CVolTD^n(t^n)$ is computed: 2D
          relation~(\ref{FusionDD}) is replaced by
	\begin{displaymath}
		\fCVolTD^n:=\frac{1}{\VCVolTD^n }\sum_{(i',j',k')\in\sigma_{i,j,k}^n}
			\VVirtTDp^n\,\fVirtTDp^n,
	\end{displaymath}

	\item 
The virtual cells and control volumes are
          moved according to the new position computed at step 1.  The
          volumes $\VVirtTD^n$, $\VCVolTD^n$, $\VCVolTD^{n+1,\ast}$ of
          the virtual cells $\VirtTD^n$, control volumes
          $\CVolTD^n(t^n)$ and $\CVolTD^n(t^{n+1})$ at times $t^n$ and
          $t^{n+1}$, respectively, are computed.

	\item The boundary condition is computed for each cut cell with the discrete form of~(\ref{BoundaryTD}), which yields:
	\begin{displaymath}
		\fTD_w(t^n,\xwTD\in\Bound,\vTD_{\iv}\in\Vin) = \rhow\MaxwTD[1,\uwTD,\Tw],
	\end{displaymath}
	where $\rhow$ is given by
	\begin{displaymath}
		\rhow = -\frac{\sum_{\vTD_{\iv}\in\Vout}(\vTD_{\iv}-\uwTD(t^n,\xwTD))\cdot\nwTD(t^n,\xwTD)
		\fCVolTD^n\Delta v_x\Delta v_y\Delta v_z}
{		\sum_{\vTD_{\iv}\in\Vin}(\vTD_{\iv}-\uwTD(t^n,\xwTD))\cdot\nwTD(t^n,\xwTD)
		\MaxwTD[1,\uwTD(t^n,\xwTD),\Tw]\Delta v_x\Delta v_y\Delta v_z}
	\end{displaymath}
	Here, $\xwTD$ is the center of $\SwTD$ and
        $\Vin=\{\vTD_{\iv}|(\vTD_{\iv}-\uwTD(t^n,\xwTD))\cdot\nwTD(t^n,\xwTD)<0\}$ and
        $\Vout=\{\vTD_{\iv}|(\vTD_{\iv}-\uwTD(t^n,\xwTD))\cdot\nwTD(t^n,\xwTD)>0\}$ are the
        sets of incoming and outgoing velocities, respectively.  
We introduce the cylindrical coordinates of
        $\xwTD=(r^n\cos\alpha^n,r^n\,\sin\alpha^n,z^n)$ in order to write the
        boundary velocity as $\uwTD(t^n,\xwTD) = r^n\dotThetaSol^n\times
        (\sin\alpha^n,\cos\alpha^n,0)$.

	\item The stress tensor is computed by 3D extension of~(\ref{StressTensorDD}):
	\begin{displaymath}\begin{split}
	\pwTD(t^n,\xwTD)&=\sum_{\vTD_{\iv}\in\Vin}(\vTD_{\iv}-\uwTD(t^n,\xwTD))\otimes
						(\vTD_{\iv}-\uwTD(t^n,\xwTD))\fBoundTD \,\Delta v_x\Delta v_y\Delta v_z\\
			&\quad +\sum_{\vTD_{\iv}\in\Vout}(\vTD_{\iv}-\uwTD(t^n,\xwTD))\otimes
						(\vTD_{\iv}-\uwTD(t^n,\xwTD))\fCVolTD^n \,\Delta v_x\Delta v_y\Delta v_z.
	\end{split}\end{displaymath}
	The rotational velocity is then computed with 
	\begin{displaymath}
	\dotThetaSol^{n+1} = \dotThetaSol^n	+  \frac{\Delta t}{J}
		\sum_{i=1}^{N_x}\sum_{j=1}^{N_y}\sum_{k=1}^{N_z}
			\left[\xwTD\times \big(\pwTD\cdot\nwTD(t^n,\xwTD) \big)\DwTD\right]
			\cdot \left[\begin{array}{c}0\\0\\1\end{array}\right],
	\end{displaymath}
	where the sum is an approximation of the torque acting on the
        vanes. Then the new wall velocity at time $t^{n+1}$ is computed.

      \item Scheme~(\ref{ExSchemeDD}) is easily extended to 3D to
        compute the average value $\fCVolDD^{n+1,\ast}$ of $f$ at time
        $t^{n+1}$ in the control volume $\CVolTD^n(t^{n+1})$: indeed,
        we write the integral form of Boltzmann
        equation~(\ref{IntBoltzmannDD}) with $V(t)=\CVolTD^n(t)$ and
        simplify the transport term just like in~(\ref{TransportDD})
        to get the first order explicit scheme
\begin{align}
		\fCVolTD^{n+1,\ast} = \frac{\VCVolTD^n}{\VCVolTD^{n+1,\ast}}\fCVolTD^n 
			&-\frac{\Delta t}{\VCVolTD^{n+1,\ast}}
			\sum_{(i',j',k')\in\sigma_{i,j,k}^n}\Big[
			\left(\FluxiTDp{+}^n-\FluxiTDp{-}^n\right)
                        \nonumber \\
			&\hspace*{113pt}
                        +\left(\FluxjTDp{+}^n-\FluxjTDp{-}^n\right)
                        \nonumber \\
			&\hspace*{113pt}	+\left(\FluxkTDp{+}^n-\FluxkTDp{-}^n	\right)
			+ \FluxwTDp \Big] \nonumber \\
			&+ \frac{\VCVolTD^n}{\VCVolTD^{n+1,\ast}}
			\frac{1}{\tau_{i,j,k,\iv}^n} \big(
                        \feqTD_{i,j,k,\iv}^n - \fCVolTD
                        \big), \label{eq-schema3D}
\end{align}
%\end{displaymath}
	where the upwind numerical flux $\FluxiTD{+}^n$, $\FluxjTD{+}^n$, $\FluxkTD{+}^n$ and $\FluxwTD^n$ are:
	\begin{displaymath}\begin{aligned}
	&\FluxiTD{+}^n := \DxTD{+} &&
		\big[\min(v_{\iv_1},0)\,\fTD_{i+1,j,k,\iv}^n+\max(v_{\iv_1},0)\,\fCVolTD^n\big]\\
	&\FluxjTD{+}^n := \DyTD{+} &&
		\big[\min(v_{\iv_1},0)\,\fTD_{i,j+1,k,\iv}^n+\max(v_{\iv_1},0)\,\fCVolTD^n\big]\\
	&\FluxkTD{+}^n := \DzTD{+} &&
		\big[\min(v_{\iv_1},0)\,\fTD_{i,j,k+1,\iv}^n+\max(v_{\iv_1},0)\,\fCVolTD^n\big]\\
	&\FluxwTD^n := \DwTD &&
		\big[\min([\vTD_{\iv}-\uwTD(t^n,\xwTD)]\cdot\nwTD(t^n,\xwTD),0)\,\fBoundTD \\
& & & 		+\max([\vTD_{\iv}-\uwTD(t^n,\xwTD)]\cdot\nwTD(t^n,\xwTD),0)\,\fCVolTD^n\big]
	\end{aligned}\end{displaymath}
    \item This average value of the distribution function is
      distributed to the virtual cells~$\VirtTD^{n+1}$ merged into
      $\CVolTD^n(t^{n+1})$ by $\fVirtTD^{n+1} := \fCVolTD^{n+1,\ast}$.
\end{enumerate}
\subsection{Octree procedure}
Full 3D simulations are computationally very expensive.  For instance,
based on the 2D computations presented in
section~\ref{subsec:crookes}, a 3D simulation of the Crookes
radiometer requires $200^3$ degrees of freedom for the space
discretization and $30^3$ for the velocity discretization.  In
addition to these \numprint{2 160}$\times 10^8$ grid points,
\numprint{10000}~time steps are expected to reach the stationary
rotational velocity of the vanes.  In conclusion, even a massively
parallel computing is not sufficient to do to this simulation in a
reasonable computational time (i.e. less than several weeks).

In order to reduce the computational time, an octree procedure is implemented.
First, a coarse Cartesian mesh is initialized and refined around boundaries.
This means that all the cells that are close enough to the boundary are divided in 8 smaller cells (or 4 in 2D).
These new cells can be divided again if necessary; and we define the
depth of a cell as the number of divisions that lead to this cell. The
splitting criterion is the following:
\begin{equation}\label{criteria}
   \text{if }\quad \phi<2\sqrt{2}\frac{\Delta x}{2^d}\quad\text{ and }\quad d\leq d_{\max}\quad\text{the cell is divided},
\end{equation}
where $\phi$ is the distance from the cell center to the nearest
boundary and $d_{\max}$ is a prescribed maximum depth.  For the 3D
simulation of the Crookes radiometer, the coarse mesh is made of
$50^3$ cells and $d_{\max}=2$.  With these parameters, the spatial
mesh only contains~\numprint{240000} cells.
%The smallest cells are therfore of the same size as the ones of a Cartesien $200^3$ mesh.
Since the boundary moves, the mesh is adapted to the new location
of the boundary at the beginning of each time iteration.  Note that it
may happen that 8 cells have to merge during this process, if
all of them no longer respect criterion~(\ref{criteria}).

When the mesh changes -- i.e when cells split or merge -- the
distribution function has to be interpolated on the new mesh.  This is
done by assuming that the distribution function is constant over a
cell and by using a standard restriction/prolongation method (by
average and 0th order interpolation).  Since every cells are cuboid,
scheme~(\ref{eq-schema3D}) can be applied by using an
appropriate data structure to access the neighboring cells of each
numerical interfaces: here, we use the standard Z-ordering which is very
efficient for that.

\subsection{Parallel implementation}
We describe here two natural strategies for a parallel
implementation of a kinetic solver with the \textit{Message Passing
  Interface} (MPI) library, and we propose our own hybrid technique.

The first method is velocity parallelization, or decomposition domain
method in the velocity space: each processor computes the distribution
function in the whole space domain, but for only a part of the
discrete velocities. This approach is for instance used
in~\cite{Aristov-Zabelok-2002}  and more
recently in~\cite{Titarev_2012}.  For a given discrete velocity, the
scheme is independent of the other velocities: then it is is used
independently by each processor, and each of them compute partial
moments (by using its own reduced set of discrete velocities). To
compute the full moments, all processors gather the sum of the partial
moments.

The second method is a more standard space domain decomposition,
see~\cite{Ilgaz-Tuncer-2007,Kolobov-2007}.  Each processor computes
the distribution function for the whole discrete velocity grid, but
only for a subdomain in the position space. To compute the numerical
fluxes across the interfaces between different subdomains, the method
requires communications: it is sufficient that each processor sends
the distribution of its interface cells to
its neighbors. We refer to~\cite{TDU2014} for a comparison of these
two strategies.

Our implementation combines these ideas. Each
processor uses the scheme on a space subdomain, for a partial set of
discrete velocities. The space domain decomposition is made on the
initial mesh, before the refinement procedure: the Cartesian structure
of this mesh makes the portioning very easy. For problems with moving
boundaries, it is difficult to ensure a good dynamic load balancing
between different processors: the number of cells of a subdomain can
change a lot due to the space refinement induced by the displacement
of the solid obstacle. In order to optimize the workload distribution,
we use a small number of subdomains. Groups of processors are
given to each subdomain, and each processor will apply the scheme for
a partial set of discrete velocities. The advantage of this technique
is that we can use a large number of processors without a two large
number of space subdomains. This makes the workload well balanced 
during the simulation for each subdomain, at least for the Crookes
radiometer presented in the following section, since there is always
one vane in each subdomain.  
For the corresponding 3D simulation, our technique is quite efficient,
since it has been made with 240 processors for a
CPU time lower than 12 hours.

Note that there is an other kind of hybrid parallelization which uses
both MPI and OpenMP libraries (see~\cite{Baranger-2014}), but this is
not what is used here.

\subsection{Numerical example : the Crookes radiometer}
\begin{figure}[!h]\begin{center}
	\raisebox{25pt}{
		\includegraphics[width=0.48\textwidth]{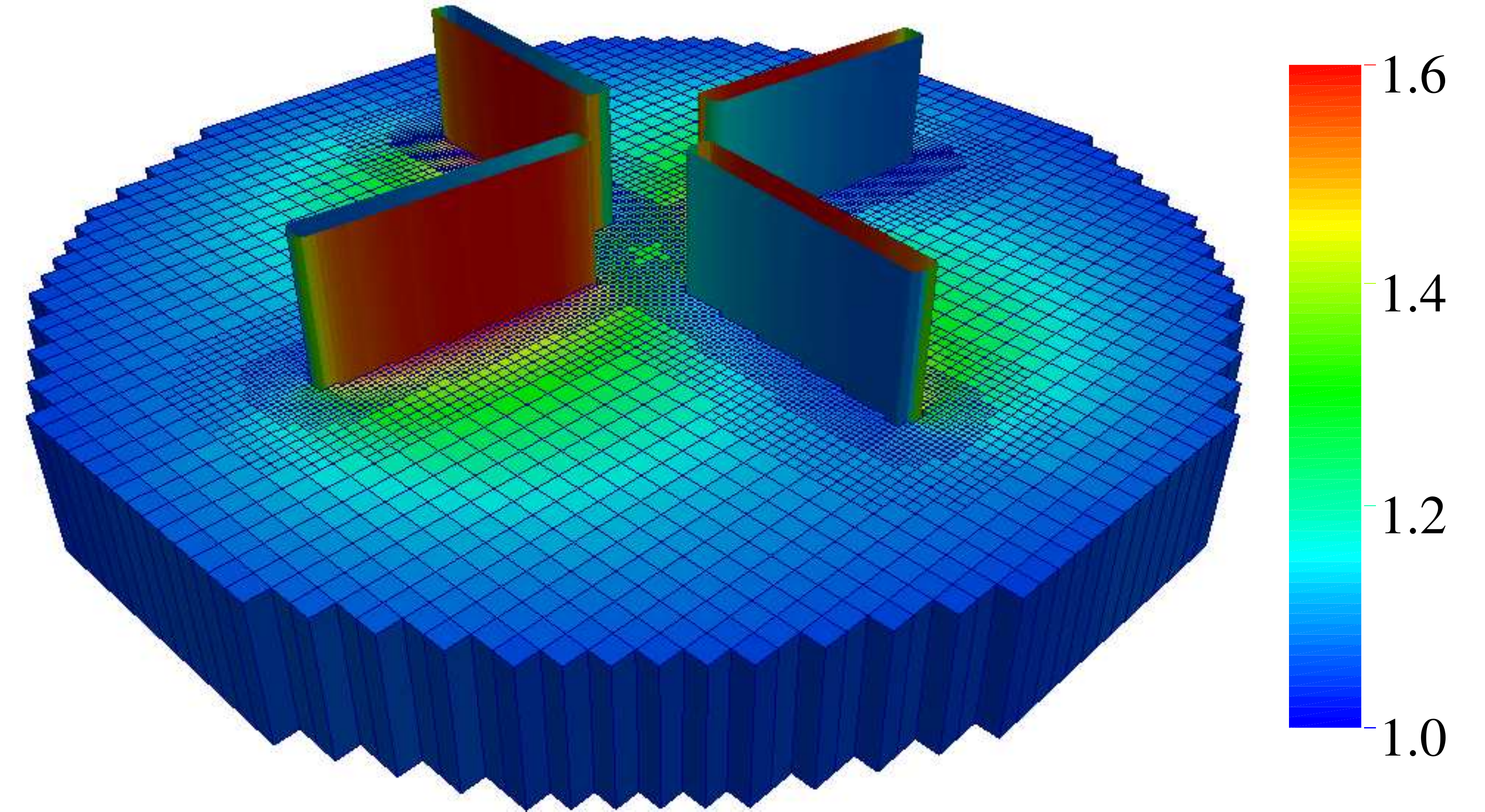}}
	\hspace{0.015\textwidth}
		\includegraphics[width=0.48\textwidth]{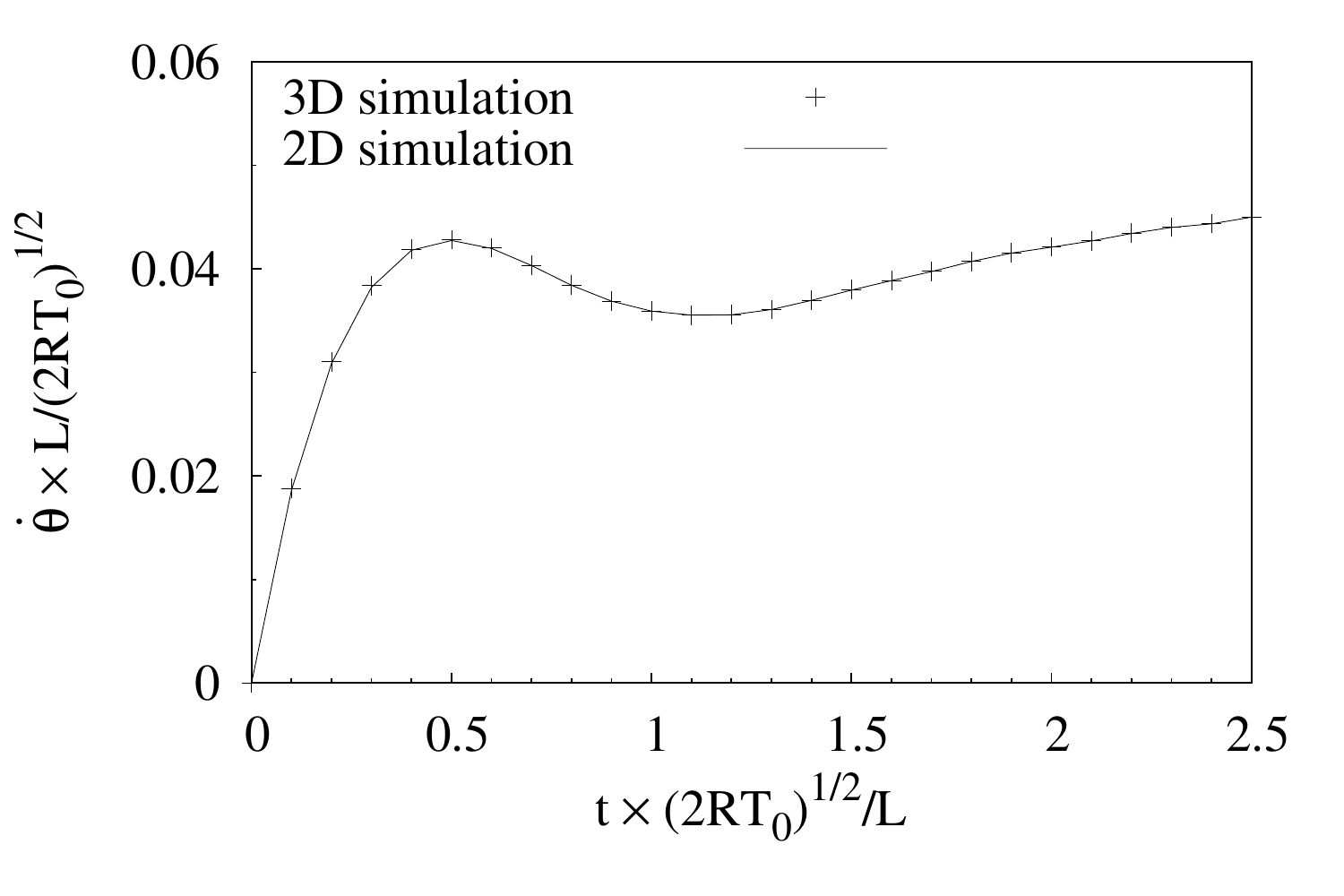}
	\caption{3D extruded radiometer: temperature field $T/T_0$ and mesh
          (left) and radial velocity profile as a function of time for 2D and 3D simulations
          (right).}
	\label{validation}
\end{center}\end{figure}

First, the implementation of the method has been checked with a 3D
simulation of the plane 2D radiometer similar to the one presented in
section~\ref{subsec:crookes}: this 2D geometry is extruded to get a
cylinder shaped radiometer (see figure~\ref{validation}, left), and
periodic boundary conditions are imposed at the upper and lower
boundaries to simulate the infinite vanes. In this case, the moment of
inertia is $97\rho_0 L^5/370$, where $L$
is the height of the extruded vanes. For a Knudsen number of 0.5, 2D
and 3D simulations give exactly the same results, as it is shown in
figure~\ref{validation}, right. For this comparison, the AMR
technique is used for both simulations, and a plane section of the
3D mesh is the same as the 2D mesh.

\begin{figure}[!h]\begin{center}
	\resizebox{0.9\textwidth}{!}{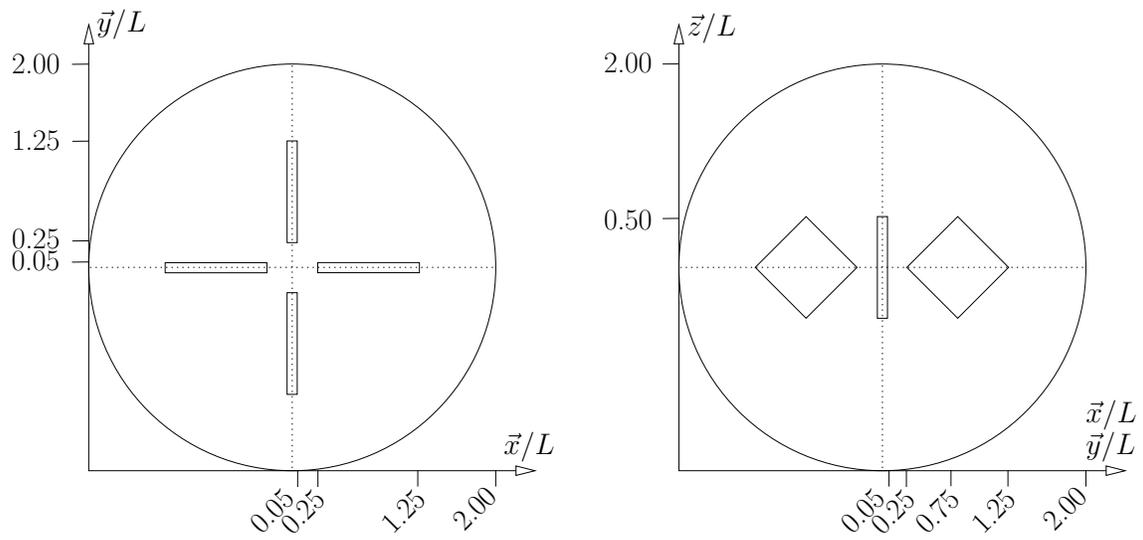}
	\caption{3D radiometer: cross section in the
          plane $xOy$ (left) and in the plane $xOz$ (right). See
          figure~\ref{radiometre3D-simu} for a 3D view.}
	\label{radiometre3D-geometrie}
\end{center}\end{figure}

Now the real 3D radiometer is made of four square shaped vanes. The
dimensions of the vanes are $L$ for the diagonal and $L/10$ for the
width. They are immersed in a sphere of radius $2L$, and their centers
are in the plane $z=0$, at a distance $0.75 L$ from the center of the
sphere. The corresponding geometry is shown in
figure~\ref{radiometre3D-geometrie}.  The moment of inertia of the
vanes is computed with a material density $\rho_0$ equal to the mass
density of the surrounding gas which gives $J= \frac{97}{750}\, \rho_0 L^5$. This
density is not realistic, but it makes the vanes faster, and it
is easier to observe their movement.

At $t=0$, the radial velocity of the vanes is zero and the temperature
$T_0$ in the domain is uniform.  The gas is governed by the BGK model,
where the relaxation time is $\tau =
\frac{\sqrt{\pi}}{2}\frac{\rho_0}{\rho}\frac{\Kn L}{\sqrt{2RT_0}}$. All
the boundary conditions are diffuse reflections with constant
temperatures: the sphere boundary is maintained at temperature $T_0$
and the white side of the vanes as well, while their black side is
maintained at temperature $2T_0$. On the edge of the vanes, the
temperature is discontinuous ($T_0$ on one part and $2T_0$ on the
other part).

\begin{figure}[h!]\begin{center}
	\includegraphics[width=0.48\textwidth]{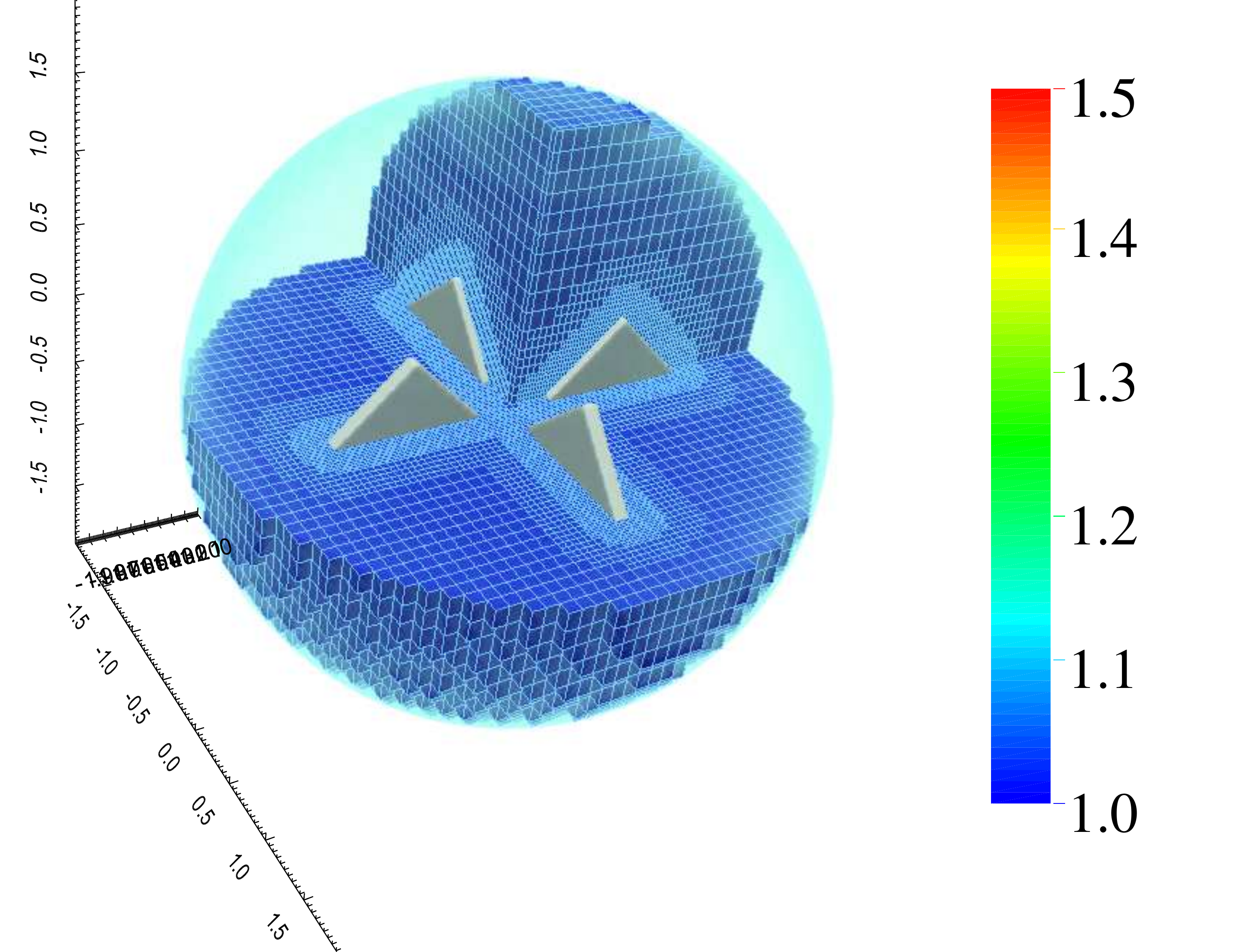}
	\hspace{0.015\textwidth}
	\includegraphics[width=0.48\textwidth]{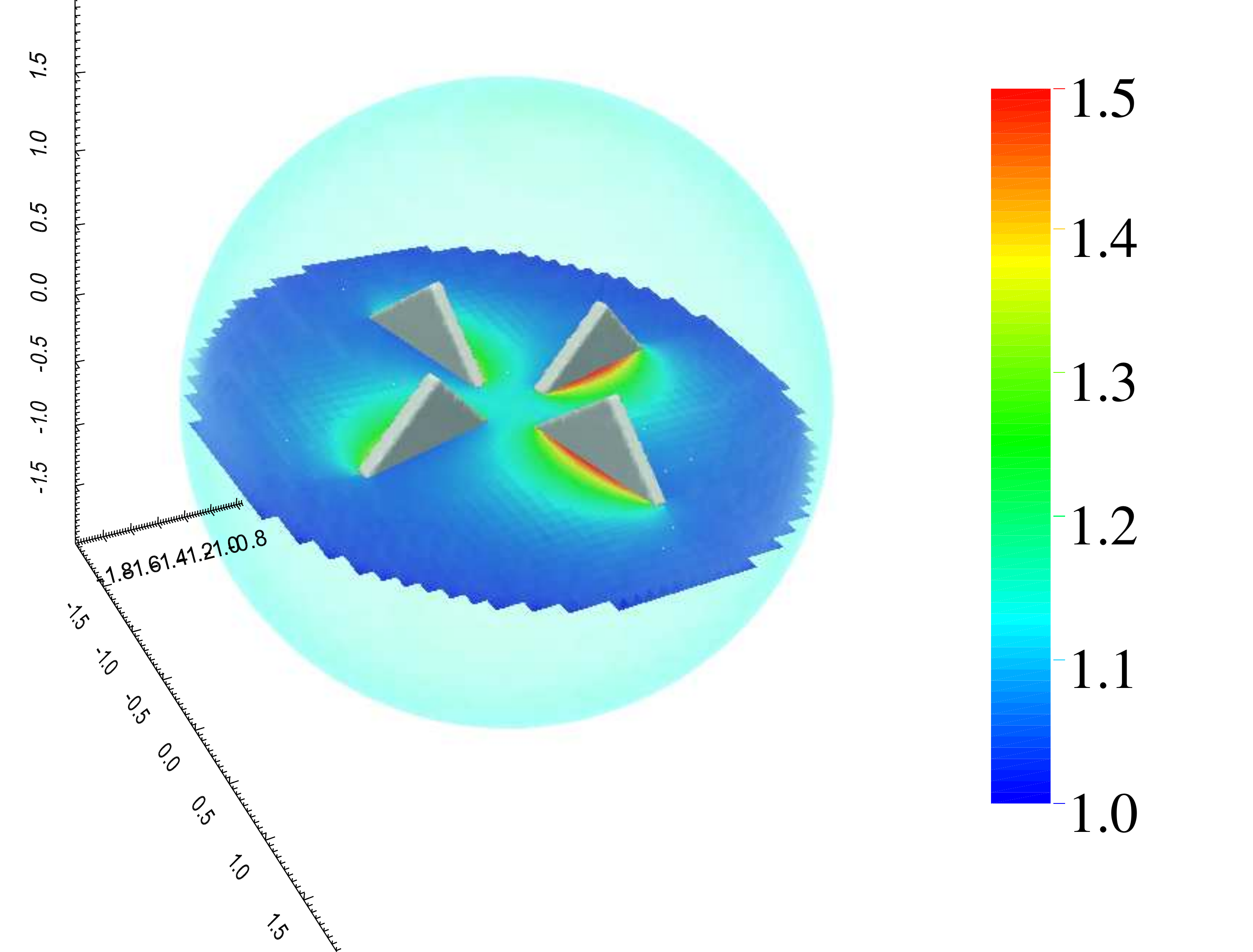}
	\includegraphics[width=0.48\textwidth]{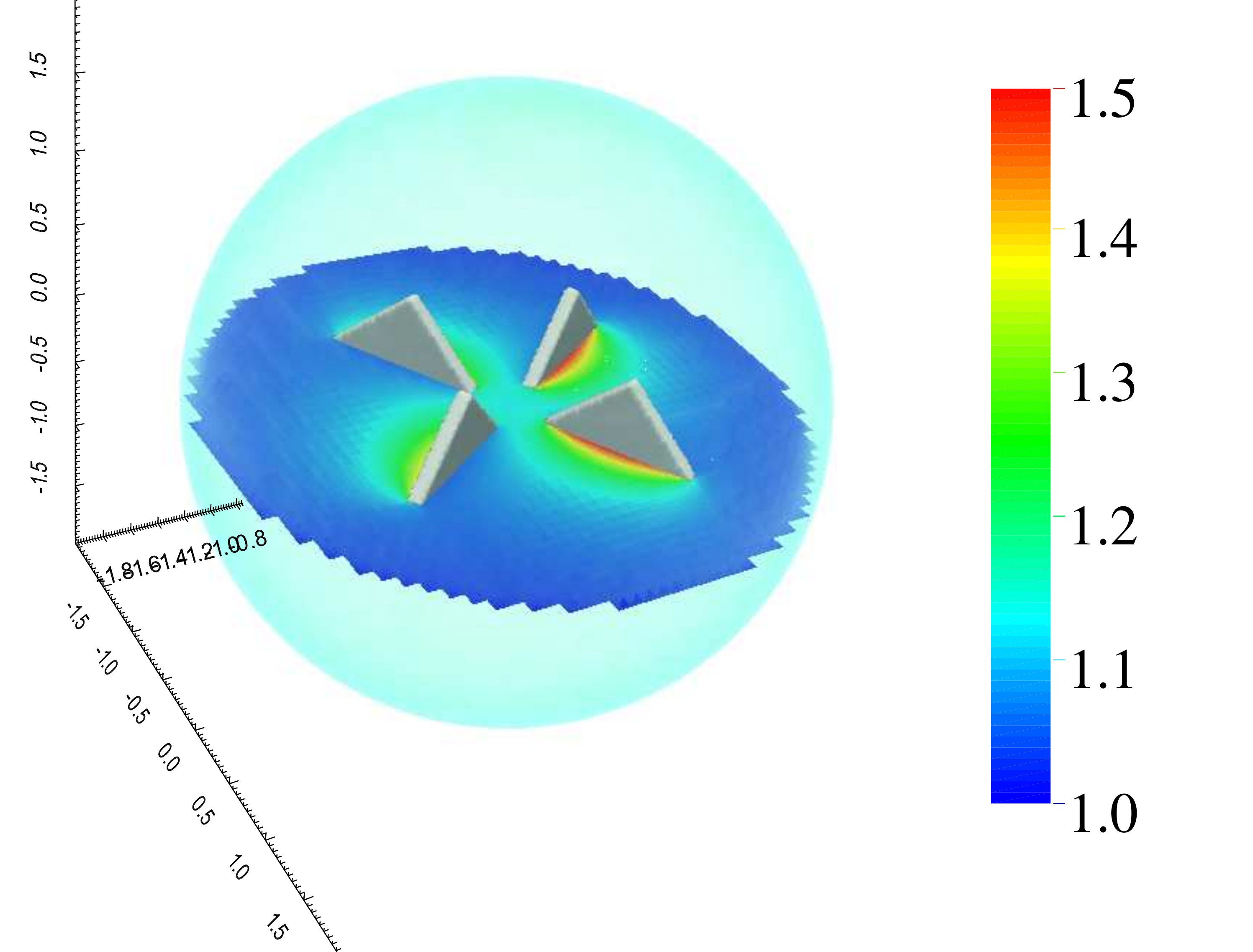}
	\hspace{0.015\textwidth}
	\includegraphics[width=0.48\textwidth]{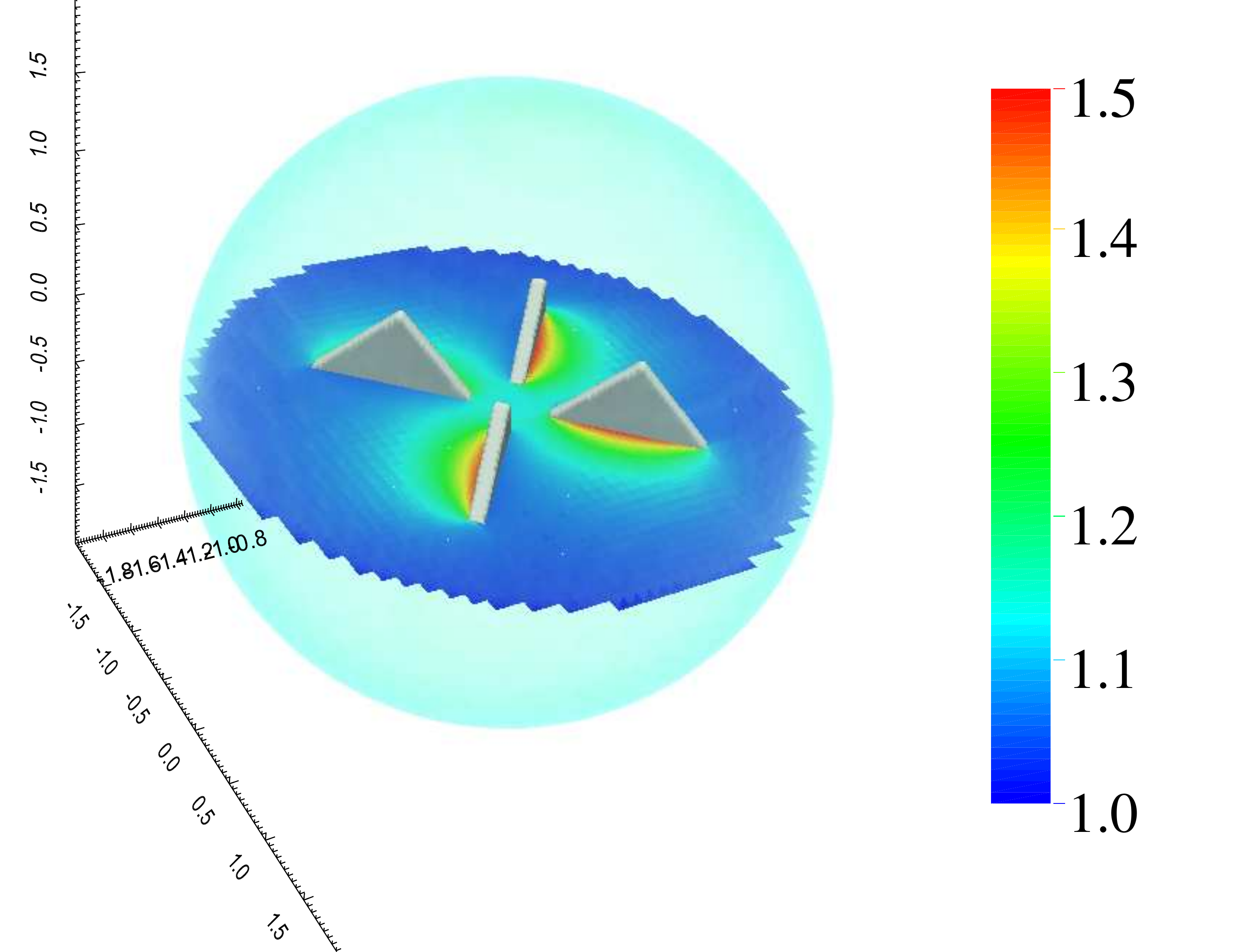}
	\caption{3D Crookes radiometer for $\Kn=0.5$: 
	the temperature field $T/T_0$ on the plane $z=0$ is shown at
        times $t\times L/\sqrt{2RT_0}=0,$
	 5, 10 and 15 (from left to right and from top to bottom). The
         mesh is shown at $t=0$ (top-left).}
	\label{radiometre3D-simu}
\end{center}\end{figure}

We also use the 2D simulation to estimate the resolution required by
the 3D computation. The difference between the results obtained with a
2D mesh of $50^2$ cells refined by the AMR technique with a maximum
depth $d_{\max}=2$ and the fine Cartesian structured mesh of $500^2$
cells is less than 5\%, which is considered as sufficiently accurate
here. Consequently, a 3D AMR mesh of $50^3$ cells with a depth of 2
should be accurate enough for a 3D simulation. This is computationally
possible, since this mesh contains \numprint{240000} cells, which is
much smaller than the equivalent Cartesian mesh of
$200^3=\numprint{8000000}$ cells. The corresponding simulation is
shown in figure~\ref{radiometre3D-simu} at different times.

\begin{figure}[h!]\begin{center}
	\subfigure[Time evolution of the radial velocities.]{
	\raisebox{2pt}{
		\includegraphics[width=0.45\textwidth]{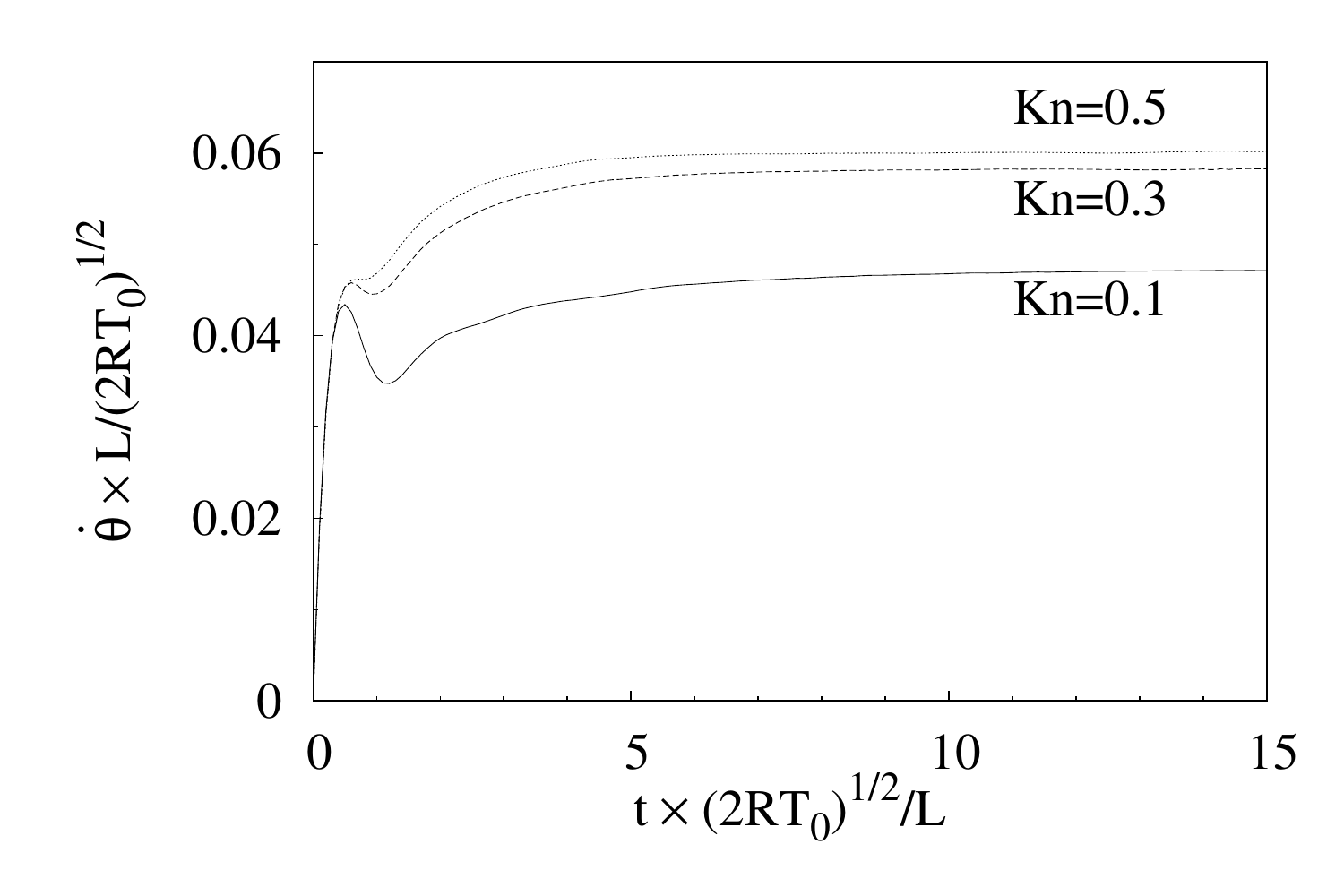}}\label{fig:3Dtheta_t}}
	\hspace{0.01\textwidth}
	\subfigure[Stationnary velocity for 2D and 3D simulations for
        three Knudsen numbers.]{	
	\includegraphics[width=0.45\textwidth]{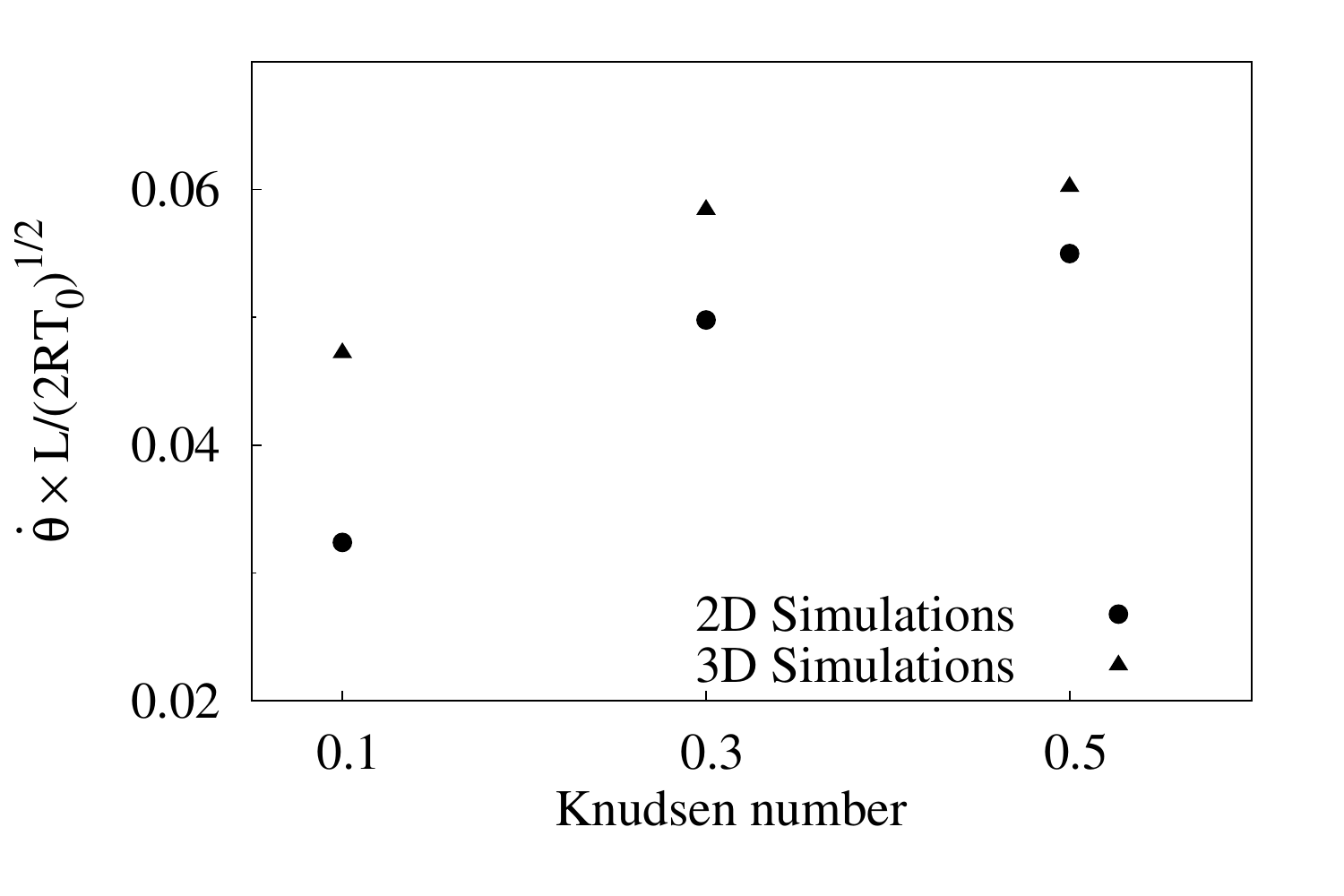}\label{fig:3Dtheta_stat}}
	\caption{Radial velocity profiles for 3D simulations with 3
          different Knudsen numbers.}	
\end{center}\end{figure}

We have made other simulations with three Knudsen numbers 0.1,
0.3, and 0.5 until the steady state is reached. In
figure~\ref{fig:3Dtheta_t} is shown the evolution of the radial
velocity: we clearly see the convergence to a constant velocity, which
is larger for $\Kn=0.5$. A comparison between 2D an 3D geometries is
shown in figure~\ref{fig:3Dtheta_stat}: the 3D vanes are clearly
faster that the 2D vanes.

Note that this test is just shown to illustrate the potential of our
method. We are not aware of any similar simulation in the literature
so far, and we are not able to present any comparison. Moreover, we do
not claim this is a realistic simulation, since the moment of inertia
of the vanes is too small, and their width is too large. However, we
do not know any experimental measures of the motion of the
radiometer. If any, it would probably be necessary to make a more
intensive simulation, since the refinement should be stronger around
the vanes that are generally very thin.

\section{Conclusion}
\label{sec:concl}
A numerical method for solving kinetic equations with moving obstacles
has been presented.  This method is an extension to the kinetic theory
of the cut cell technique used in computational fluid dynamics.  The
main advantage of this algorithm is that it combines the simplicity of
the Cartesian grid based methods to the accuracy of the body fitted
methods, which ensures exact mass conservation.
% The accuracy of the method has been proved by many comparisons with
% results coming from the literature.
The method is easily extended to 3D flows, and its accuracy has been
proved with the simulation of a Crookes radiometer.  Another advantage
of our approach is a simple and generic treatment of all kinds of cut
cells. This is essential, especially for 3D problems in which there
are many different kinds of cut cells.

Our goal is now to improve the accuracy of our method by using a
second order scheme. Since the mesh is Cartesian, the main difficulty
is to approximate the gradient of the distribution function on the cut
cells with enough accuracy.  Such an extension has already been done
(see~\cite{these_Guillaume}), but while it works well for non moving
obstacles, it is not efficient enough for general problems. An other
perspective is the validation of the 3D algorithm, in particular by
using experimental data.

\paragraph{Acknowledgments.} Experiments presented in this paper were
carried out using the Pla\-FRIM experimental testbed, being developed
under the Inria PlaFRIM development action with support from LABRI and
IMB and other entities: Conseil R\'egional d'Aquitaine, FeDER,
Universit\'e de Bordeaux and CNRS (see {\tt
  https://plafrim.bordeaux.inria.fr/}). Computer time for this study
was also provided by the computing facilities MCIA (M\'esocentre de Calcul
Intensif Aquitain) of the Universit\'e de Bordeaux and of the Universit\'e
de Pau et des Pays de l'Adour. This study has been carried
out in the frame of “the Investments for the future” Programme IdEx
Bordeaux – CPU (ANR-10-IDEX-03-02).

%-----------------------------%
%                             %
%          APPENDIX           %
%                              %
%-----------------------------%
\newpage
\appendix

\section{Computation of the cell geometric parameters}
\label{LevelSet}

We describe below how these parameters are computed for 2D
problems. There is no specific difficulty to extend our algorithms to
3D problems: however, the formula are a bit long, and to
shorten the paper, this extension is left to the reader.

The computation is made with the following four steps: 
\begin{enumerate}
\item identification of the type of each cell (solid, gas, our cut cell);
\item for each virtual cell, computation of the lenghts of its five
  edges, and computation of its normal vector of its fifth edge
  (the one which is a linear approximation of the solid boundary);
\item identification of the cut cells that have to merge, and
  computation of the set $\sigma_{i,j}$;
\item computation of the area of each virtual cells and corresponding
  control volumes.
\end{enumerate}

Before we describe these steps in the following sections, we point out
that a level-set signed distance function $\Lvset :
\xDD\rightarrow\Lvset(\xDD)$ is systematically used. It gives the
shortest distance between a point $\xDD$ in the computational domain
and the solid boundary. This function is negative if $\xDD$ is inside
the solid, and positive if it is inside the gas, and hence the zero
level-set of $\Lvset$ is the solid boundary. In our algorithm, the
values of $\Lvset$ are computed analytically at each node
$\xDD_{i\pm\frac{1}{2},j\pm\frac{1}{2}}$ of the Cartesian grid, and we
set
$\Lvset_{i\pm\frac{1}{2},j\pm\frac{1}{2}}:=\Lvset(\xDD_{i\pm\frac{1}{2},j\pm\frac{1}{2}})$.
Since solid boundaries move, it is necessary to update these values at
each time step.

\begin{figure}[!h]
\begin{center}\begin{tikzpicture}
	\draw[very thick, fill=gray] (0.05,0.05) circle (0.2) node[left] 
		{$\left.\begin{array}{r} \xeDD{-}{-}\\ \LvsDD{-}{-}>0\\ \IDD{-}{-}=1 \end{array}\right]$};
	\draw[very thick, fill=gray] (4.95,4.95) circle (0.2) node[right] 
		{$\left[\begin{array}{l} \xeDD{+}{+}\\ \LvsDD{+}{+}<0\\ \IDD{+}{+}=0 \end{array}\right.$};
	\draw[very thick, fill=gray] (4.95,0.05) circle (0.2) node[right] 
		{$\left[\begin{array}{l} \xeDD{+}{-}\\ \LvsDD{+}{-}>0\\ \IDD{+}{-}=1 \end{array}\right.$};
	\draw[very thick, fill=gray] (0.05,4.95) circle (0.2) node[left] 
		{$\left.\begin{array}{r} \xeDD{-}{+}\\ \LvsDD{-}{+}>0\\ \IDD{-}{+}=1 \end{array}\right]$};
	\draw[very thick, fill=gray] (2.75,4.75) rectangle (3.15,5.15) node[above left] %{$\xwjDD{+}$};
		{$\left.\begin{array}{r} \xwjDD{+}\\ \IDD{-}{+}\IDD{+}{+}=0\end{array}\right]\,$};	
	\draw[very thick, fill=gray] (4.75,2.05) rectangle (5.15,2.45) node[right] %{$\xwiDD{+}$};
		{$\left[\begin{array}{l} \xwiDD{+}\\ \IDD{-}{+}\IDD{+}{+}=0\end{array}\right.$};
	\draw (2.5,0) node[below]{$\IDD{-}{-}\IDD{-}{+}=1$};
	\draw (2.5,0.05) node[above]{$\LyDD{-}$};
	\draw (1.4,4.95) node[below]{$\LyDD{+}$};
	\draw (0.05,2.5) node[right]{$\LxDD{-}$};
	\draw (4.95,1.3) node[left]{$\LxDD{+}$};
	\draw (4,4) node{$\LwDD$};
	\draw[dashed] (0,0) rectangle (5,5);
	\draw[double] (3,6) arc (-10:90:-3.5);
	\draw (0.1,0.1) -- (0.1,4.9) -- (2.95,4.9) -- (4.9,2.25) -- (4.9,0.1) -- cycle ;
\end{tikzpicture}
\caption{Summary of the notations used in Appendix~\ref{LevelSet}: the
  cell $\Omega_{i,j}$ is shown with a dotted line, its corresponding
  virtual cell $\VirtDD$ with a solid line, and the solid boundary with a double
  line.}
\label{parametres}
\end{center}
\end{figure}
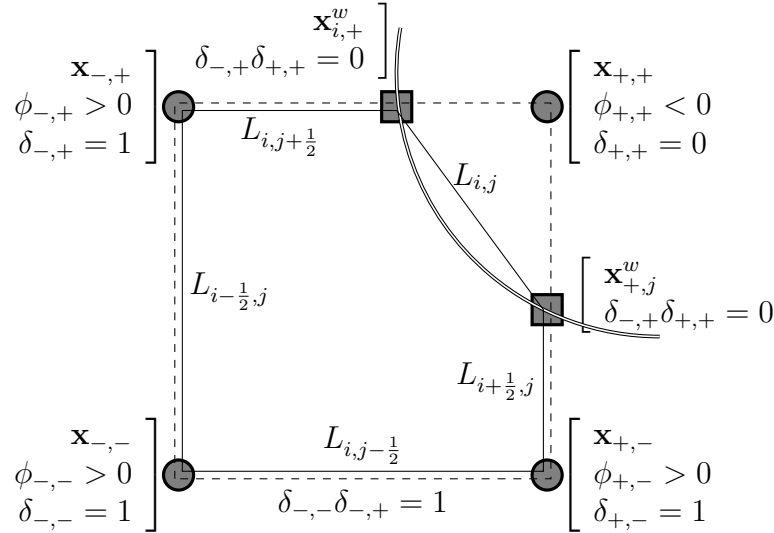

From now on, we consider a single cell~$\Omega_{i,j}$, and in order to
simplify the notations, indices $i$ and $j$ of variables defined at
the vertices of this cell will be omitted. For instance, the values
$\xDD_{i+\frac{1}{2},j+\frac{1}{2}}$ and
$\Lvset_{i+\frac{1}{2},j+\frac{1}{2}}$ will be denoted by
$\xeDD{+}{+}$ and $\phi_{+,+}$.  All the notations used here are shown
in figure~\ref{parametres}.

%=========
   \subsection{Identification of cell types}
%=========
It is clear that $\Omega_{i,j}$ is a gas cell if the values of
$\Lvset$ at its four vertices are positive. At the contrary, it is a
solid cell if these values are all negative. Finally, if these values
have different signs, the cell is cut by the solid boundary. These
three types are identified by looking at 
%   La fonction distance permet de déterminer le type d'une maille
%    $\Omega_{i,j}$ avec aisance.  En effet, si la valeur de cette
%    fonction aux quatre sommet de la maille sont positive, alors tous
%    ces points sont dans le gaz et la maille est gazeuse.
%    Réciproquement, la maille est solide quand la fonction distance est
%    négative sur chacun des sommets.  Enfin, lorsque la fonction
%    distance n'est pas du même signe sur l'ensemble des sommets de la
%    maille signifie que certain sommets sont dans le gaz et d'autres
%    dans le solide : la maille est donc traversée par la frontière
%    solide.  Numériquement, ces trois cas de figure sont déterminé de
%    la façon suivante :
% \begin{displaymath}\begin{aligned}
% 	\min > 0 &\Leftrightarrow &&\Omega_{i,j}\text{ est une maille de gaz}.\\
% 	\max < 0 &\Leftrightarrow &&\Omega_{i,j}\text{ est une maille solide}.\\
% 	\max > 0 \text{ et } \min < 0 
% 		&\Leftrightarrow &&\Omega_{i,j}\text{ est une maille coupée}.
% \end{aligned}\end{displaymath}
%où 
the sign of $m:=\min(\LvsDD{-}{-},
\LvsDD{+}{-},\LvsDD{+}{+},\LvsDD{-}{+})$ and $M:=
\max(\LvsDD{-}{-}, \LvsDD{+}{-},\LvsDD{+}{+},\LvsDD{-}{+})$: 
\begin{displaymath}\begin{aligned}
	m > 0 &\Leftrightarrow &&\Omega_{i,j}\text{ is a gas cell},\\
	M < 0 &\Leftrightarrow &&\Omega_{i,j}\text{ is a solid cell},\\
	M > 0 \text{ and } m < 0 
		&\Leftrightarrow &&\Omega_{i,j}\text{ is a cut cell}.
\end{aligned}\end{displaymath}
Of course, if a cell is crossed by a solid object thinner than $\Delta
x$, it will be considered as a gas cell: it is important that the
Cartesian grid is fine enough to resolve all the solid objects.

Finally, the case $\LvsDD{\pm}{\pm}=0$ is complex and is avoided by using the
modification $\Lvset:=\text{sign}(\Lvset)\times\max(\mid\Lvset\mid, 10^{-10}\Delta
x)$.  This means that grid points that are exactly on the solid
boundary are numerically considered as moved on a distance of
$10^{-10}\Delta x$. This has no influence on the accuracy of the
results.

%=========
   \subsection{Lengths of edges and normal vector of the virtual cell}
%=========

   To each of the four edges of $\Omega_{i,j}$ are associated the four
   points $\xwiDD{\pm}$ et $\xwjDD{\pm}$ defined as follows. If it is
   a cut cell, two of these points are intersection points of an edge
   with the solid boundary, and the two others are not used by the
   algorithm. When $\Omega_{i,j}$ is a gas or solid cell, none of
   these four points is used.  These points are defined by linear
   approximations
\begin{displaymath}
\xwiDD{\pm} = \left[\begin{aligned}
	&x_{i\pm\frac{1}{2}} \\
	&y_{j-\frac{1}{2}} - \Delta y 
		\frac{ \LvsDD{\pm}{-} }{\LvsDD{\pm}{+} - \LvsDD{\pm}{-}}
\end{aligned}\right]\quad\text{and}\quad
\xwjDD{\pm} = \left[\begin{aligned}
	&x_{i-\frac{1}{2}}- \Delta x
		\frac{ \LvsDD{-}{\pm} }{\LvsDD{+}{\pm} - \LvsDD{-}{\pm}}\\
	&y_{j\pm\frac{1}{2}}
\end{aligned}\right].
\end{displaymath}

To compute the lengths of the edges of $\Omega_{i,j}$, we use the
symbol $\IDD{\pm}{\pm}$ that is 1 if $\xeDD{\pm}{\pm}$ is inside the
gas, and 0 if it is inside the solid. This value is given by
$\IDD{\pm}{\pm} = \max(\LvsDD{\pm}{\pm},0)/|\LvsDD{\pm}{\pm}|$.
Also note that looking at the products $\IDD{\pm}{-}\IDD{\pm}{+}$ (or
$\IDD{-}{\pm}\IDD{+}{\pm}$) tells us if left and right edges
(respectively, upper and lower edges) are crossed by the solid
boundary. This notation allows us to easily write the four Cartesian
edge lengths of the virtual cell $\VirtDD$:
\begin{displaymath}\begin{aligned}
	&\DxDD{\pm}=\IDD{\pm}{-}\|\xwiDD{\pm}-\xeDD{\pm}{-}\|
		+\IDD{\pm}{+}\|\xeDD{\pm}{+}-\xwiDD{\pm}\|
		+\IDD{\pm}{-}\IDD{\pm}{+}\|\xeDD{\pm}{+}-\xeDD{\pm}{-}\|, \\
	&\DyDD{\pm}=\IDD{-}{\pm}\|\xwjDD{\pm}-\xeDD{-}{\pm}\|
		+\IDD{+}{\pm}\|\xeDD{+}{\pm}-\xwjDD{\pm}\|
		+\IDD{-}{\pm}\IDD{+}{\pm}\|\xeDD{+}{\pm}-\xeDD{-}{\pm}\|.
\end{aligned}\end{displaymath}
The length $\DwDD$ of the fifth edge (that fit with the solid boundary) is the norm
of the vector defined by the two solid boundary/edge intersection
points.

Finally, the normal vector to this edge is computed by a Green
formula, which naturally gives the correct outward direction:
\begin{displaymath}
	\nDD = \frac{1}{\DwDD}\left( \DxDD{+}-\DxDD{-}\right)\left[\begin{array}{c}1\\0\end{array}\right]
		+ \frac{1}{\DwDD}\left(\DyDD{+}-\DyDD{-}\right) \left[\begin{array}{c}0\\1\end{array}\right].
\end{displaymath}

%=========
   \subsection{Computation of $\sigma_{i,j}$}
   \label{subsec:sigmaij}
%=========
%

We consider one virtual cell $\VirtDD$. The set $\sigma_{i,j}$
collects the indices of all the virtual cells that merge with
$\VirtDD$. In this set, let us denote by $(i',j')$ the indices of the
unique virtual cell whose center $\xDD_{i',j'}$ is inside the gas,
that is to say
\begin{equation}\label{ii}
	\frac{1}{4}(\LvsDD{-}{-}+\LvsDD{+}{-}+\LvsDD{+}{+}+\LvsDD{-}{+})>0.
\end{equation}
We point out that the rule given in section~\ref{subsubsec:CV} ensures
that this "master" virtual cell is unique. 
Its indices $(i',j')$ are determined by the following algorithm
\begin{equation}\label{fusion}\begin{aligned}
&(i',j') := (i,j)\\
&\text{while (\ref{ii}) is false, do} \\
	&\qquad M :=\max(\DxDD{+},\DxDD{-},\DyDD{+},\DyDD{-})\\
	&\qquad \text{if }\DxDD{+}=M \text{ then } (i',j') := (i'+1,j')\\
	&\qquad \text{if }\DxDD{-}=M \text{ then } (i',j') := (i'-1,j')\\
	&\qquad \text{if }\DyDD{+}=M \text{ then } (i',j') := (i',j'+1)\\
	&\qquad \text{if }\DyDD{-}=M \text{ then } (i',j') := (i',j'-1)\\
&\text{end while}
\end{aligned}\end{equation}

Finally, note that the set $\sigma_{i,j}$ is not really computed:
practically, we compute the numerical fluxes across the edges of each
virtual cells, and these fluxes are directly added to the fluxes of
the master cell of indices $(i',j')$. \\ \\
%

%=========
   \subsection{Virtual cell and control volume areas}
%=========

The area of a virtual cell can be computed with the length of its
edges and the coordinates of its vertices by a Green formula: 
\begin{displaymath}\begin{aligned}
	\SVirtDD = \int_{\VirtDD} \dxDD = \frac{1}{2} \int_{\VirtDD}\nabla\cdot\xDD\dxDD 
	 &= \frac{1}{4}\DxDD{+}
		\left(\IDD{+}{+} \xeDD{+}{+} + \IDD{+}{-} \xeDD{+}{-}+\IwDD{+}{+}{+}{-} \xwiDD{+} \right) 
		\cdot \left[\begin{array}{c}1\\0\end{array}\right] \\
	&+\frac{1}{4} \DxDD{-}
		\left(\IDD{-}{+} \xeDD{-}{+}+\IDD{-}{-} \xeDD{-}{-}+\IwDD{-}{+}{-}{-} \xwiDD{-} \right) 
		\cdot \left[\begin{array}{c}1\\0\end{array}\right] \\
	&+\frac{1}{4} \DyDD{+} 
		\left(\IDD{+}{+} \xeDD{+}{+}+\IDD{-}{+} \xeDD{-}{+}+\IwDD{+}{+}{-}{+} \xwjDD{+} \right) 
		\cdot \left[\begin{array}{c}0\\1\end{array}\right] \\
	&+\frac{1}{4} \DyDD{-}
		\left(\IDD{-}{-} \xeDD{-}{-}+\IDD{+}{-} \xeDD{+}{-}+\IwDD{-}{-}{+}{-} \xwjDD{-} \right) 
		\cdot \left[\begin{array}{c}0\\1\end{array}\right] \\
	&+\frac{1}{4} \DwDD
		\left( \IwDD{+}{+}{+}{-} \xwiDD{+}+\IwDD{-}{+}{-}{-} \xwiDD{-}+\IwDD{+}{+}{-}{+} \xwjDD{+}+\IwDD{-}{-}{+}{-} \xwjDD{-} \right)
		\cdot\nDD 
\end{aligned}\end{displaymath}
Here, we used $\IwDD{\pm}{+}{\pm}{-}:=1-\IDD{\pm}{+}\IDD{\pm}{-}$
which is 1 if left or right edges are crossed by the solid boundary,
and 0 else. This formula is nothing but the sum on each edge of the
dot product between its normal vector and the vector pointing to the
center of the edge, multiplied by the length edge.

Finally, the area of the control volume is computed and stored in the
master cell of indices $(i',j')$ defined by algorithm~(\ref{fusion}),
with the following loop along all the cut cells: 
\begin{displaymath}\begin{aligned}
&\SCVolDD=0\text{ for every cells.}\\
&\text{For all }(i,j)\text{ do :}\\
&\qquad \text{compute }(i',j')\text{ with algorithm~(\ref{fusion})}.\\
&\qquad \SCVolDDp := \SCVolDDp + \SVirtDD\\
&\text{end do}
\end{aligned}\end{displaymath} 
%

%-----------------------------%
%                             %
%            IMAGES           %
%                             %
%-----------------------------%
\clearpage
\clearpage
\clearpage

\clearpage
\clearpage

\clearpage

\clearpage
\clearpage
\clearpage
\clearpage
\clearpage
\clearpage

\clearpage
\clearpage

\bibliographystyle{plain}
\bibliography{biblio}
\end{document}

%% file: Categories_1.pstex_t
\begin{picture}(0,0)%
\includegraphics{Categories_1.pdf}%
\end{picture}%
\setlength{\unitlength}{3947sp}%
\begingroup\makeatletter\ifx\SetFigFont\undefined%
\gdef\SetFigFont#1#2#3#4#5{%
  \reset@font\fontsize{#1}{#2pt}%
  \fontfamily{#3}\fontseries{#4}\fontshape{#5}%
  \selectfont}%
\fi\endgroup%
\begin{picture}(12024,9024)(589,-8773)
\put(3001,-6361){\makebox(0,0)[lb]{\smash{{\SetFigFont{29}{34.8}{\rmdefault}{\mddefault}{\updefault}{\color[rgb]{0,0,0}Gas}%
}}}}
\put(10201,-8161){\makebox(0,0)[lb]{\smash{{\SetFigFont{29}{34.8}{\rmdefault}{\mddefault}{\updefault}{\color[rgb]{0,0,0}Solid 1}%
}}}}
\put(5401,-3361){\makebox(0,0)[lb]{\smash{{\SetFigFont{29}{34.8}{\rmdefault}{\mddefault}{\updefault}{\color[rgb]{0,0,0}Solid 2}%
}}}}
\end{picture}%

%% file: Categories_2bis.pstex_t
\begin{picture}(0,0)%
\includegraphics{Categories_2bis.pdf}%
\end{picture}%
\setlength{\unitlength}{3947sp}%
\begingroup\makeatletter\ifx\SetFigFont\undefined%
\gdef\SetFigFont#1#2#3#4#5{%
  \reset@font\fontsize{#1}{#2pt}%
  \fontfamily{#3}\fontseries{#4}\fontshape{#5}%
  \selectfont}%
\fi\endgroup%
\begin{picture}(12024,9024)(589,-8773)
\end{picture}%

%% file: CutCells_1.pstex_t
\begin{picture}(0,0)%
\includegraphics{CutCells_1}%
\end{picture}%
\setlength{\unitlength}{4144sp}%
\begingroup\makeatletter\ifx\SetFigFont\undefined%
\gdef\SetFigFont#1#2#3#4#5{%
  \reset@font\fontsize{#1}{#2pt}%
  \fontfamily{#3}\fontseries{#4}\fontshape{#5}%
  \selectfont}%
\fi\endgroup%
\begin{picture}(5762,5290)(1192,-5822)
\put(6616,-3346){\rotatebox{270.0}{\makebox(0,0)[lb]{\smash{{\SetFigFont{25}{30.0}{\rmdefault}{\mddefault}{\updefault}{\color[rgb]{0,0,0}$\DxDD{+}$}%
}}}}}
\put(3736,-3316){\makebox(0,0)[lb]{\smash{{\SetFigFont{25}{30.0}{\rmdefault}{\mddefault}{\updefault}{\color[rgb]{0,0,0}$\VirtDD$}%
}}}}
\put(2926,-871){\makebox(0,0)[lb]{\smash{{\SetFigFont{25}{30.0}{\rmdefault}{\mddefault}{\updefault}{\color[rgb]{0,0,0}$\DyDD{+}$}%
}}}}
\put(5221,-1636){\rotatebox{306.0}{\makebox(0,0)[lb]{\smash{{\SetFigFont{25}{30.0}{\rmdefault}{\mddefault}{\updefault}{\color[rgb]{0,0,0}$\DwDD$}%
}}}}}
\put(2836,-5686){\makebox(0,0)[lb]{\smash{{\SetFigFont{25}{30.0}{\rmdefault}{\mddefault}{\updefault}{\color[rgb]{0,0,0}$\DyDD{-}=\Delta x$}%
}}}}
\put(1531,-4201){\rotatebox{90.0}{\makebox(0,0)[lb]{\smash{{\SetFigFont{25}{30.0}{\rmdefault}{\mddefault}{\updefault}{\color[rgb]{0,0,0}$\DxDD{-}=\Delta y$}%
}}}}}
\end{picture}%

%% file: CutCells_2.pstex_t
\begin{picture}(0,0)%
\includegraphics{CutCells_2}%
\end{picture}%
\setlength{\unitlength}{4144sp}%
\begingroup\makeatletter\ifx\SetFigFont\undefined%
\gdef\SetFigFont#1#2#3#4#5{%
  \reset@font\fontsize{#1}{#2pt}%
  \fontfamily{#3}\fontseries{#4}\fontshape{#5}%
  \selectfont}%
\fi\endgroup%
\begin{picture}(5447,5200)(1327,-5732)
\put(4051,-871){\makebox(0,0)[lb]{\smash{{\SetFigFont{25}{30.0}{\rmdefault}{\mddefault}{\updefault}{\color[rgb]{0,0,0}$\DyDD{+}$}%
}}}}
\put(3922,-2980){\rotatebox{308.0}{\makebox(0,0)[lb]{\smash{{\SetFigFont{25}{30.0}{\rmdefault}{\mddefault}{\updefault}{\color[rgb]{0,0,0}$\DwDD$}%
}}}}}
\put(4951,-2311){\makebox(0,0)[lb]{\smash{{\SetFigFont{25}{30.0}{\rmdefault}{\mddefault}{\updefault}{\color[rgb]{0,0,0}$\VirtDD$}%
}}}}
\put(6436,-1996){\rotatebox{270.0}{\makebox(0,0)[lb]{\smash{{\SetFigFont{25}{30.0}{\rmdefault}{\mddefault}{\updefault}{\color[rgb]{0,0,0}$\DxDD{+}$}%
}}}}}
\put(1666,-4201){\rotatebox{90.0}{\makebox(0,0)[lb]{\smash{{\SetFigFont{25}{30.0}{\rmdefault}{\mddefault}{\updefault}{\color[rgb]{1,1,1}$\DxDD{-}=0$}%
}}}}}
\put(3151,-5596){\makebox(0,0)[lb]{\smash{{\SetFigFont{25}{30.0}{\rmdefault}{\mddefault}{\updefault}{\color[rgb]{1,1,1}$\DyDD{-}=0$}%
}}}}
\end{picture}%

%% file: CutCells_3.pstex_t
\begin{picture}(0,0)%
\includegraphics{CutCells_3}%
\end{picture}%
\setlength{\unitlength}{4144sp}%
\begingroup\makeatletter\ifx\SetFigFont\undefined%
\gdef\SetFigFont#1#2#3#4#5{%
  \reset@font\fontsize{#1}{#2pt}%
  \fontfamily{#3}\fontseries{#4}\fontshape{#5}%
  \selectfont}%
\fi\endgroup%
\begin{picture}(5537,5245)(1237,-5732)
\put(2611,-4381){\makebox(0,0)[lb]{\smash{{\SetFigFont{25}{30.0}{\rmdefault}{\mddefault}{\updefault}{\color[rgb]{0,0,0}$\DyDD{-}$}%
}}}}
\put(3061,-3031){\makebox(0,0)[lb]{\smash{{\SetFigFont{25}{30.0}{\rmdefault}{\mddefault}{\updefault}{\color[rgb]{0,0,0}$\VirtDD$}%
}}}}
\put(5131,-3751){\rotatebox{76.0}{\makebox(0,0)[lb]{\smash{{\SetFigFont{25}{30.0}{\rmdefault}{\mddefault}{\updefault}{\color[rgb]{0,0,0}$\DwDD$}%
}}}}}
\put(2206,-5596){\makebox(0,0)[lb]{\smash{{\SetFigFont{25}{30.0}{\rmdefault}{\mddefault}{\updefault}{\color[rgb]{1,1,1}$\DyDD{-}$}%
}}}}
\put(3151,-826){\makebox(0,0)[lb]{\smash{{\SetFigFont{25}{30.0}{\rmdefault}{\mddefault}{\updefault}{\color[rgb]{0,0,0}$\DyDD{+}$}%
}}}}
\put(1576,-4336){\rotatebox{90.0}{\makebox(0,0)[lb]{\smash{{\SetFigFont{25}{30.0}{\rmdefault}{\mddefault}{\updefault}{\color[rgb]{0,0,0}$\DxDD{-}=\Delta y$}%
}}}}}
\put(6301,-3166){\makebox(0,0)[lb]{\smash{{\SetFigFont{25}{30.0}{\rmdefault}{\mddefault}{\updefault}{\color[rgb]{1,1,1}$\DyDD{+}$}%
}}}}
\put(6436,-2266){\rotatebox{270.0}{\makebox(0,0)[lb]{\smash{{\SetFigFont{25}{30.0}{\rmdefault}{\mddefault}{\updefault}{\color[rgb]{1,1,1}$\DxDD{+}=0$}%
}}}}}
\end{picture}%

%% file: test_ccn_cvn.pstex_t
\begin{picture}(0,0)%
\includegraphics{test_ccn_cvn}%
\end{picture}%
\setlength{\unitlength}{4144sp}%
\begingroup\makeatletter\ifx\SetFigFont\undefined%
\gdef\SetFigFont#1#2#3#4#5{%
  \reset@font\fontsize{#1}{#2pt}%
  \fontfamily{#3}\fontseries{#4}\fontshape{#5}%
  \selectfont}%
\fi\endgroup%
\begin{picture}(13964,4558)(-1769,-4157)
\put(1621,-1456){\makebox(0,0)[lb]{\smash{{\SetFigFont{25}{30.0}{\rmdefault}{\mddefault}{\updefault}{\color[rgb]{0,0,0}$\VirtDD^n$}%
}}}}
\put(8011,-1321){\makebox(0,0)[lb]{\smash{{\SetFigFont{25}{30.0}{\rmdefault}{\mddefault}{\updefault}{\color[rgb]{0,0,0}$\CVolDD^n(t_n)$}%
}}}}
\put(3781,-2851){\makebox(0,0)[lb]{\smash{{\SetFigFont{25}{30.0}{\rmdefault}{\mddefault}{\updefault}{\color[rgb]{0,0,0}$\xDD_{i,j}$}%
}}}}
\put(361, 29){\makebox(0,0)[lb]{\smash{{\SetFigFont{25}{30.0}{\rmdefault}{\mddefault}{\updefault}{\color[rgb]{0,0,0}$\Omega_g$}%
}}}}
\put(4591,-4021){\makebox(0,0)[lb]{\smash{{\SetFigFont{25}{30.0}{\rmdefault}{\mddefault}{\updefault}{\color[rgb]{0,0,0}$\Omega_s$}%
}}}}
\put(-1754,-1996){\makebox(0,0)[lb]{\smash{{\SetFigFont{25}{30.0}{\rmdefault}{\mddefault}{\updefault}{\color[rgb]{0,0,0}Time $t^n$}%
}}}}
\put(10531,-2851){\makebox(0,0)[lb]{\smash{{\SetFigFont{25}{30.0}{\rmdefault}{\mddefault}{\updefault}{\color[rgb]{0,0,0}$\xDD_{i,j}$}%
}}}}
\put(7111, 29){\makebox(0,0)[lb]{\smash{{\SetFigFont{25}{30.0}{\rmdefault}{\mddefault}{\updefault}{\color[rgb]{0,0,0}$\Omega_g$}%
}}}}
\put(11341,-4021){\makebox(0,0)[lb]{\smash{{\SetFigFont{25}{30.0}{\rmdefault}{\mddefault}{\updefault}{\color[rgb]{0,0,0}$\Omega_s$}%
}}}}
\put(10531,-1051){\makebox(0,0)[lb]{\smash{{\SetFigFont{25}{30.0}{\rmdefault}{\mddefault}{\updefault}{\color[rgb]{0,0,0}$\xDD_{i,j+1}$}%
}}}}
\end{picture}%

%% file: test_ccnpun_cvnpun.pstex_t
\begin{picture}(0,0)%
\includegraphics{test_ccnpun_cvnpun}%
\end{picture}%
\setlength{\unitlength}{4144sp}%
\begingroup\makeatletter\ifx\SetFigFont\undefined%
\gdef\SetFigFont#1#2#3#4#5{%
  \reset@font\fontsize{#1}{#2pt}%
  \fontfamily{#3}\fontseries{#4}\fontshape{#5}%
  \selectfont}%
\fi\endgroup%
\begin{picture}(14009,4736)(-1769,-4335)
\put(4681,-4021){\makebox(0,0)[lb]{\smash{{\SetFigFont{25}{30.0}{\rmdefault}{\mddefault}{\updefault}{\color[rgb]{0,0,0}$\Omega_s$}%
}}}}
\put(451, 29){\makebox(0,0)[lb]{\smash{{\SetFigFont{25}{30.0}{\rmdefault}{\mddefault}{\updefault}{\color[rgb]{0,0,0}$\Omega_g$}%
}}}}
\put(3871,-2851){\makebox(0,0)[lb]{\smash{{\SetFigFont{25}{30.0}{\rmdefault}{\mddefault}{\updefault}{\color[rgb]{0,0,0}$\xDD_{i,j}$}%
}}}}
\put(11386,-4021){\makebox(0,0)[lb]{\smash{{\SetFigFont{25}{30.0}{\rmdefault}{\mddefault}{\updefault}{\color[rgb]{0,0,0}$\Omega_s$}%
}}}}
\put(7156, 29){\makebox(0,0)[lb]{\smash{{\SetFigFont{25}{30.0}{\rmdefault}{\mddefault}{\updefault}{\color[rgb]{0,0,0}$\Omega_g$}%
}}}}
\put(10576,-2851){\makebox(0,0)[lb]{\smash{{\SetFigFont{25}{30.0}{\rmdefault}{\mddefault}{\updefault}{\color[rgb]{0,0,0}$\xDD_{i,j}$}%
}}}}
\put(-1754,-1996){\makebox(0,0)[lb]{\smash{{\SetFigFont{25}{30.0}{\rmdefault}{\mddefault}{\updefault}{\color[rgb]{0,0,0}Time $t^{n+1}$}%
}}}}
\put(7921,-1141){\makebox(0,0)[lb]{\smash{{\SetFigFont{25}{30.0}{\rmdefault}{\mddefault}{\updefault}{\color[rgb]{0,0,0}$\CVolDD^n(t_{n+1})$}%
}}}}
\put(1801,-1546){\makebox(0,0)[lb]{\smash{{\SetFigFont{25}{30.0}{\rmdefault}{\mddefault}{\updefault}{\color[rgb]{0,0,0}$\VirtDD^{n+1}$}%
}}}}
\end{picture}%

%% file: test_cvpunnpun.pstex_t
\begin{picture}(0,0)%
\includegraphics{test_cvpunnpun}%
\end{picture}%
\setlength{\unitlength}{4144sp}%
\begingroup\makeatletter\ifx\SetFigFont\undefined%
\gdef\SetFigFont#1#2#3#4#5{%
  \reset@font\fontsize{#1}{#2pt}%
  \fontfamily{#3}\fontseries{#4}\fontshape{#5}%
  \selectfont}%
\fi\endgroup%
\begin{picture}(6989,4736)(5656,-4335)
\put(11791,-4021){\makebox(0,0)[lb]{\smash{{\SetFigFont{25}{30.0}{\rmdefault}{\mddefault}{\updefault}{\color[rgb]{0,0,0}$\Omega_s$}%
}}}}
\put(7561, 29){\makebox(0,0)[lb]{\smash{{\SetFigFont{25}{30.0}{\rmdefault}{\mddefault}{\updefault}{\color[rgb]{0,0,0}$\Omega_g$}%
}}}}
\put(10981,-2851){\makebox(0,0)[lb]{\smash{{\SetFigFont{25}{30.0}{\rmdefault}{\mddefault}{\updefault}{\color[rgb]{0,0,0}$\xDD_{i,j}$}%
}}}}
\put(8461,-3256){\makebox(0,0)[lb]{\smash{{\SetFigFont{25}{30.0}{\rmdefault}{\mddefault}{\updefault}{\color[rgb]{0,0,0}$\xDD_{i-1,j}$}%
}}}}
\put(5671,-1456){\makebox(0,0)[lb]{\smash{{\SetFigFont{25}{30.0}{\rmdefault}{\mddefault}{\updefault}{\color[rgb]{0,0,0}$\CVolDD^{n+1}(t_{n+1})$}%
}}}}
\end{picture}%

%% file: steps.pstex_t
\begin{picture}(0,0)%
\includegraphics{steps}%
\end{picture}%
\setlength{\unitlength}{4144sp}%
\begingroup\makeatletter\ifx\SetFigFont\undefined%
\gdef\SetFigFont#1#2#3#4#5{%
  \reset@font\fontsize{#1}{#2pt}%
  \fontfamily{#3}\fontseries{#4}\fontshape{#5}%
  \selectfont}%
\fi\endgroup%
\begin{picture}(13995,3216)(418,-2344)
\put(676,614){\makebox(0,0)[lb]{\smash{{\SetFigFont{20}{24.0}{\rmdefault}{\mddefault}{\updefault}{\color[rgb]{0,0,0}Iteration $n$}%
}}}}
\put(11926,614){\makebox(0,0)[lb]{\smash{{\SetFigFont{20}{24.0}{\rmdefault}{\mddefault}{\updefault}{\color[rgb]{0,0,0}Iteration $n+1$}%
}}}}
\put(8911,-61){\makebox(0,0)[lb]{\smash{{\SetFigFont{20}{24.0}{\rmdefault}{\mddefault}{\updefault}{\color[rgb]{0,0,0}Boundary at $t^{n+1}$}%
}}}}
\put(2251,-61){\makebox(0,0)[lb]{\smash{{\SetFigFont{20}{24.0}{\rmdefault}{\mddefault}{\updefault}{\color[rgb]{0,0,0}Boundary at $t^n$}%
}}}}
\put(3241,-961){\makebox(0,0)[lb]{\smash{{\SetFigFont{20}{24.0}{\rmdefault}{\mddefault}{\updefault}{\color[rgb]{0,0,0}Merging}%
}}}}
\put(6931,-961){\makebox(0,0)[lb]{\smash{{\SetFigFont{20}{24.0}{\rmdefault}{\mddefault}{\updefault}{\color[rgb]{0,0,0}Scheme}%
}}}}
\put(10531,-961){\makebox(0,0)[lb]{\smash{{\SetFigFont{20}{24.0}{\rmdefault}{\mddefault}{\updefault}{\color[rgb]{0,0,0}Update}%
}}}}
\end{picture}%

%% file: appearing.pstex_t
\begin{picture}(0,0)%
\includegraphics{appearing}%
\end{picture}%
\setlength{\unitlength}{4144sp}%
\begingroup\makeatletter\ifx\SetFigFont\undefined%
\gdef\SetFigFont#1#2#3#4#5{%
  \reset@font\fontsize{#1}{#2pt}%
  \fontfamily{#3}\fontseries{#4}\fontshape{#5}%
  \selectfont}%
\fi\endgroup%
\begin{picture}(13995,3216)(418,-2344)
\put(676,614){\makebox(0,0)[lb]{\smash{{\SetFigFont{20}{24.0}{\rmdefault}{\mddefault}{\updefault}{\color[rgb]{0,0,0}Iteration $n$}%
}}}}
\put(11926,614){\makebox(0,0)[lb]{\smash{{\SetFigFont{20}{24.0}{\rmdefault}{\mddefault}{\updefault}{\color[rgb]{0,0,0}Iteration $n+1$}%
}}}}
\put(3241,-961){\makebox(0,0)[lb]{\smash{{\SetFigFont{20}{24.0}{\rmdefault}{\mddefault}{\updefault}{\color[rgb]{0,0,0}Merging}%
}}}}
\put(6931,-961){\makebox(0,0)[lb]{\smash{{\SetFigFont{20}{24.0}{\rmdefault}{\mddefault}{\updefault}{\color[rgb]{0,0,0}Scheme}%
}}}}
\put(10531,-961){\makebox(0,0)[lb]{\smash{{\SetFigFont{20}{24.0}{\rmdefault}{\mddefault}{\updefault}{\color[rgb]{0,0,0}Update}%
}}}}
\put(9181,-61){\makebox(0,0)[lb]{\smash{{\SetFigFont{20}{24.0}{\rmdefault}{\mddefault}{\updefault}{\color[rgb]{0,0,0}Boundary at $t^{n+1}$}%
}}}}
\put(2791,-61){\makebox(0,0)[lb]{\smash{{\SetFigFont{20}{24.0}{\rmdefault}{\mddefault}{\updefault}{\color[rgb]{0,0,0}Boundary at $t^n$}%
}}}}
\end{picture}%

%% file: disappearing.pstex_t
\begin{picture}(0,0)%
\includegraphics{disappearing}%
\end{picture}%
\setlength{\unitlength}{4144sp}%
\begingroup\makeatletter\ifx\SetFigFont\undefined%
\gdef\SetFigFont#1#2#3#4#5{%
  \reset@font\fontsize{#1}{#2pt}%
  \fontfamily{#3}\fontseries{#4}\fontshape{#5}%
  \selectfont}%
\fi\endgroup%
\begin{picture}(13995,3216)(418,-2344)
\put(676,614){\makebox(0,0)[lb]{\smash{{\SetFigFont{20}{24.0}{\rmdefault}{\mddefault}{\updefault}{\color[rgb]{0,0,0}Iteration $n$}%
}}}}
\put(11926,614){\makebox(0,0)[lb]{\smash{{\SetFigFont{20}{24.0}{\rmdefault}{\mddefault}{\updefault}{\color[rgb]{0,0,0}Iteration $n+1$}%
}}}}
\put(3241,-961){\makebox(0,0)[lb]{\smash{{\SetFigFont{20}{24.0}{\rmdefault}{\mddefault}{\updefault}{\color[rgb]{0,0,0}Merging}%
}}}}
\put(6931,-961){\makebox(0,0)[lb]{\smash{{\SetFigFont{20}{24.0}{\rmdefault}{\mddefault}{\updefault}{\color[rgb]{0,0,0}Scheme}%
}}}}
\put(10531,-961){\makebox(0,0)[lb]{\smash{{\SetFigFont{20}{24.0}{\rmdefault}{\mddefault}{\updefault}{\color[rgb]{0,0,0}Update}%
}}}}
\put(2071,-61){\makebox(0,0)[lb]{\smash{{\SetFigFont{20}{24.0}{\rmdefault}{\mddefault}{\updefault}{\color[rgb]{0,0,0}Boundary at $t^n$}%
}}}}
\put(9901,-61){\makebox(0,0)[lb]{\smash{{\SetFigFont{20}{24.0}{\rmdefault}{\mddefault}{\updefault}{\color[rgb]{0,0,0}Boundary at $t^{n+1}$}%
}}}}
\end{picture}%

%% file: Taguchi_1.pstex_t
\begin{picture}(0,0)%
\includegraphics{Taguchi_1}%
\end{picture}%
\setlength{\unitlength}{4144sp}%
\begingroup\makeatletter\ifx\SetFigFont\undefined%
\gdef\SetFigFont#1#2#3#4#5{%
  \reset@font\fontsize{#1}{#2pt}%
  \fontfamily{#3}\fontseries{#4}\fontshape{#5}%
  \selectfont}%
\fi\endgroup%
\begin{picture}(4837,2386)(2869,-3074)
\put(3646,-871){\makebox(0,0)[lb]{\smash{{\SetFigFont{14}{16.8}{\rmdefault}{\mddefault}{\updefault}{\color[rgb]{0,0,0}$T_0$}%
}}}}
\put(6481,-1411){\makebox(0,0)[lb]{\smash{{\SetFigFont{14}{16.8}{\rmdefault}{\mddefault}{\updefault}{\color[rgb]{0,0,0}$D/10$}%
}}}}
\put(3636,-2996){\makebox(0,0)[lb]{\smash{{\SetFigFont{14}{16.8}{\rmdefault}{\mddefault}{\updefault}{\color[rgb]{0,0,0}$T_0$}%
}}}}
\put(5681,-2556){\makebox(0,0)[lb]{\smash{{\SetFigFont{14}{16.8}{\rmdefault}{\mddefault}{\updefault}{\color[rgb]{0,0,0}$4D$}%
}}}}
\put(7691,-1936){\makebox(0,0)[lb]{\smash{{\SetFigFont{14}{16.8}{\rmdefault}{\mddefault}{\updefault}{\color[rgb]{0,0,0}$4D$}%
}}}}
\put(7246,-1936){\makebox(0,0)[lb]{\smash{{\SetFigFont{14}{16.8}{\rmdefault}{\mddefault}{\updefault}{\color[rgb]{0,0,0}$D$}%
}}}}
\put(4996,-1591){\makebox(0,0)[lb]{\smash{{\SetFigFont{14}{16.8}{\rmdefault}{\mddefault}{\updefault}{\color[rgb]{0,0,0}$2T_0$}%
}}}}
\put(4591,-1591){\makebox(0,0)[lb]{\smash{{\SetFigFont{14}{16.8}{\rmdefault}{\mddefault}{\updefault}{\color[rgb]{0,0,0}$T_0$}%
}}}}
\end{picture}%

%% file: Radiometer_1.pstex_t
\begin{picture}(0,0)%
\includegraphics{Radiometer_1}%
\end{picture}%
\setlength{\unitlength}{4144sp}%
\begingroup\makeatletter\ifx\SetFigFont\undefined%
\gdef\SetFigFont#1#2#3#4#5{%
  \reset@font\fontsize{#1}{#2pt}%
  \fontfamily{#3}\fontseries{#4}\fontshape{#5}%
  \selectfont}%
\fi\endgroup%
\begin{picture}(4524,4524)(4489,-5923)
\put(8191,-1906){\makebox(0,0)[lb]{\smash{{\SetFigFont{14}{16.8}{\rmdefault}{\mddefault}{\updefault}{\color[rgb]{0,0,0}$T_0$}%
}}}}
\put(6871,-4771){\makebox(0,0)[lb]{\smash{{\SetFigFont{14}{16.8}{\rmdefault}{\mddefault}{\updefault}{\color[rgb]{0,0,0}$T_c$}%
}}}}
\put(7661,-4336){\makebox(0,0)[lb]{\smash{{\SetFigFont{14}{16.8}{\rmdefault}{\mddefault}{\updefault}{\color[rgb]{0,0,0}$R$}%
}}}}
\put(8216,-3706){\makebox(0,0)[lb]{\smash{{\SetFigFont{14}{16.8}{\rmdefault}{\mddefault}{\updefault}{\color[rgb]{0,0,0}$l$}%
}}}}
\put(5131,-3191){\makebox(0,0)[lb]{\smash{{\SetFigFont{14}{16.8}{\rmdefault}{\mddefault}{\updefault}{\color[rgb]{0,0,0}$L$}%
}}}}
\put(6431,-4766){\makebox(0,0)[lb]{\smash{{\SetFigFont{14}{16.8}{\rmdefault}{\mddefault}{\updefault}{\color[rgb]{0,0,0}$T_h$}%
}}}}
\end{picture}%

%% file: pump_1.pstex_t
\begin{picture}(0,0)%
\includegraphics{pump_1}%
\end{picture}%
\setlength{\unitlength}{4144sp}%
\begingroup\makeatletter\ifx\SetFigFont\undefined%
\gdef\SetFigFont#1#2#3#4#5{%
  \reset@font\fontsize{#1}{#2pt}%
  \fontfamily{#3}\fontseries{#4}\fontshape{#5}%
  \selectfont}%
\fi\endgroup%
\begin{picture}(4977,4977)(2013,-6149)
\put(6286,-2001){\makebox(0,0)[lb]{\smash{{\SetFigFont{14}{16.8}{\rmdefault}{\mddefault}{\updefault}{\color[rgb]{0,0,0}$r$}%
}}}}
\put(3756,-2656){\makebox(0,0)[lb]{\smash{{\SetFigFont{14}{16.8}{\rmdefault}{\mddefault}{\updefault}{\color[rgb]{0,0,0}$r$}%
}}}}
\put(5266,-3591){\makebox(0,0)[lb]{\smash{{\SetFigFont{14}{16.8}{\rmdefault}{\mddefault}{\updefault}{\color[rgb]{0,0,0}$4r$}%
}}}}
\end{picture}%

%% file: pump_2.pstex_t
\begin{picture}(0,0)%
\includegraphics{pump_2}%
\end{picture}%
\setlength{\unitlength}{4144sp}%
\begingroup\makeatletter\ifx\SetFigFont\undefined%
\gdef\SetFigFont#1#2#3#4#5{%
  \reset@font\fontsize{#1}{#2pt}%
  \fontfamily{#3}\fontseries{#4}\fontshape{#5}%
  \selectfont}%
\fi\endgroup%
\begin{picture}(10017,9438)(1246,-8587)
\put(1276,-5551){\makebox(0,0)[lb]{\smash{{\SetFigFont{20}{24.0}{\rmdefault}{\mddefault}{\updefault}{\color[rgb]{0,0,0}$-0.16$}%
}}}}
\put(1261,-8191){\makebox(0,0)[lb]{\smash{{\SetFigFont{20}{24.0}{\rmdefault}{\mddefault}{\updefault}{\color[rgb]{0,0,0}$-0.4$}%
}}}}
\put(3601,-8491){\makebox(0,0)[lb]{\smash{{\SetFigFont{20}{24.0}{\rmdefault}{\mddefault}{\updefault}{\color[rgb]{0,0,0}$-0.24$}%
}}}}
\put(1426,-1936){\makebox(0,0)[lb]{\smash{{\SetFigFont{20}{24.0}{\rmdefault}{\mddefault}{\updefault}{\color[rgb]{0,0,0}$0.16$}%
}}}}
\put(1996,-8476){\makebox(0,0)[lb]{\smash{{\SetFigFont{20}{24.0}{\rmdefault}{\mddefault}{\updefault}{\color[rgb]{0,0,0}$-0.4$}%
}}}}
\put(6661,-8491){\makebox(0,0)[lb]{\smash{{\SetFigFont{20}{24.0}{\rmdefault}{\mddefault}{\updefault}{\color[rgb]{0,0,0}$0$}%
}}}}
\put(9151,-8491){\makebox(0,0)[lb]{\smash{{\SetFigFont{20}{24.0}{\rmdefault}{\mddefault}{\updefault}{\color[rgb]{0,0,0}$0.24$}%
}}}}
\put(10906,-8491){\makebox(0,0)[lb]{\smash{{\SetFigFont{20}{24.0}{\rmdefault}{\mddefault}{\updefault}{\color[rgb]{0,0,0}$0.4$}%
}}}}
\put(1501,599){\makebox(0,0)[lb]{\smash{{\SetFigFont{20}{24.0}{\rmdefault}{\mddefault}{\updefault}{\color[rgb]{0,0,0}$0.4$}%
}}}}
\put(2701,-3661){\makebox(0,0)[lb]{\smash{{\SetFigFont{20}{24.0}{\rmdefault}{\mddefault}{\updefault}{\color[rgb]{0,0,0}Inflow}%
}}}}
\put(9901,-3661){\makebox(0,0)[lb]{\smash{{\SetFigFont{20}{24.0}{\rmdefault}{\mddefault}{\updefault}{\color[rgb]{0,0,0}Outflow}%
}}}}
\end{picture}%

%% file: radiometre3D-2.pstex_t
\begin{picture}(0,0)%
\includegraphics{radiometre3D-2}%
\end{picture}%
\setlength{\unitlength}{4144sp}%
\begingroup\makeatletter\ifx\SetFigFont\undefined%
\gdef\SetFigFont#1#2#3#4#5{%
  \reset@font\fontsize{#1}{#2pt}%
  \fontfamily{#3}\fontseries{#4}\fontshape{#5}%
  \selectfont}%
\fi\endgroup%
\begin{picture}(12356,5856)(-1093,-6190)
\put(8219,-6168){\rotatebox{45.0}{\makebox(0,0)[lb]{\smash{{\SetFigFont{20}{24.0}{\rmdefault}{\mddefault}{\updefault}{\color[rgb]{0,0,0}0.05}%
}}}}}
\put(8553,-6175){\rotatebox{45.0}{\makebox(0,0)[lb]{\smash{{\SetFigFont{20}{24.0}{\rmdefault}{\mddefault}{\updefault}{\color[rgb]{0,0,0}0.25}%
}}}}}
\put(9044,-6173){\rotatebox{45.0}{\makebox(0,0)[lb]{\smash{{\SetFigFont{20}{24.0}{\rmdefault}{\mddefault}{\updefault}{\color[rgb]{0,0,0}0.75}%
}}}}}
\put(9666,-6171){\rotatebox{45.0}{\makebox(0,0)[lb]{\smash{{\SetFigFont{20}{24.0}{\rmdefault}{\mddefault}{\updefault}{\color[rgb]{0,0,0}1.25}%
}}}}}
\put(10452,-6136){\rotatebox{45.0}{\makebox(0,0)[lb]{\smash{{\SetFigFont{20}{24.0}{\rmdefault}{\mddefault}{\updefault}{\color[rgb]{0,0,0}2.00}%
}}}}}
\put(5471,-2785){\makebox(0,0)[lb]{\smash{{\SetFigFont{20}{24.0}{\rmdefault}{\mddefault}{\updefault}{\color[rgb]{0,0,0}0.50}%
}}}}
\put(5471,-1051){\makebox(0,0)[lb]{\smash{{\SetFigFont{20}{24.0}{\rmdefault}{\mddefault}{\updefault}{\color[rgb]{0,0,0}2.00}%
}}}}
\put(2028,-6175){\rotatebox{45.0}{\makebox(0,0)[lb]{\smash{{\SetFigFont{20}{24.0}{\rmdefault}{\mddefault}{\updefault}{\color[rgb]{0,0,0}0.25}%
}}}}}
\put(3141,-6171){\rotatebox{45.0}{\makebox(0,0)[lb]{\smash{{\SetFigFont{20}{24.0}{\rmdefault}{\mddefault}{\updefault}{\color[rgb]{0,0,0}1.25}%
}}}}}
\put(3927,-6136){\rotatebox{45.0}{\makebox(0,0)[lb]{\smash{{\SetFigFont{20}{24.0}{\rmdefault}{\mddefault}{\updefault}{\color[rgb]{0,0,0}2.00}%
}}}}}
\put(-1074,-2995){\makebox(0,0)[lb]{\smash{{\SetFigFont{20}{24.0}{\rmdefault}{\mddefault}{\updefault}{\color[rgb]{0,0,0}0.25}%
}}}}
\put(-1074,-1915){\makebox(0,0)[lb]{\smash{{\SetFigFont{20}{24.0}{\rmdefault}{\mddefault}{\updefault}{\color[rgb]{0,0,0}1.25}%
}}}}
\put(-1078,-3245){\makebox(0,0)[lb]{\smash{{\SetFigFont{20}{24.0}{\rmdefault}{\mddefault}{\updefault}{\color[rgb]{0,0,0}0.05}%
}}}}
\put(-1074,-1063){\makebox(0,0)[lb]{\smash{{\SetFigFont{20}{24.0}{\rmdefault}{\mddefault}{\updefault}{\color[rgb]{0,0,0}2.00}%
}}}}
\put(1688,-6168){\rotatebox{45.0}{\makebox(0,0)[lb]{\smash{{\SetFigFont{20}{24.0}{\rmdefault}{\mddefault}{\updefault}{\color[rgb]{0,0,0}0.05}%
}}}}}
\put(-134,-601){\makebox(0,0)[lb]{\smash{{\SetFigFont{20}{24.0}{\rmdefault}{\mddefault}{\updefault}{\color[rgb]{0,0,0}$\vec{y}/L$}%
}}}}
\put(4366,-5281){\makebox(0,0)[lb]{\smash{{\SetFigFont{20}{24.0}{\rmdefault}{\mddefault}{\updefault}{\color[rgb]{0,0,0}$\vec{x}/L$}%
}}}}
\put(10801,-5281){\makebox(0,0)[lb]{\smash{{\SetFigFont{20}{24.0}{\rmdefault}{\mddefault}{\updefault}{\color[rgb]{0,0,0}$\vec{y}/L$}%
}}}}
\put(6391,-646){\makebox(0,0)[lb]{\smash{{\SetFigFont{20}{24.0}{\rmdefault}{\mddefault}{\updefault}{\color[rgb]{0,0,0}$\vec{z}/L$}%
}}}}
\put(10801,-4921){\makebox(0,0)[lb]{\smash{{\SetFigFont{20}{24.0}{\rmdefault}{\mddefault}{\updefault}{\color[rgb]{0,0,0}$\vec{x}/L$}%
}}}}
\end{picture}%